\documentclass[nonblindrev]{informs3noheader}

\OneAndAHalfSpacedXI

\usepackage[margin=1in]{geometry}

\usepackage{endnotes}
\let\footnote=\endnote
\let\enotesize=\normalsize
\def\notesname{Endnotes}%
\def\makeenmark{\hbox to1.275em{\theenmark.\enskip\hss}}
\def\enoteformat{
    \rightskip0pt\leftskip0pt\parindent=1.275em
    \leavevmode\llap{\makeenmark}
}
\usepackage[bb=dsserif]{mathalpha}
\usepackage{bm}
\usepackage{bbm}
\usepackage{booktabs}
\usepackage{colortbl}
\usepackage{multirow}
\usepackage{multicol}
\usepackage[tight-spacing=true]{siunitx}
\usepackage{braket}
\usepackage{algorithm}
\usepackage{algorithmicx}
\usepackage{algpseudocode}
\usepackage{mathtools}
\usepackage{caption}
\usepackage{subcaption}
\usepackage{lscape}
\usepackage[shortlabels,inline]{enumitem}
\usepackage[flushleft]{threeparttable}
\usepackage{longtable}

\usepackage{thmtools,thm-restate} 

\usepackage{pifont}

\usepackage{tikz}
\usetikzlibrary{arrows.meta}
\usetikzlibrary{shapes.geometric,fit}
\usetikzlibrary{positioning,calc,shapes.misc}
\usetikzlibrary{decorations.pathreplacing}
\usetikzlibrary{decorations.text}
\usepackage{pgfplots}
\pgfplotsset{width=10cm,compat=1.18}
\usepgfplotslibrary{external}
\usepgfplotslibrary{fillbetween}

\definecolor{darkred}{rgb}{0.5,0,0}
\definecolor{darkblue}{rgb}{0,0.2,0.5}
\definecolor{darkpurple}{rgb}{0.3,0.1,0.3}
\definecolor{darkgreen}{rgb}{0.1,0.4,0.2}

\usepackage[hypertexnames=false]{hyperref}
\hypersetup{
    colorlinks=true,
    linkcolor=blue,
    citecolor=blue,
    filecolor=black,      
    urlcolor=black,
}

\def\R{\ensuremath\mathbb{R}}

\def\calA{\ensuremath\mathcal{A}}

\def\calN{\ensuremath\mathcal{N}}

\def\calP{\ensuremath\mathcal{P}}
\def\calQ{\ensuremath\mathcal{Q}}
\def\calR{\ensuremath\mathcal{R}}
\def\calS{\ensuremath\mathcal{S}}

\def\calV{\ensuremath\mathcal{V}}

\def\bi{\ensuremath\boldsymbol{i}}

\def\bw{\ensuremath\boldsymbol{w}}

\def\bz{\ensuremath\boldsymbol{z}}

\def\bnu{\ensuremath\boldsymbol{\nu}}

\def\bmu{\ensuremath\boldsymbol{\mu}}

\def\bkappa{\ensuremath\boldsymbol{\kappa}}

\def\st{\text{s.t.}}

\definecolor{mygreen}    {RGB}{0,90,0}
\definecolor{myblue}     {RGB}{0,51,140}
\definecolor{myorange}   {RGB}{238,118,0}
\definecolor{myred}      {RGB}{126,0,0}
\definecolor{mygray}     {RGB}{100,100,105}
\definecolor{mygrayblue} {RGB}{0,128,128}
\definecolor{mygraygreen}{RGB}{128,128,0}
\definecolor{DarkPurple}     {RGB}{142, 36, 170}
\definecolor{LightPurple}    {RGB}{57, 130, 7}

\newcommand{\Depots}[0]{\ensuremath \mathcal{V}_{\text{D}}}
\newcommand{\Custs}[0]{\ensuremath \mathcal{V}_{\text{C}}}
\newcommand{\Stations}[0]{\ensuremath \mathcal{V}_{\text{S}}}

\def\OPT{\texttt{OPT}}
\def\EVRP{\texttt{EVRP}}
\def\OPTLP{\overline{\texttt{OPT}}}
\def\EVRPLP{\overline{\texttt{EVRP}}}

\newcommand{\snode}[2][]{
    \ensuremath{n_{\text{start}}^{#1}\left({#2}\right)}
}
\newcommand{\enode}[2][]{
    \ensuremath{n_{\text{end}}^{#1}\left({#2}\right)}
}
\newcommand{\nnode}{
    \ensuremath{n_{\text{next}}}
}
\newcommand{\serve}[2]{
    \ensuremath{\gamma_{#1}^{#2}}
}

\newcommand{\cbar}[1]{\ensuremath \overline{c} \left( #1 \right)}

\newcommand{\Ttravel}[1]{\ensuremath T^{\text{travel}} \left( #1 \right)}
\newcommand{\TendMin}[1]{\ensuremath T^{\text{end}}_{\text{min}} \left( #1 \right)}
\newcommand{\TstartMax}[1]{\ensuremath T^{\text{start}}_{\text{max}} \left( #1 \right)}
\newcommand{\Teff}[1]{\ensuremath T^{\text{eff}} \left( #1 \right)}
\newcommand{\BendMin}[1]{\ensuremath B^{\text{end}}_{\text{min}} \left( #1 \right)}
\newcommand{\BendMax}[1]{\ensuremath B^{\text{end}}_{\text{max}} \left( #1 \right)}

\newcommand{\ginv}[1]{\ensuremath g^{-1}(#1)}
\newcommand{\gT}[2]{\ensuremath g_T\left(#1, \, #2\right)}
\newcommand{\gC}[2]{\ensuremath g_C\left(#1, \, #2\right)}
\newcommand{\Pos}[1]{\ensuremath \left( #1 \right)^+}

\newcommand{\fracpart}[1]{\ensuremath \left\langle #1 \right\rangle}

\newcommand{\Tend}[1]{\ensuremath T^\text{end}_{#1}}
\newcommand{\Bend}[1]{\ensuremath B^\text{end}_{#1}}
\newcommand{\Tstart}[1]{\ensuremath T^\text{start}_{#1}}
\newcommand{\Bstart}[1]{\ensuremath B^\text{start}_{#1}}

\newcommand{\ind}[1]{\ensuremath \mathbbm{1}\left( #1 \right)}

\newcommand{\Sgen}[0]{
    \ensuremath{\calS_{\text{gen}}}
}
\newcommand{\Squeue}[0]{
    \ensuremath{\calS_{\text{queue}}}
}
\newcommand{\Sresult}[0]{
    \ensuremath{\calS_{\text{result}}}
}

\newcommand{\Pnone}[0]{
    \ensuremath{\calP_{\text{all}}}
}
\newcommand{\Pelem}[0]{
    \ensuremath{\calP_{\text{elem}}}
}
\newcommand{\Pnew}[0]{
    \ensuremath{\calP_{\text{new}}}
}

\newcommand{\Pgen}[0]{
    \ensuremath{\calP_{\text{gen}}}
}
\newcommand{\Pqueue}[0]{
    \ensuremath{\calP_{\text{queue}}}
}
\newcommand{\Pinit}[0]{
    \ensuremath{\calP_{\text{init}}}
}
\newcommand{\Presult}[0]{
    \ensuremath{\calP_{\text{result}}}
}

\DeclarePairedDelimiter\floor{\lfloor}{\rfloor}

\usepackage{natbib}
\bibpunct[, ]{(}{)}{,}{a}{}{,}%
\def\bibfont{\small}%

\usepackage{mdframed}
\mdfsetup{skipabove=22pt,skipbelow=22pt}
\newmdenv[
  topline=true,
  bottomline=true,
  rightline=true,
  leftline=true,
  innertopmargin=11pt,
  linewidth=.75pt
]{exbox}

\TheoremsNumberedThrough
\ECRepeatTheorems
\usepackage{placeins}
\usepackage{tcolorbox}

\EquationsNumberedThrough

\begin{document}

\RUNAUTHOR{Jacquillat and Lo}

\RUNTITLE{Subpath-Based Column Generation for Electric Vehicle Routing Problems}

\TITLE{Subpath-Based Column Generation for \\ Electric Vehicle Routing Problems}

\ARTICLEAUTHORS{%
    \AUTHOR{Alexandre Jacquillat and Sean Lo}
    \AFF{
        Operations Research Center and Sloan School of Management, Massachusetts Institute of Technology,  Cambridge, MA.
    }
}

\ABSTRACT{%
Motivated by widespread electrification targets, this paper studies an Electric Vehicle Routing Problem with Time Windows and Nonlinear Charging (EVRPTWNL) that jointly optimizes routing decisions and charging decisions given vehicle capacities, time windows and nonlinear charging. We develop a column generation scheme with a subpath-based label-setting algorithm that decomposes the pricing problem into two phases: (i) generating subpaths between charging stations, and (ii) combining subpaths into paths while optimizing charging decisions in between. We formalize a domination framework to establish the convergence and exactness of the algorithm, and prove that the methodology can solve a range of EVRP variants (e.g., with vehicle capacities, time windows, and nonlinear charging) and relaxation-tightening strategies (e.g., \textit{ng}-relaxations and subset-row cuts). Computational results show improvements over path-based benchmarks in both computational time and solution quality, especially when time windows become wider, when vehicles can perform multiple tasks on a single charge and when vehicles still need to recharge several times across the planning horizon. Ultimately, the methodology can scale to otherwise intractable instances with up to 100 customers, thereby enhancing fleet management capabilities across electrified logistics areas.
}

\KEYWORDS{
vehicle routing,
sustainable operations,
column generation,
dynamic programming
} 

\maketitle

\vspace{-12pt}
\section{Introduction}

The share of electricity in energy use is projected to rise from 20\% to nearly 30\% by 2030 due to the deployment of technologies such as electric vehicles, industrial robots and heat pumps \citep{iea}. Electrification can mitigate the reliance on high-cost energy sources, but charging requirements and reduced asset utilization can also hinder adoption. Thus, large-scale electrification necessitates dedicated analytics and optimization tools to efficiently and reliably deploy electrified technologies into operating systems and processes.

As part of this overarching challenge, this paper studies a capacitated, multi-depot Electric Vehicle Routing Problem with Time Windows and Nonlinear Charging (EVRPTWNL). This problem optimizes the operations of electric vehicles (or machines) that need to visit a set of customers (or perform tasks) within pre-defined time windows and with finite capacities. This general environment can capture practical requirements arising in several domains, including:

\begin{example}[Middle-mile distribution]
In regional distribution, logistics providers transport goods between consolidation centers, depots, and retail facilities across multi-city service areas. Electric medium-duty trucks often cannot complete extended multi-stop routes on a single charge, and therefore need to recharge at logistics facilities during loading and unloading activities or at shared public infrastructure \citep{icct_2022_charging}.
\end{example}
 
\begin{example}[Industrial logistics]
Industrial vehicles circulate repeatedly to support tightly coupled production or infrastructure activities across geographically dispersed industrial sites, such as container ports, batching plants, construction sites, distribution centers, and recycling locations. \citep{boysen_energy_2025}. Electrified equipment can reduce total costs of ownership but high-payload requirements require charging between loading, staging, or waiting periods.
\end{example}

\begin{example}[Energy-constrained last-mile distribution]
Delivery requests materialize across a service region with varying timing requirements. Emerging last-mile deliveries based on light electric vehicles (e.g., bicycles, scooters, drones, robots) operate under strict load and energy limitations \citep{boysen_last-mile_2021}. Each trip consumes a significant portion of battery capacity, requiring intermediate charging at charging pads or micro-depots \citep{boston2025cargo}.
\end{example}

\begin{example}[Robotic fleet management]
Autonomous electric robots are increasingly deployed over spatially distributed environments, spanning baggage handling and aircraft servicing at airports \citep{changi2026autonomous}, deliveries of medical supplies and samples at hospitals \citep{valner_scalable_2022}, and field monitoring or harvesting in agriculture \citep{harik_design_2021}. These systems operate over long horizons with finite battery capacity, requiring frequent recharging between tasks.
\end{example}

Across these applications, fleet operators must coordinate the assignment of tasks to vehicles, the sequence of tasks performed by each vehicle, and the charging schedule to sustain throughput while avoiding service delays and preventing energy shortages. Accordingly, the EVRPTWNL jointly optimizes discrete routing decisions and continuous charging decisions. By design, the model can capture complexities arising in various domains, such as heterogeneous travel or switching costs, vehicle capacities, time window requirements, and nonlinear battery charging dynamics. In particular, the charging decisions introduce an additional degree of freedom in that multiple charging schedules can be feasible for the same same routing decisions, as depicted in Figure~\ref{fig.intro_example}.

\begin{figure}[h!]
    \centering
    \begin{subfigure}[t]{0.48\textwidth}
    \centering
    \begin{tikzpicture}
        \begin{axis}[
            axis equal image,
            width=9cm,
            xmin=0, xmax=28,
            ymin=0, ymax=17,
            xlabel={\small Time},
            ylabel={\small Charge},
            axis lines=middle,
            xtick=\empty,
            ytick=\empty,
        ]
        \node at (1,14) {1};
        \node at (3,11) {2};
        \node at (6,9.5) {3};
        \node at (9,5) {4};
        \node at (14,7) {5};
        \node at (15,4) {6};
        \node at (21.5,5) {7};
        \node at (23.5,0.7) {8};
        
        \addplot [black, thick, no marks] coordinates {
            (0,14) 
            (2,11)
            (3,9.5)
            (5,9.5)
            (8,5)
        };
        \addplot [thick, dashed, domain=8:10] {5*x/2 - 15};
        \addplot [black, thick, no marks] coordinates {
            (10,10) 
            (12,7) 
            (13,7) 
            (14,4)
        };
        \addplot [thick, dashed, domain=14:16] {5*x/2 - 31};
        \addplot [black, thick, no marks] coordinates {
            (16,9)
            (20.5,5)
            (23,0)
        };

        \addplot [black, dotted, no marks] coordinates {
            (10,10) 
            (11,12.5)
            (13,9.5) 
            (14,6.5)
            (15,9)
            (19.5,5)
            (22,0)
        };

        \addplot [only marks, mark=square, fill=black] coordinates {
            (0,14) 
            (8,5)
            (10,10)
            (14,4)
            (16,9)
            (23,0)
        };
        \addplot [only marks, mark=*, fill=black] coordinates {
            (2,11)
            (5,9.5)
            (13,7) 
            (20.5,5)
        };
        \end{axis}
    \end{tikzpicture}
    \end{subfigure}\hfill
    \begin{subfigure}[t]{0.48\textwidth}
    \centering
    \begin{tikzpicture}
        \begin{axis}[
            axis equal image,
            width=9cm,
            xmin=0, xmax=28,
            ymin=0, ymax=17,
            xlabel={\small Time},
            ylabel={\small Charge},
            axis lines=middle,
            xtick=\empty,
            ytick=\empty,
        ]
        \node at (1,14) {1};
        \node at (3,11) {2};
        \node at (6,9.5) {3};
        \node at (9,5) {4};
        \node at (12,9.5) {5};
        \node at (13,6.5) {6};
        \node at (20.5,5) {7};
        \node at (22.5,0.7) {8};
        
        \addplot [black, thick, no marks] coordinates {
            (0,14) 
            (2,11)
            (3,9.5)
            (5,9.5)
            (8,5)
        };
        \addplot [thick, dashed, domain=8:11] {5*x/2 - 15};
        \addplot [black, thick, no marks] coordinates {
            (11,12.5) 
            (13,9.5) 
            (14,6.5)
        };
        \addplot [thick, dashed, domain=14:15] {5*x/2 - 28.5};
        \addplot [black, thick, no marks] coordinates {
            (15,9)
            (19.5,5)
            (22,0)
        };

        \addplot [black, dotted, no marks] coordinates {
            (10,10) 
            (12,7) 
            (13,7) 
            (14,4)
            (16,9)
            (20.5,5)
            (23,0)
        };

        \addplot [only marks, mark=square, fill=black] coordinates {
            (0,14) 
            (8,5)
            (11,12.5) 
            (14,6.5)
            (15,9)
            (22,0)
        };
        \addplot [only marks, mark=*, fill=black] coordinates {
            (2,11)
            (5,9.5)
            (13,9.5)
            (19.5,5)
        };
        \end{axis}
    \end{tikzpicture}
    \end{subfigure}
    \caption{
    Sample sequences of routing and charging decisions, with customers (squares), charging stations or depots (circles), charging actions (dashed curves); travel (solid lines) and waiting due to time windows (solid horizontal lines). The thin dotted lines reproduce the alternative solution on each figure.}
    \label{fig.intro_example}
    \vspace{-10pt}
\end{figure}

We formulate the EVRPTWNL via a path-based set-partitioning model. The model assigns each vehicle to a path that encapsulates a sequence of customer visits and charging stations---but no charging schedule. For any sequence, there may be an infinite number of charging combinations that would guarantee feasibility in view of battery requirements and time windows. For instance, all solutions between the two depicted in Figure~\ref{fig.intro_example} could be considered as part of the same path.

This paper develops a column generation algorithm armed with a novel two-phase, subpath-based label-setting algorithm for the pricing problem. Column generation iterates between a master problem that seeks a path-based solution, and a pricing problem that seeks a path of negative reduced cost or proves that none exists. In the EVRPTWNL, the pricing problem features an NP-hard elementary resource-constrained shortest path structure, and is typically solved via label-setting algorithms that construct paths iteratively from depot to depot (see Section~\ref{sec.literature}). Instead, we develop an exact subpath-based label-setting that further decomposes the pricing problem into subpaths, defined as a sequence of customer visits between charging stations. We find that the algorithm can result in significant computational improvements; along with a master problem heuristic, it yields high-quality solutions and tight optimality gaps in large-scale problems, outperforming path-based label-setting benchmarks. Specifically, this paper makes the following contributions:
\begin{enumerate}[1.]
    \item We develop a new column generation algorithm via a subpath-based label-setting algorithm that decomposes the (path-based) pricing problem into two phases. In the first phase, the algorithm extends subpaths along arcs between depots or charging stations. In the second phase, it extends paths along subpaths, optimizing charging decisions in-between. This algorithm provides a new decomposition scheme for column generation in hierarchical settings with discrete routing decisions and continuous resource management decisions.

    \item We show the exactness of our subpath-based label-setting algorithm via a novel domination framework. We demonstrate that existing label-setting concepts are not sufficient, leading to arbitrarily large optimality gaps. We therefore design a new domination framework relying on traditional forward domination concepts to extend subpaths along arcs, but on forward \textit{and} backward domination concepts to extend paths along subpaths (Property~\ref{property.domination_ours}). We prove that these new domination properties are sufficient to establish the exactness of the subpath-based label-setting algorithm in the pricing problem (Theorem~\ref{thm.exact_PP}). 
    This result is enabled by the separation of routing and charging decisions in our two-phase subpath-based label-setting procedure, which circumvents the infinite number of possible charging sequences.
    
    \item We apply our subpath-based label-setting algorithm to handle a range of electric vehicle routing variants and relaxation-tightening strategies. Specifically, we show how our algorithm can handle vehicle capacities, time windows, nonlinear charging dynamics, as well as \textit{ng}-relaxations and subset-row cuts. For each one, we define appropriate label-setting resources, along with resource extension functions and dominance criteria. The main technical result is to prove that these elements satisfy the subpath-level and path-level domination criteria elicited in Property~\ref{property.domination_ours}. This implies the exactness of the subpath-based label-setting algorithm in the pricing problem, hence the exactness of the column generation procedure for the corresponding relaxations of the EVRPTWNL (Theorem~\ref{thm.EVRPTWNL_correctness}, Propositions~\ref{prop.ngroute_correctness} and~\ref{prop.SRCs_correctness} and Corollaries~\ref{cor.cg_finite_correct}--\ref{cor.iterative_cg_correctness}).

    \item We perform extensive computational benchmarking of our two-phase, subpath-based label-setting algorithm against traditional path-based label-setting benchmarks, using both publicly available instances and new instances inspired by our motivating examples. Results show that our subpath-based label-setting algorithm can significantly accelerate the pricing subproblem, and thereby provide improvements in downstream computation time and solution quality. These benefits become stronger as time windows become wider, battery capacities become smaller, paths become longer (many customers within the planning horizon) and paths are decomposed into more subpaths (charging stations visited frequently). This regime echoes our motivating examples in which limited battery capacities, high payload requirements or spatial dispersion may require multiple charging stops within the planning horizon. Ultimately, the proposed method can scale to otherwise intractable instances with up to 100 customers, thereby enhancing the fleet management capabilities in electrified logistics.
\end{enumerate}

This paper is organized as follows. Section~\ref{sec.literature} reviews the literature on electric vehicle routing and subpath-based decomposition. Section~\ref{sec.model} formulates the EVRPTWNL and the column generation scheme. Section~\ref{sec.labelsetting} describes our two-phase subpath-based label-setting algorithm, applies it to the EVRPTWNL, and proves its correctness. Section~\ref{sec.tighter} augments the algorithm to handle relaxation-tightening strategies---\textit{ng}-relaxations and subset-row cuts---and retrieve a feasible solution via a master problem heuristic. Section~\ref{sec.numerical_results} presents computational results. Section~\ref{sec.conclusion} concludes.

\section{Literature review}
\label{sec.literature}

One thread of the literature on electrified logistics deals with the strategic problem of locating charging stations based on users' routing choices \citep{arslan2019branch}, traffic congestion \citep{kinay2021charging}, car-sharing \citep{brandstatter2020location}, electricity markets \citep{he2013optimal}, and battery swapping \citep{mak2013infrastructure,qi2023scaling}. Another branch optimizes routing operations for a vehicle given the availability of charging stations \citep{sweda2017adaptive}, speed-dependent operations \citep{nejad2016optimal}, or queuing at capacitated charging stations \citep{kullman2021electric}. In-between, our paper falls into the literature on the electric vehicle routing problem (EVRP).

Canonical vehicle routing problems with time windows \citep{kallehauge2005vehicle} and capacitated vehicles \citep{ralphs2003capacitated} link discrete routing decisions and continuous decisions (e.g., timing or load upon each customer visit). In these problems, the continuous dynamics are determined by the routing sequence. The EVRP features an extra degree of freedom to determine where, when and for how long to charge. Early approaches employed clustering-based heuristics \citep{erdogan2012green} or assumptions that vehicles always charge to full \citep{schneider2014electric}. Recent variants include speed-dependent consumption, nonlinear charging \citep{felipe2014heuristic,goeke2015routing,montoya2017electric,fernandez2022arc}, capacitated charging stations \citep{froger2022electric}, public transit \citep{devos2022electric}, and dial-a-ride \citep{molenbruch2023electric}.

Exact EVRP methodologies rely on set-partitioning formulations and column generation, using path-based variables encoding a vehicle's full sequence of routing and charging decisions (Figure~\ref{fig.tikz_path}). To generate path-based variables, the pricing problem features a strongly NP-hard elementary resource-constrained shortest-path structure \citep{dror1994note}. It is typically solved by label-setting algorithms with dedicated domination criteria to encode charging decisions (Figure~\ref{fig.tikz_path_monodirectional}). \cite{desaulniers2016exact} developed a bidirectional label-setting algorithm for EVRP with time windows, using three resources to capture charging requirements. Other resources were used to model the effective range of vehicles under battery swapping \citep{andelmin2017exact} and vehicles' state of charge between customer visits \citep{parmentier2023electric}. \cite{nafstad_branch-price-and-cut_2025} captured non-linear charging with pointwise function comparisons due to heterogeneous charging schedules.

The main bottleneck lies in the exponential growth of the number of paths. Bidirectional label-setting approaches extend paths forward from the source and backward from the sink, as in Figure~\ref{fig.tikz_path_bidirectional} \citep{righini2006symmetry,desaulniers2016exact}. Our two-phase subpath-based label-setting algorithm provides an alternative, hierarchical decomposition of the pricing problem. This relates to the approach from \cite{gschwind2019stabilized} in multi-commodity split-delivery vehicle routing, who first determined bundles of commodities and then constructed routes; and the one from \cite{andelmin2017exact} and \cite{parmentier2023electric} in EVRP, who first generated non-dominated charging arcs between consecutive customers separated by a charging station and then generated paths. In contrast, our algorithm generates subpaths independently between depots or charging stations, and then combines them into paths while determining charging actions (Figure~\ref{fig.tikz_subpath}--\ref{fig.tikz_path_subpath}). By breaking down a path into smaller subpaths, this approach further quells the rate of exponential growth in the pricing problem: rather than scaling exponentially with the number of customers (or with approximately \textit{half} of it for bidirectional label-setting), our approach scales exponentially with \textit{the number of customers per subpath} and \textit{the number of candidate subpaths per path}.

\tikzset{
    sourcesink/.style={rectangle,draw=black,fill=purple,thick,minimum size=3.5mm},
    charger/.style={rectangle,draw=black,fill=yellow,thick,minimum size=3.5mm},
    customer/.style={circle,draw=black,thick,minimum size=3.5mm},
    >={Stealth},
    legend/.style={right,inner sep=2mm,text height=1.5ex,text depth=0.25ex}
}
\begin{figure}[h!]
    \centering
    \caption{Our two-phase subpath-based label-setting algorithm versus path-based label-setting benchmarks.}
    \label{fig.tikz_compare}
    \begin{subfigure}{\textwidth}
        \centering
        \begin{tikzpicture}
        [
            scale=0.7,
            auto,
            pin distance=2mm,
        ]
            \node[sourcesink]   (d1)    at ( 0.0, 0.0) {};
            \node[customer]     (c11)    at ( 1.5,-0.5) {};
            \node[customer]     (c12)    at ( 2.0, 1.0) {};
            \node[customer]     (c13)    at ( 3.0, 0.0) {};
            \node[charger]  
                                (r1)    at ( 4.0, 1.0) {};
            \node[customer]     (c21)    at ( 5.0, 0.0) {};
            \node[customer]     (c22)    at ( 5.5, 1.5) {};
            \node[customer]     (c23)    at ( 6.5, 0.5) {};
            \node[customer]     (c24)    at ( 6.5,-1.0) {};
            \node[charger]  
                                (r2)    at ( 8.0,-1.0) {};
            \node[customer]     (c31)    at ( 9.0, 0.5) {};
            \node[customer]     (c32)    at (10.0,-0.5) {};
            \node[customer]     (c33)   at (10.5, 1.5) {};
            \node[customer]     (c34)   at (11.0, 0.5) {};
            \node[charger]      (r3)    at (12.0, 0.5) {};
            \node[customer]     (c41)   at (12.5,-0.5) {};
            \node[customer]     (c42)   at (13.5, 1.5) {};
            \node[customer]     (c43)   at (14.0, 0.0) {};
            \node[customer]     (c44)   at (15.0,-1.0) {};
            
            \node[customer]     (c45)   at (15.0, 1.5) {};
            \node[sourcesink]   (d2)    at (16.0, 0.0) {};
        
            \draw[-]        (d1)    to [bend right=15] 
                                (c11);
            \draw[-]        (c11)   to [bend right=15]
                                (c13);
            \draw[-]        (c13)   to (r1);
            \draw[-]        (r1)    to [bend left=15] 
                                (c22);
            \draw[-]        (c22)   to [bend left=30] 
                                (c23);
            \draw[-]        (c23)   to [out=270,in=135]
                                (c24);
            \draw[-]        (c24)  to [bend right=30] 
                                (r2);
            \draw[-]        (r2)    to [out=45,in=225]
                                (c31); 
            \draw[-]        (c31)   to [bend left=30] 
                                (c33);
            \draw[-]        (c33)   to [bend left=30] 
                                (r3);
            \draw[-]        (r3)    
                                to [bend right=15] (c41); 
            \draw[-]        (c41)
                                to [bend right=30] (c44);
            \draw[->]       (c44)
                                to [bend right=15] (d2);

            \node[sourcesink]   (d0)    at (17.5, 1.4) {};
            \node[charger]      (r0)    at (17.5, 0.7) {};
            \node[customer]     (c0)    at (17.5, 0.0) {};
            \node[legend,xshift=3pt]   (d0)    at (d0) {\small Depot};
            \node[legend,xshift=3pt]   (r0)    at (r0) {\small Charging station};
            \node[legend,xshift=3pt]   (c0)    at (c0) {\small Customer};
        \end{tikzpicture}%
        \caption{\small A path starting and ending at depots, visiting customers and recharging periodically in between.}
        \label{fig.tikz_path}
    \end{subfigure}%
    \newline
    \vspace{1em}
    \begin{subfigure}{\textwidth}
        \centering
        \begin{tikzpicture}
        [
            scale=0.8,
            auto,
            pin distance=2mm,
            inner sep=1mm,
        ]
            \node[sourcesink]   (d1)    at ( 0.0, 0.0) {};
            \node[customer]     (c11)    at ( 1.5,-0.5) {};
            \node[customer]     (c12)    at ( 2.0, 1.0) {};
            \node[customer]     (c13)    at ( 3.0, 0.0) {};
            \node[charger]  
                                (r1)    at ( 4.0, 1.0) {};
            \node[customer]     (c21)    at ( 5.0, 0.0) {};
            \node[customer]     (c22)    at ( 5.5, 1.5) {};
            \node[customer]     (c23)    at ( 6.5, 0.5) {};
            \node[customer]     (c24)    at ( 6.5,-1.0) {};
            \node[charger]  
                                (r2)    at ( 8.0,-1.0) {};
            \node[customer]     (c31)    at ( 9.0, 0.5) {};
            \node[customer]     (c32)    at (10.0,-0.5) {};
            \node[customer]     (c33)   at (10.5, 1.5) {};
            \node[customer]     (c34)   at (11.0, 0.5) {};
            \node[charger]      (r3)    at (12.0, 0.5) {};
            \node[customer]     (c41)   at (12.5,-0.5) {};
            \node[customer]     (c42)   at (13.5, 1.5) {};
            \node[customer]     (c43)   at (14.0, 0.0) {};
            \node[customer]     (c44)   at (15.0,-1.0) {};
            
            \node[customer]     (c45)   at (15.0, 1.5) {};
            \node[sourcesink]   (d2)    at (16.0, 0.0) {};
        
            \draw[-]        (d1)    to [bend right=15] 
                                (c11);
            \draw[-]        (c11)   to [bend right=15]
                                (c13);
            \draw[-]        (c13)   to (r1);
            \draw[-]        (r1)    to [bend left=15] 
                                (c22);
            \draw[-]        (c22)   to [bend left=30] 
                                (c23);
            \draw[-]        (c23)   to [out=270,in=135]
                                (c24);
            \draw[-]        (c24)  to [bend right=30] 
                                (r2);
            \draw[->]       (r2)    to [out=45,in=225] 
                                node [inner sep=1mm] {$p$} 
                                (c31); 
            \draw[->,dashed](c31)   to [bend left=15] 
                                node [swap,inner sep=1mm] {$a$}
                                (c33); 
            \draw[->,dashed](c31)   to [bend left=45] 
                                node [swap,inner sep=1mm] {$a'$}
                                (c32);
        \end{tikzpicture}%
        \caption{\small Monodirectional path-based label-setting.}
        \label{fig.tikz_path_monodirectional}
    \end{subfigure}%
    \newline
    \vspace{1em}
    \begin{subfigure}{\textwidth}
        \centering
        \begin{tikzpicture}
        [
            scale=0.8,
            auto,
            pin distance=2mm,
            inner sep=1mm,
        ]
            \node[sourcesink]   (d1)    at ( 0.0, 0.0) {};
            \node[customer]     (c11)    at ( 1.5,-0.5) {};
            \node[customer]     (c12)    at ( 2.0, 1.0) {};
            \node[customer]     (c13)    at ( 3.0, 0.0) {};
            \node[charger]  
                                (r1)    at ( 4.0, 1.0) {};
            \node[customer]     (c21)    at ( 5.0, 0.0) {};
            \node[customer]     (c22)    at ( 5.5, 1.5) {};
            \node[customer]     (c23)    at ( 6.5, 0.5) {};
            \node[customer]     (c24)    at ( 6.5,-1.0) {};
            \node[charger]  
                                (r2)    at ( 8.0,-1.0) {};
            \node[customer]     (c31)    at ( 9.0, 0.5) {};
            \node[customer]     (c32)    at (10.0,-0.5) {};
            \node[customer]     (c33)   at (10.5, 1.5) {};
            \node[customer]     (c34)   at (11.0, 0.5) {};
            \node[charger]      (r3)    at (12.0, 0.5) {};
            \node[customer]     (c41)   at (12.5,-0.5) {};
            \node[customer]     (c42)   at (13.5, 1.5) {};
            \node[customer]     (c43)   at (14.0, 0.0) {};
            \node[customer]     (c44)   at (15.0,-1.0) {};
            
            \node[customer]     (c45)   at (15.0, 1.5) {};
            \node[sourcesink]   (d2)    at (16.0, 0.0) {};
        
            \draw[-]        (d1)    to [bend right=15] 
                                (c11);
            \draw[-]        (c11)   to [bend right=15]
                                (c13);
            \draw[-]        (c13)   to (r1);
            \draw[-]        (r1)    to [bend left=15] 
                                (c22);
            \draw[-]        (c22)   to [bend left=30]
                                node [swap] {$\overrightarrow{p}$} 
                                (c23);
            \draw[->,dashed](c23)   to [out=270,in=135]
                                node [swap,inner sep=1mm] {$\overrightarrow{a}$}
                                (c24);
            \draw[->,dashed](c23)   to [out=270,in=135]
                                node [inner sep=1mm] {$\overrightarrow{a}'$}
                                (r2);

            \draw[<-,dashed](c31)   
                                to [bend left=30] 
                                node [inner sep=1mm] {$\overleftarrow{a}$}
                                (c33);
            \draw[<-,dashed](c32)   
                                to [bend left=30] 
                                node [swap,inner sep=1mm] {$\overleftarrow{a}'$}
                                (c33);
            \draw[-]        (c33)   
                                to [bend left=30] (r3);                                
            \draw[-]        (r3)
                                to [bend right=15] (c41); 
            \draw[-]        (c41)
                                to [bend right=30] (c44);
            \draw[-]        (c44)
                                to [bend right=15] 
                                node [swap] {$\overleftarrow{p}$} 
                                (d2);
        \end{tikzpicture}%
        \caption{\small Bi-directional path-based label-setting: extending forward partial paths and backward partial paths.}
        \label{fig.tikz_path_bidirectional}
    \end{subfigure}%
    \vspace{1em}
    \newline
    \begin{subfigure}{0.3\textwidth}
        \centering
        \begin{tikzpicture}
        [
            scale=0.8,
            auto,
            pin distance=2mm,
        ]
            \node[sourcesink]   (d1)    at ( 0.0, 0.0) {};
            \node[customer]     (c11)   at ( 1.5,-0.5) {};
            \node[customer]     (c12)   at ( 2.0, 1.0) {};
            \node[customer]     (c13)   at ( 3.0, 0.0) {};
            \node[charger]      (r1)    at ( 4.0, 1.0) {};
        
            \draw[-]            (d1) to [bend right=15] 
                                    (c11);
            \draw[->]           (c11) to [bend right=15] 
                                    node [swap,inner sep=1mm] {$s$}
                                    (c13);
            \draw[->,dashed]    (c13) to [out=90,in=0] 
                                    node [swap,inner sep=1mm] {$a'$} 
                                    (c12);
            \draw[->,dashed]    (c13) to node [swap,inner sep=1mm] {$a$} 
                                    (r1);
        \end{tikzpicture}%
        \caption{\small Our two-phase label-setting: generating subpaths.}
        \label{fig.tikz_subpath}
    \end{subfigure}%
    \begin{subfigure}{0.70\textwidth}
        \centering
        \begin{tikzpicture}
        [
            scale=0.7,
            auto,
            pin distance=2mm,
        ]
            
            \node[sourcesink]   (d1)    at ( 0.0, 0.0) {};
            \node[customer]     (c11)   at ( 1.5,-0.5) {};
            \node[customer]     (c12)   at ( 2.0, 1.0) {};
            \node[customer]     (c13)   at ( 3.0, 0.0) {};
            \node[charger,pin=120:$\tau_1$]  
                                (r1)    at ( 4.0, 1.0) {};
            \node[customer]     (c21)   at ( 5.0, 0.0) {};
            \node[customer]     (c22)   at ( 5.5, 1.5) {};
            \node[customer]     (c23)   at ( 6.5, 0.5) {};
            \node[customer]     (c24)   at ( 6.5,-1.0) {};
            \node[charger,pin=120:$\tau_2$,pin=300:$\tau'_2$]  
                                (r2)    at ( 8.0,-1.0) {};
            \node[customer]     (c31)   at ( 9.0, 0.5) {};
            \node[customer]     (c32)   at (10.0,-0.5) {};
            \node[customer]     (c33)   at (10.5, 1.5) {};
            \node[customer]     (c34)   at (11.0, 0.5) {};
            \node[charger]      (r3)    at (12.0, 0.5) {};
            \node[customer]     (c41)   at (12.5,-0.5) {};
            \node[customer]     (c42)   at (13.5, 1.5) {};
            \node[customer]     (c43)   at (14.0, 0.0) {};
            \node[customer]     (c44)   at (15.0,-1.0) {};
            
            \node[customer]     (c45)   at (15.0, 1.5) {};
            \node[sourcesink]   (d2)    at (16.0, 0.0) {};
        
            \draw[-]        (d1)  
                                to [bend right=15] 
                                (c11);
            \draw[-]        (c11)   
                                to [bend right=15]
                                node [swap,inner sep=1mm] {$s_1$}
                                (c13);
            \draw[->]       (c13)  to (r1);
            \draw[-]        (r1)  to [bend left=15] (c22);
            \draw[-]        (c22)  to [bend left=30] (c23);
            \draw[-]        (c23)  to [out=270,in=135]
                                node [swap,inner sep=1mm] {$s_2$}
                                (c24);
            \draw[->]       (c24)  to [bend right=30] (r2);
            \draw[thick,decorate,decoration={brace,amplitude=10pt,mirror}] 
                (0,-0.7) -- (8,-1.7)
                node [midway,swap,inner sep=2mm,yshift=-3pt] {$p = \{s_1, s_2\}$};
            \draw[-,dashed] (r2)  to [out=45,in=225] 
                                node [swap,inner sep=1mm] {$s_3$} 
                                (c31); 
            \draw[-,dashed] (c31)  to [bend left=30] (c34);
            \draw[->,dashed](c34) to [bend right=15] (r3);
            \draw[-,dashed] (r2)  to [out=0,in=240] 
                                node [swap,inner sep=1mm] {$s'_3$}
                                (c32);
            \draw[-,dashed] (c32)  to [out=75,in=180] (c34);
            \draw[->,dashed](c34) to [bend left=30] (r3);
        \end{tikzpicture}%
        \caption{\small Our two-phase label-setting: combining subpaths into paths.}
        \label{fig.tikz_path_subpath}
    \end{subfigure}%
\end{figure}

At the same time, this approach introduces methodological challenges to guarantee convergence and exactness. In response, we formalize new subpath- and path-based domination properties. In particular, we prove that Property~\ref{property.domination_ours}\ref{property.domination_pss} is critical to propagate domination between subpaths ($s_1 \succeq s_2$) into domination between partial paths ($p \oplus s_1 \succeq p \oplus s_2$). We also show that this property is automatically satisfied when resources satisfy monotonicity conditions (such as the cost and duration of refuel paths in \cite{andelmin2017exact}, and the weight and cost of multiedges in \cite{gschwind2019stabilized}) but necessary in our context to handle non-additive resources capturing time windows and \textit{ng}-sets. Altogether, this paper contributes a novel subpath-based decomposition of the pricing problem for column generation with hierarchical structure, and a new rigorous domination framework establishing the exactness of the resulting subpath-based label-setting algorithm.

\paragraph{Extended formulations.}
Our decomposition approach also relates to subpath-based extended formulations in combinatorial optimization. In pickup-and-delivery, \cite{alyasiry2019exact} and \cite{zhang2023routing} optimized over subpaths encapsulating sequences of pickups and dropoffs between points where the vehicle is empty. \cite{rist2021new} optimized over subpaths encapsulating sequences of consecutive pickups or dropoffs. \cite{cummings2024deviated} optimized flexible segments between mandatory stops. 
Recent papers applied column generation to generate subpath-based variables iteratively \citep{hasan2021benefits,rist2022column,cummings2024deviated,jacquillat2022optimizing}. In contrast, our methodology relies on a path-based formulation and adds path-based variables iteratively, but further decomposes the pricing problem into subpaths. 
\cite{rossum_adaptive_2026} also apply a subpath-based pricing problem decomposition, and avoids enumerating non-dominated subpaths through a bucket-based approach.

In our context, a subpath-based formulation could rely on variables encoding the sequence of customer visits between charging stations \citep{lee_exact_2021}, which would require linking constraints to optimize charging decisions in the master problem and complicate the column generation algorithm. An alternative one would encode the sequence of customer visits between charging stations as well as charging times, with flow balance constraints combining subpaths into paths in the master problem. Such a formulation would mitigate the exponential growth in the number of variables but its benefits may be limited if subpaths cannot be easily shared across paths \citep{ceselli_branch-and-cut-and-price_2021,bezzi_route-based_2023}. For instance, subpath $s_1$ could end at a charging station with 80\% battery whereas subpath $s_2$ could start with 70\% battery, and alleviating such mismatches may require extensive column generation iterations. Instead, we retain a path-based formulation, leading to a more complex pricing problem---and motivating our subpath-based decomposition algorithm.

\section{Model formulation and column generation preliminaries}
\label{sec.model}

We address a capacitated, multi-depot Electric Vehicle Routing Problem with Time Windows and Nonlinear Charging (EVRPTWNL). We define the problem in Section~\ref{subsec.statement} and a set partitioning formulation in Section~\ref{subsec.formulation}; we also provide an explicit arc-based benchmark in~\ref{app.arc_formulation}. Section~\ref{ssec.CG} concludes by outlining a column generation algorithm and highlighting its challenges to manage the size of the pricing problem and tighten the relaxation, which motivate our methodology.

\subsection{Problem statement}
\label{subsec.statement}

We consider a fleet of $K$ homogeneous electric vehicles with the same battery $B$, the same capacity $D$, the same travel times and costs, and the same battery charging and depletion dynamics. We represent operations on a directed graph $(\calV, \calA)$ over a time horizon $T$. Nodes are partitioned into sets of depots $\Depots$, customers $\Custs$, and charging stations $\Stations$, such that $\calV = \Custs \cup \Depots \cup \Stations$. 
Arcs allow self-loops only at depots: $\calA = \Set{(i, j) | i, j \in \Depots} \cup \Set{(i, j) | i \neq j \in \calV}$.
Each arc $(i, j) \in \calA$ is associated with a travel time $t_{i,j}$, a cost $c_{i,j}$, and a battery utilization $b_{i,j}$ (all strictly positive when $i \neq j$). Each customer has a demand $d_i$ and needs to be served within a time window $[\alpha_i, \beta_i]$. To extend our notation to depots and charging stations $i \in \Depots \cup \Stations$, we let $d_i = 0$, $\alpha_i=0$ and $\beta_i=T$.

Each vehicle starts in a depot in $\Depots$ with full charge, visits customers in $\Custs$, recharges in charging stations in $\Stations$, and ends in a depot. We let $v^{\text{start}}_j$ denote the number of vehicles starting at each depot $j$. We also impose a minimum number of vehicles $v^{\text{end}}_j$ ending in each $j \in \Depots$, for instance to enable maintenance requirements or facilitate future operations \citep{gopalan_aircraft_1998,hemmelmayr_2013_waste,rahimi-vahed_2015_fleet,birolini2023day,jacquillat2026iterative}. This choice can capture single-depot instances from \cite{desaulniers2016exact}, multi-depot instances where each vehicle ends at its starting depot, and more general multi-depot instances in which each vehicle can start and end in any depot. We also assume that charging stations are uncapacitated.

A charging schedule $g: \R_+ \rightarrow [0, B]$ defines the charge gained as a function of charging time charging from a depleted battery. We assume that $g$ is concave and increasing, with $g(0) = 0$ \citep{lam_branch-and-cut-and-price_2022,lera-romero_branch-cut-and-price_2023,nafstad_branch-price-and-cut_2025}. We define $\gT{b_1}{b_2}=\ginv{b_2} - \ginv{b_1}$ as the amount of time required to charge from $b_1$ to $b_2$. We also define the end state of charge $\gC{b}{t}= g(\ginv{b} + t)$ as a function of the current charge $b$ and the charging time $t$.

The EVRPTWNL seeks a path for each vehicle to minimize travel costs, while ensuring that vehicle loads lie below capacity, that all customers are served within their time windows, and that no vehicle runs out of battery. An arc-based integer optimization formulation is provided in~\ref{app.arc_formulation}.

\subsection{Path-based set-partitioning formulation}
\label{subsec.formulation}

The core complexity of the EVRPTWNL is to link discrete routing decisions with continuous decisions governing the load, time and charge upon each customer visit. In particular, the charging dynamics introduce an extra degree of freedom to determine where, when and for how long to charge. The arc-based formulation in~\ref{app.arc_formulation} manages these interdependencies via ``big-$M$'' coupling constraints, which induce weak relaxations and hinder the scalability of branch-and-cut algorithms. Instead, we define a path-based formulation using Dantzig-Wolfe decomposition principles.

We formalize a path in Definition~\ref{def.path} as the sequence of routing decisions for a vehicle through the planning horizon. These decisions uniquely determine the load upon each customer visit, which determines path feasibility in view of the vehicle capacity constraints. However, the charging schedule and the resulting timing of each customer visit are not entirely governed by the routing decisions. In turn, a path is feasible if \textit{there exists} a charging schedule such that the vehicle does not run out of battery and can visit all assigned customers within their time windows.

\begin{definition}[Path]
    \label{def.path}
    A \textit{path} $p$ is defined by a node sequence $N(p) = (n_0, n_1, n_2, \dots, n_m)$ such that $(n_0, n_1), (n_1, n_2), \dots, (n_{m-1}, n_m) \in \calA$, 
    with starting node 
    $\snode{p} = n_0 \in \Depots$, 
    intermediate nodes 
    $n_1, \dots, n_{m-1} \in \Custs \cup \Stations$, 
    and ending node 
    $\enode{p} = n_m \in \Depots$.
    The binary parameters $\sigma_{i,\text{start}}^p$ and  $\sigma_{i,\text{end}}^p$ denote whether path $p$ starts and ends at depot $i \in \Depots$ respectively.
    The parameter $\serve{i}{p}$ stores the number of times path $p$ visits customer $i\in\Custs$: $\serve{i}{p} = | \Set{k \in \{0, \dots, m\} | n_k = i} |$.
    The parameter $D(p)$ stores total load:
    $D(p) = \sum_{k=0}^m d_{n_{k}}$.
    Its cost is the sum of travel costs of its constituent arcs: $c^p := \sum_{\ell=0}^{m-1} c_{n_\ell,n_{\ell+1}}$. A \textit{partial path} relaxes the condition $\enode{p} \in \Depots$.

    Path $p$ with node sequence $N(p) = (n_0, n_1, n_2, \dots, n_m)$ is feasible if there exist charging times $\Set{\tau_k | n_k \in \Stations }$ and of corresponding timing variables $t_k$ and charge variables $b_k$ for $k\in\{0, \dots, m\}$, defined recursively as follows, such that
    \begin{enumerate*}[(i)]
        \item 
        \label{def.path_feasible_load}
        $D(p) \leq D$ (vehicle capacity);
        \item 
        \label{def.path_feasible_charge}
        $b_k \in [0, B]$ for all $k\in\{0, \dots, m\}$ (battery); and
        \item 
        \label{def.path_feasible_time}
        $t_{k} \in [\alpha_{n_{k}}, \beta_{n_{k}}]$ for all $k \in \{0, \dots, m\}$ (time windows).
    \end{enumerate*}
    \begin{alignat}{3}
        \label{eq.propagation_time}
        t_0 = 0,
        \quad 
        \text{ and } \forall \ k \in \{0, \dots, m-1\}: 
        t_{k+1} & = \begin{cases}
            \max \{ t_k + \tau_k + t_{n_k,n_{k+1}}, \alpha_{n_{k+1}} \}
            & \text{ if } n_k \in \Stations
            \\
            \max \{ t_k + t_{n_k,n_{k+1}}, \alpha_{n_{k+1}} \}
            & \text{ otherwise.}
        \end{cases}
        \\
        \label{eq.propagation_charge}
        b_0 = B,
        \quad 
        \text{ and } \forall \ k \in \{0, \dots, m-1\}: 
        b_{k+1} & = \begin{cases}
            \gC{b_k}{\tau_k} - b_{n_k,n_{k+1}} 
            & \text{ if } n_k \in \Stations
            \\
            b_k - b_{n_k,n_{k+1}}
            & \text{ otherwise.}
        \end{cases}
    \end{alignat}
\end{definition}

We let $z^p$ be a binary decision variable denoting if path $p \in \calP$ is selected. The EVRPTWNL minimizes costs (Equation~\eqref{eq.path.obj}) 
while enforcing vehicles' start and end location requirements (Equations~\eqref{eq.path.start_depot} and~\eqref{eq.path.end_depot})
and ensuring that each customer is served once (Equation~\eqref{eq.path.serve_customers}). 
We refer to this problem as $\EVRP(\calP)$, to its optimum as $\OPT(\calP)$, to its linear relaxation as $\EVRPLP(\calP)$, and to its linear bound as $\OPTLP(\calP)$. Note that this problem is parametrized by a set of feasible paths $\calP$, which can depend on the state-space relaxations considered (e.g. \textit{ng}-routes in Section~\ref{ssec.ng}).
\begin{alignat}{3}
    \min \quad 
    & \sum_{p \in \calP} c^p z^p
    \label{eq.path.obj}
    \\
    \st \quad 
    & \sum_{\substack{
        p \in \calP
    }} \sigma_{j,\text{start}}^p z^p = v^{\text{start}}_j
    &\quad& \forall \ j \in \Depots
    \label{eq.path.start_depot}
    \\
    & \sum_{\substack{
        p \in \calP
    }} \sigma_{j,\text{end}}^p z^p \geq v^{\text{end}}_j
    &\quad& \forall \ j \in \Depots
    \label{eq.path.end_depot}
    \\
    & \sum_{p \in \calP} \serve{i}{p} z^p = 1
    &\quad& \forall \ i \in \Custs
    \label{eq.path.serve_customers}
    \\
    & z^p \in \{0, 1\}, 
    \ \forall \ p \in \calP
    \label{eq.path.nonneg_integer}
\end{alignat}

\subsection{Column generation algorithm}
\label{ssec.CG}

Column generation (Algorithm~\ref{alg.CG} in \ref{app.cg}) solves the linear relaxation $\EVRPLP(\calP)$. The master problem generates a feasible relaxation solution based on a subset of path-based variables, stored in $\calP_{\ell}$ at iteration $\ell$. The pricing problem seeks a set $\Pnew$ of variables with negative reduced cost or proves that none exists; it minimizes the reduced cost of variable $z^p$ across all paths $p\in\calP$, defined as:
\begin{equation}\label{eq.RC}
    \cbar{p}
    := c^p
    - \kappa_{\snode{p}}
    - \mu_{\enode{p}}
    - \sum_{i \in \Custs} \serve{i}{p} \nu_i,
\end{equation}
where $\bkappa$, $\bmu$, and $\bnu$ denote the dual variables associated with Equations~\eqref{eq.path.start_depot},~\eqref{eq.path.end_depot} and~\eqref{eq.path.serve_customers}, respectively.

Column generation converges to the optimal relaxation bound $\OPTLP(\calP)$ due to the finite number of path-based variables. This hinges on the definition of a path as a sequence of nodes that excludes charging decisions (Definition~\ref{def.path}), which retrieves a discrete path-based structure. In turn, the main difficulty in the EVRPTWNL lies in the pricing problem, which needs to jointly optimize over the routing decisions and identify a charging schedule to maintain feasibility in view of the battery and time window constraints. To circumvent the exponential growth in the number of paths, we propose a subpath-based label-setting algorithm that first generates subpaths between charging stations and then combines subpaths into paths while optimizing charging decisions in between.

Upon termination of column generation, we augment the subpath-based label-setting algorithm tighten the relaxation bound via adaptive \textit{ng}-relaxations \citep{baldacci2011new} and subset-row cuts \citep{jepsen2008subset}, which requires new label-setting resources, new resource extension functions, and new dominance criteria. We also retrieve a feasible primal solution by solving the terminal master problem with integrality constraints. 
The overall framework (described in Algorithm~\ref{alg.AdaptiveColumnGenerationWithCuts} in Section~\ref{sec.tighter}) is evaluated heuristically because no branching scheme is implemented. However, it can be easily embedded into a branch-and-price-and-cut scheme 
\citep{barnhart1998branch} to recover an exact  solution using the same branching rules as in \cite{desaulniers2016exact}, or into alternative heuristic approaches \citep{eveborn_branch-price-and-cut_2026}. Nonetheless, our computational results yield provably high-quality solutions upon termination with narrow optimality gaps.
\section{A subpath-based pricing problem decomposition}
\label{sec.labelsetting}

The pricing problem in Step 2 of Algorithm~\ref{alg.CG} seeks a shortest path from depot to depot subject to path feasibility requirements, with modified edge costs $\{\overline{c}_{i,j}\:(i,j) \in \calA\}$ given as follows. By construction, the reduced cost $\cbar{p}$ of path $p$ is the sum of $\overline{c}_{i,j}$ over the arcs in $N(p)$, with:
\begin{equation}
    \overline{c}_{i,j} = 
    c_{i,j}
    - \ind{i \in \Depots} \kappa_i
    - \ind{j \in \Depots} \mu_j
    - \ind{j \in \Custs} \nu_j
\end{equation}

The pricing problem can be solved with label-setting algorithms by extending partial paths along arcs, and applying dominance criteria to prune partial paths that can no longer be extended into a path of minimal reduced cost \citep{irnich2005shortest}. The extension is defined as follows.
\begin{definition}[Arc extensions of paths]
    \label{def.arc_extensions_paths}
    Let $p$ be a partial path. 
    An arc $a = (n, n') \in \mathcal{A}$ is an arc extension of $p$ if $\enode{p}=n$. 
    The extended partial path is denoted $p \oplus a$.
\end{definition}
A domination criterion would ideally encode the fact that partial path $p_2$ is dominated by $p_1$ (written $p_1 \succeq p_2$) if every possible completion $q$ of $p_2$ could instead be applied to $p_1$ and lead to a weakly lower reduced cost. Label-setting algorithms relax that requirement by merely ensuring that any arc extension of $p_2$ remains dominated by the same extension of $p_1$. Thus, the following domination property is sufficient to establish the correctness of (path-based) label-setting algorithms.
\begin{property}[Domination for paths]
    \label{property.domination_benchmark}
    Let $p_1, p_2$ be feasible partial paths with $p_1 \succeq p_2$, and $a \in \mathcal{A}$ extend $p_1$ and $p_2$. Suppose $p_2 \oplus a$ is feasible. Then $p_1 \oplus a$ is feasible too, and $p_1 \oplus a \succeq p_2 \oplus a$.
\end{property}

\cite{desaulniers2016exact} and \cite{nafstad_branch-price-and-cut_2025} develop path-based label-setting algorithms (Figure~\ref{fig.tikz_path_monodirectional}) for the EVRPTW and the EVRPTWNL, respectively, by associating each partial path $p$ with a set of resources, propagating them along arcs via resource extension functions (REFs), and defining dominance criteria from pairwise resource comparisons. The main element to prove the correctness of the algorithm is to show that these resources and REFs satisfy Property~\ref{property.domination_benchmark}.

Instead, our methodology first generates subpaths between charging stations and then combines them into paths while optimizing charging decisions in between. As we shall see, a naive extension of Property~\ref{property.domination_benchmark} no longer guarantees optimality of the pricing problem, and can lead to arbitrarily large optimality gaps (Proposition~\ref{prop.arbitrarily_large}). In response, we define new domination properties for subpaths and paths in our subpath-based label-setting algorithm (Section~\ref{ssec.labelsetting_2LS}). We apply this domination framework to handle the EVRPTWNL elements---namely, vehicle capacities, time windows, and nonlinear charging---by defining corresponding resources, REFs and domination criteria at the subpath level in the first phase of the algorithm (Section~\ref{ssec.labelsetting_subpath}) and at the path level in the second phase (Section~\ref{ssec.labelsetting_path}). The main technical results in Section~\ref{ssec.labelsetting_correctness} show that these resources satisfy the domination properties, and that these properties are sufficient to guarantee the exactness of our subpath-based label-setting algorithm in the pricing problem. All proofs are provided in~\ref{app.proofs_4}.

\subsection{From subpaths to paths: domination properties}
\label{ssec.labelsetting_2LS}

We define a subpath as the customer visits between two consecutive non-customer nodes (depots or charging stations). Thus, a subpath is constituted by a sequence of arcs, and each path is constituted by a sequence of subpaths along with charging decisions in between. 

\begin{definition}[Subpath]
    \label{def.subpath}
    A subpath $s$ is defined by a node sequence $N(s) = (n_0, n_1, n_2, \dots, n_m)$ such that $(n_0, n_1), (n_1, n_2), \dots, (n_{m-1}, n_m) \in \calA$, 
    with starting node $\snode{s} = n_0 \in \Depots \cup \Stations$, 
    intermediate nodes $n_1, \dots, n_{m-1} \in \Custs$,
    and ending node $\enode{s} = n_{m} \in \Depots \cup \Stations$.
    The parameter $\serve{i}{s}$ stores the number of times subpath $s$ visits customer $i\in\Custs$: $\serve{i}{s} = | \Set{k \in \{0, \dots, m\} | n_k = i} |$. 
    The parameter $D(s)$ stores total load: $D(s) = \sum_{k=0}^m d_{n_k}$.
    A \textit{partial subpath} relaxes the condition $\enode{s} \in \Depots \cup \Stations$.

    Subpath $s$ with node sequence $N(s) = (n_0, n_1, n_2, \dots, n_m)$ 
    is feasible, if 
    there exist $(t_0, b_0)$ such that, defining timing variables $t_k$ and charge variables $b_k$ for $k \in \{1, \dots, m\}$ recursively from $t_0$ and $b_0$ according to Equations~\eqref{eq.propagation_time}--\eqref{eq.propagation_charge}, we have
    \begin{enumerate*}[(i)]
        \item $D(s) \leq D$ (vehicle capacity);
        \item $b_k \in [0, B]$ for all $k \in \{0, \dots, m\}$ (battery); and
        \item $t_k \in [\alpha_{n_k}, \beta_{n_k}]$ for all $k \in \{0, \dots, m\}$ (time windows).
    \end{enumerate*}
\end{definition}

The algorithm first extends partial subpaths along arcs, and then extends partial paths along subpaths. These extensions are written as follows.
\begin{definition}[Arc extensions of subpaths]
    \label{def.arc_extensions_subpaths}
    Let $s$ be a partial subpath. 
    An arc $a = (n, n') \in \mathcal{A}$ is an arc extension of $s$ if $\enode{s}=n$. 
    The extended partial subpath is denoted $s \oplus a$.
\end{definition}
\begin{definition}[Subpath extensions of paths]
    \label{def.subpath_extensions_paths}
    Let $p$ be a partial path. 
    A subpath $s$ is a subpath extension of $p$ if $\enode{p} = \snode{s} \in \Stations$. 
    The extended partial path is denoted $p \oplus s$.
\end{definition}

Per Bellman's principle of optimality, if a path $p = s_1 \oplus \dots \oplus s_m$ is non-dominated, so are each of its constituting subpaths $s_1, \dots, s_m$ and each partial path $p_k = s_1 \oplus \dots \oplus s_k$ (Lemmas~\ref{lemma.subpath_truncate_nondom} and~\ref{lemma.path_truncate_nondom}). This motivates our approach to find non-dominated subpaths and extend them toward paths; we write $s_1 \succeq s_2$ and $p_1 \succeq p_2$ to indicate that partial subpath $s_1$ dominates $s_2$ and that partial path $p_1$ dominates $p_2$. We define the following three domination properties, which provide the fundamental criteria to prune dominated partial subpaths and dominated partial paths.
\begin{property}[Dominance for paths and subpaths]
    \label{property.domination_ours}
    \hfill
    \begin{enumerate}[(i),noitemsep,topsep=0pt]
        \item 
        \label{property.domination_ssa}
        Let $s_1, s_2$ be feasible partial subpaths with $s_1 \succeq s_2$, and let $a \in \mathcal{A}$ be an arc extending $s_1$ and $s_2$. Suppose $s_2 \oplus a$ is feasible. Then $s_1 \oplus a$ is feasible too, and $s_1 \oplus a \succeq s_2 \oplus a$.
        \item 
        \label{property.domination_pps}
        Let $p_1, p_2$ be feasible partial paths with $p_1 \succeq p_2$, and let $s$ be a subpath extending $p_1$ and $p_2$. Suppose $p_2 \oplus s$ is feasible. Then $p_1 \oplus s$ is feasible too, and $p_1 \oplus s \succeq p_2 \oplus s$.
        \item 
        \label{property.domination_pss}
        Let $s_1, s_2$ be feasible subpaths with $s_1 \succeq s_2$, both extending a feasible partial path $p$. Suppose $p \oplus s_2$ is feasible. Then $p \oplus s_1$ is feasible too, and $p \oplus s_1 \succeq p \oplus s_2$.
    \end{enumerate}
\end{property}

These domination properties augment Property~\ref{property.domination_benchmark} in our subpath-based label-setting algorithm. Properties~\ref{property.domination_ours}\ref{property.domination_ssa} and \ref{property.domination_ours}\ref{property.domination_pps} generalize Property~\ref{property.domination_benchmark} by propagating domination patterns forward along arc extensions of subpaths and along subpath extensions of paths. Property~\ref{property.domination_ours}\ref{property.domination_pss} is new to our subpath-based label-setting algorithm to propagate domination patterns ``backward'' along subpath extensions of paths. In the context of EVRPTWNL, our main contribution is to develop resources, REFs and domination criteria that satisfy Property~\ref{property.domination_ours} (Theorem~\ref{thm.EVRPTWNL_correctness}, Corollary~\ref{cor.cg_finite_correct}). In the context of the subpath-based label-setting algorithm, our main contribution is to demonstrate that Property~\ref{property.domination_ours} is sufficient to solve the pricing problem to optimality (Theorem~\ref{thm.exact_PP}). In particular, the new backward domination property (Property~\ref{property.domination_ours}\ref{property.domination_pss}) is necessary to establish the exactness of the methodology, and merely extending forward domination concepts (Properties~\ref{property.domination_ours}\ref{property.domination_ssa} and \ref{property.domination_ours}\ref{property.domination_pps}) could lead to arbitrarily large errors in the pricing problem (Proposition~\ref{prop.arbitrarily_large}).

\subsection{First phase: generating non-dominated subpaths}
\label{ssec.labelsetting_subpath}

In the first phase of the algorithm, partial subpaths are iteratively extended along arcs from one non-customer node (charging station or depot) to another. Definition~\ref{def.subpath_resources} specifies the subpath-level resources: the first one tracks the reduced cost; the second one handles charging requirements; the third one handles capacity constraints; and the last three handle time window constraints. Note that our subpath-based decomposition enables us to focus on routing decisions between charging actions, so our first-phase subpath-level resources do not consider charging decisions and the non-linear charging dynamics, which will be treated in the second phase of the algorithm.

\begin{definition}[Subpath resources]
    \label{def.subpath_resources}
    For a partial subpath $s$ with node sequence $(n_0, \dots, n_m)$, we define the following resources:
    \begin{enumerate}[(i),noitemsep]
        \item 
        \label{def.subpath_resources_rc}
        $\cbar{s} = \sum_{\ell=0}^{m-1} \overline{c}_{n_\ell,n_{\ell+1}}$: reduced cost contribution;
        \item 
        \label{def.subpath_resources_B}
        $B(s) = \sum_{\ell=0}^{m-1} b_{n_\ell,n_{\ell+1}}$: charge consumption;
        \item 
        \label{def.subpath_resources_load}
        $D(s) = \sum_{\ell=0}^m d_{n_\ell}$: load served to customers;
        \item 
        \label{def.subpath_resources_Ttravel}
        $\Ttravel{s} = \sum_{\ell=0}^{m-1} t_{n_\ell,n_{\ell+1}}$: the travel time;
        \item 
        \label{def.subpath_resources_TendMin}
        $\TendMin{s}$: earliest end time given time windows; it is defined with the recursive definition:
        \begin{align}
            T^E_0 = \alpha_{n_0}; 
            \quad 
            \forall \ \ell \in \{0, \dots, m-1\}, \
            T^E_{\ell+1} = \max \left\{ 
                T^E_\ell + t_{n_\ell,n_{\ell+1}}, \
                \alpha_{n_{\ell+1}} 
            \right\};
            \quad
            \TendMin{s}=T^E_m.
        \end{align}
        \item 
        \label{def.subpath_resources_TstartMax}
        $\TstartMax{s}$: latest start time given time windows; it is defined with the recursive definition:
        \begin{align}
            T^L_m = \beta_{n_m};
            \quad 
            \forall \ \ell \in \{0, \dots, m-1\}, \
            T^L_\ell = \min \left\{ 
                T^L_{\ell+1} - t_{n_{\ell},n_{\ell+1}}, \
                \beta_{n_\ell}
            \right\};
            \quad
            \TstartMax{s}=T^L_0.
        \end{align}
    \end{enumerate}
\end{definition}

The initial label for $s$ with node sequence $(n_0,)$ is $\left( \cbar{s}, B(s), D(s), \Ttravel{s}, \TendMin{s}, \TstartMax{s} \right) = (0, 0, 0, 0, \alpha_{n_0}, \beta_{n_0})$. Below, Proposition~\ref{prop.subpath_REFs} updates the resources recursively along arc extensions. The main challenge is define a forward update for $\TstartMax{s}$ consistent with backward updates in Definition~\ref{def.subpath_resources} (Proposition~\ref{prop.subpath_REFs}\ref{prop.subpath_REFs_TstartMax}).

\begin{proposition}
    \label{prop.subpath_REFs}
    Let $s$ be a feasible partial subpath and $a = (n, n')$ an arc extension. The following holds:
    \begin{enumerate*}[(i),noitemsep]
        \item 
        \label{prop.subpath_REFs_rc}
        $\cbar{s \oplus a} = \cbar{s} + \overline{c}_{n,n'}$;
        \item 
        \label{prop.subpath_REFs_load}
        $D(s \oplus a) = D(s) + d_{n'}$;
        \item 
        \label{prop.subpath_REFs_B}
        $B(s \oplus a) = B(s) + b_{n,n'}$;
        \item 
        \label{prop.subpath_REFs_Ttravel}
        $\Ttravel{s \oplus a} = \Ttravel{s} + t_{n,n'}$;
        \item 
        \label{prop.subpath_REFs_TendMin}
        $\TendMin{s \oplus a} = \max \left\{ \TendMin{s} + t_{n,n'}, \alpha_{n'} \right\}$; and
        \item 
        \label{prop.subpath_REFs_TstartMax}
        $\TstartMax{s \oplus a} = \min \left\{ \TstartMax{s}, \beta_{n'} - \Ttravel{s \oplus a} \right\}$.
    \end{enumerate*}
    This extension $s \oplus a$ is feasible if:
    \begin{enumerate*}[(a),noitemsep]
        \item 
        \label{prop.subpath_REFs_feasible_load}
        $D(s \oplus a) \leq D$ (vehicle capacity);
        \item 
        \label{prop.subpath_REFs_feasible_B}
        $B(s \oplus a) \leq B$ (battery); and
        \item 
        \label{prop.subpath_REFs_feasible_TW}
        $\TendMin{s \oplus a} \leq \beta_{n'}$ and $\TstartMax{s \oplus a} \geq \alpha_{n_0}$ (time windows).
    \end{enumerate*}
\end{proposition}

Note that each subpath is defined by a sequence of customer visits but not by its departure time. Each subpath can therefore be moved to start or end at different times given by time window requirements---different subpaths induce different sets of feasible start times. In turn, we define subpath domination criteria that tracks reduced costs, vehicle loads, battery capacities but also timing resources. As indicated in the following definition, the correct domination criteria involve non-linear transformations of the timing resources from Definition~\ref{def.subpath_resources}, via the function $Q(\cdot)$.
\begin{definition}[Subpath domination]
    \label{def.subpath_domination}
    Let $s_1$, $s_2$ be feasible partial subpaths starting and ending at the same nodes. If the following hold:
    \begin{enumerate*}[(i),noitemsep]
        \item 
        \label{def.subpath_domination_rc}
        $\cbar{s_1} \leq \cbar{s_2}$;
        \item 
        \label{def.subpath_domination_load}
        $D(s_1) \leq D(s_2)$;
        \item 
        \label{def.subpath_domination_B}
        $B(s_1) \leq B(s_2)$;
        \item 
        \label{def.subpath_domination_TendMin}
        $\TendMin{s_1} \leq \TendMin{s_2}$;
        \item 
        \label{def.subpath_domination_TstartMax}
        $\TstartMax{s_1} \geq \TstartMax{s_2}$; and
        \item 
        \label{def.subpath_domination_Q}
        $Q(s_1) \leq Q(s_2)$, where $Q(s) := \max \left\{ \Ttravel{s}, \TendMin{s} - \TstartMax{s} \right\}$;
    \end{enumerate*}
    then $s_1 \succeq s_2$.
\end{definition}

We interpret criteria \ref{def.subpath_domination_TendMin}--\ref{def.subpath_domination_Q} by letting $f_s(t)$ denote the end time of subpath $s$ if it starts at time $t \geq 0$ (given time windows). The results below show that the full function $f_s(\cdot)$ is encapsulated through the three resources $\Ttravel{s}$, $\TendMin{s}$, and $\TstartMax{s}$, and the domination criteria express that one subpath can end earlier than another for any possible start time. Therefore, although the path-based resources from \cite{desaulniers2016exact} do not apply in our setting, we have derived a compact, low-dimensional representation of time windows in our subpath-based label-setting algorithm that determines the feasibility of each subpath and characterizes its possible end times.
\begin{lemma}
    \label{lemma.subpath_fs}
    Let $s$ be a feasible partial subpath, and $P(s) = \min \left\{ 
            \TstartMax{s}, 
            \TendMin{s} - \Ttravel{s}
        \right\}$. Then:
    \begin{align}
        \label{eq.subpath_fs}
        f_s(t) & = \max \left\{ 
            \TendMin{s},
            Q(s) + t
        \right\} = \begin{cases}
            \TendMin{s}
            & t \in [0, P(s)]
            \\
            Q(s) + t
            & t \in [P(s), \TstartMax{s}]
        \end{cases}
    \end{align}
\end{lemma}

\begin{proposition}
    \label{prop.subpath_fs_domination}
    Definition~\ref{def.subpath_domination}\ref{def.subpath_domination_TendMin}--\ref{def.subpath_domination_Q} is equivalent to:
    \begin{equation}
        \text{dom}(f_{s_2}) \subseteq \text{dom}(f_{s_1}),
        \quad \text{and} \quad 
        \forall t \in \text{dom}(f_{s_2}): 
        \ 
        f_{s_1}(t) \leq f_{s_2}(t).
    \end{equation}
\end{proposition}

These temporal dynamics are shown in Figure~\ref{fig.subpath_fs}. Figure~\ref{fig.subpath_fs_resources} shows that the function $f_s$ is piecewise linear with breakpoint $P(s)$, and is fully characterized by the three resources from Definition~\ref{def.subpath_domination}\ref{def.subpath_domination_TendMin}--\ref{def.subpath_domination_Q}. 
Figure~\ref{fig.subpath_fs_lemma} illustrates the domination patterns elicited in Proposition~\ref{prop.subpath_fs_domination}. 
Figure~\ref{fig.subpath_fs_Q} shows that domination patterns based on the three resources from Definition~\ref{def.subpath_resources}\ref{def.subpath_resources_Ttravel}--\ref{def.subpath_resources_TstartMax} alone are insufficient, and require the non-linear transformation encapsulated in the function $Q$ from Definition~\ref{def.subpath_domination}.

\begin{figure}[h!]
    \centering
    \begin{subfigure}[t]{0.32\textwidth}
    \centering
    \begin{tikzpicture}
        \begin{axis}[
            axis equal image,
            height=8cm,
            xmin=0, xmax=12,
            ymin=0, ymax=18,
            xlabel={\small Start time},
            ylabel={\small Current time},
            axis lines=middle,
            minor tick num=4,
        ]

        \addplot [black, thick, mark=*] coordinates {(0,8) (3,8) (9,14)};
        \node[anchor=north west, inner sep=2pt] at (axis cs:9,14) {$f_{s_1}$};

        \addplot [black, dashed] coordinates {(0,5) (3,8)};
        \node[anchor=north west, inner sep=2pt] at (axis cs:0,7) {\small $\Ttravel{s_1} = 5$};
    
        \addplot [black, dashed] coordinates {(9,0) (9,24)};
        \node[anchor=east, inner sep=2pt] at (axis cs:9,4) {\small $\TstartMax{s_1} = 9$};
    
        \addplot [black, dashed] coordinates {(0,8) (16,8)};
        \node[anchor=south west, inner sep=2pt] at (axis cs:0,8) {\small $\TendMin{s_1} = 8$};
    
        \end{axis}
    \end{tikzpicture}
    \caption{\small
    Illustration of Lemma~\ref{lemma.subpath_fs}.}
    \label{fig.subpath_fs_resources}
    \end{subfigure}\hfill
    \begin{subfigure}[t]{0.32\textwidth}
    \centering
    \begin{tikzpicture}
        \begin{axis}[
            axis equal image,
            height=8cm,
            xmin=0, xmax=12,
            ymin=0, ymax=18,
            xlabel={\small Start time},
            ylabel={\small Current time},
            axis lines=middle,
            minor tick num=4,
        ]
    
        \addplot [black, thick, mark=*] coordinates {(0,12) (5,12)};
        \node[anchor=south west, inner sep=2pt] at (axis cs:5,12) {$f_{s_2}$};

        \addplot [black, thick, mark=*] coordinates {(0,8) (3,8) (9,14)};
        \node[anchor=north west, inner sep=2pt] at (axis cs:9,14) {$f_{s_1}$};

        \end{axis}
    \end{tikzpicture}
    \caption{\small
    Illustration of Proposition~\ref{prop.subpath_fs_domination}.}
    \label{fig.subpath_fs_lemma}
    \end{subfigure}\hfill
    \begin{subfigure}[t]{0.32\textwidth}
    \centering
    \begin{tikzpicture}
        \begin{axis}[
            axis equal image,
            height=8cm,
            xmin=0, xmax=12,
            ymin=0, ymax=18,
            xlabel={\small Start time},
            ylabel={\small Current time},
            axis lines=middle,
            minor tick num=4,
        ]
    
        \addplot [black, thick, mark=*] coordinates {(0,12) (5,12)};
        \node[anchor=south west, inner sep=2pt] at (axis cs:5,12) {$f_{s_2}$};

        \addplot [black, dashed] coordinates {(0,7) (5,12)};
        \node[anchor=west, inner sep=2pt] at (axis cs:0,7) {\small $Q(s_2)$};

        \addplot [black, thick, dotted] coordinates {(5,12) (8.5,12) (10.5,14)};
        \addplot [black, dashed] coordinates {(0,3.5) (8.5,12)};
        \node[anchor=west, inner sep=2pt] at (axis cs:0,3.5) {\small $\Ttravel{s_2}$};

        \addplot [black, thick, mark=*] coordinates {(0,8) (3,8) (9,14)};
        \node[anchor=north west, inner sep=2pt] at (axis cs:9,14) {$f_{s_1}$};
        
        \addplot [black, dashed] coordinates {(0,5) (3,8)};
        \node[anchor=west, inner sep=2pt] at (axis cs:0,5) {\small $Q(s_1) = \Ttravel{s_1}$};

        \end{axis}
    \end{tikzpicture}
    \caption{\small
    Necessity of Definition \ref{def.subpath_domination}\ref{def.subpath_domination_Q}. }
    \label{fig.subpath_fs_Q}
    \end{subfigure}
    \caption{Mapping between $f_s$ and subpath resources. $(\TendMin{s_2}, \TstartMax{s_2}, \Ttravel{s_2}) = (12, 5, 3.5)$ and $(\TendMin{s_1}, \TstartMax{s_1}, \Ttravel{s_1}) = (8, 9, 5)$. (a) Construction of $f_{s_1}$ from $\Ttravel{s_1}$, $\TendMin{s_1}$ and $\TstartMax{s_1}$; (b) Visualization that Definition~\ref{def.subpath_domination}\ref{def.subpath_domination_TendMin}--\ref{def.subpath_domination_Q} is equivalent to $f_{s_2}\geq f_{s_1}$; (c) Necessity of the non-linear transformation $Q$: $f_{s_1} \leq f_{s_2}$ on the domain of $f_{s_2}$ although $\Ttravel{s_1} > \Ttravel{s_2}$.}
    \label{fig.subpath_fs}
\end{figure}

Algorithm~\ref{alg.FindNonDominatedSubpaths} presents the first-phase label-setting procedure to find non-dominated subpaths. Starting at a non-customer node, it iteratively extends partial subpaths along arcs while checking feasibility and pruning dominated partial subpaths, until reaching a non-customer node.

\begin{algorithm} [h!]
\caption{\textsc{FindNonDominatedSubpaths}.}\small
\label{alg.FindNonDominatedSubpaths}
\begin{algorithmic}
\item \textbf{Initialization:} $\Sgen = \varnothing$; store in $\Squeue$ all single-node partial subpaths starting at nodes $i \in \Depots \cup \Stations$.

\begin{itemize}[left=0pt..1em]
    \item[] \textbf{Step 1.} 
    Select $s \in \Squeue$ with smallest end time: $s \in \argmin \Set{ \TendMin{s} | s \in \Squeue}$. 
    Remove $s$ from $\Squeue$, and add it to $\Sgen$. 
    If $s$ is a (complete) subpath, \texttt{GOTO} Step 3. Else, \texttt{GOTO} Step 2.
    \item[] \textbf{Step 2.} For each arc extension $a \in \calA$ of $s$:
    \begin{enumerate}[label*=\arabic*.]
        \item If $s \oplus a$ is not feasible, or dominated by some $s' \in \Squeue \cup \Sgen$, \texttt{CONTINUE}.
        \item Otherwise, remove any $s' \in \Squeue$ such that $s \oplus a \succeq s'$, and add $s \oplus a$ to $\Squeue$.
    \end{enumerate}
    \item[] \textbf{Step 3.} If $\Squeue = \varnothing$, \texttt{STOP}: return $\Sresult = \set{s \in \Sgen|s \text{ is a subpath} }$. Otherwise, \texttt{GOTO} Step 1.
\end{itemize}
\end{algorithmic}
\end{algorithm}

\subsection{Second phase: combining subpaths into paths}
\label{ssec.labelsetting_path}

In the second phase, partial paths are iteratively extended along non-dominated subpaths from one depot to another. The main challenge involves determining timing resources and charge resources that can be easily computed recursively via REFs, and that can be leveraged to prune infeasible and dominated paths via appropriate domination criteria. We formalize the set of path-level resources in the following definition: the first one tracks the reduced cost; the second one handles capacity constraints; the last two characterize the ``best'' ending time and ending charge of a partial path, which requires optimizing the charging sequence for a given sequence of subpaths.

\begin{definition}[Path resources]
    \label{def.path_resources}
    Let $p$ be a partial path with subpath sequence $(s_1, \dots, s_m)$, where each subpath $s_k$ has node sequence $(n_{k,0}, n_{k,1}, \dots, n_{k,m_k})$.
    We define the following resources:
    \begin{enumerate}[(i),noitemsep]
        \item 
        \label{def.path_resources_rc}
        $
            \cbar{p} = \sum_{k=1}^m \sum_{\ell=0}^{m_k-1} \overline{c}_{n_{k,\ell},n_{k,\ell+1}}
        $: reduced cost contribution;
        
        \item 
        \label{def.path_resources_load}
        $
            D(p) = \sum_{k=1}^m \sum_{\ell=0}^{m_k-1} d_{n_{k,\ell}} 
        $: load served to customers;
        \item $\TendMin{p}$: earliest arrival time, defined as the solution of the following problem $\text{CS}^T(p)$;
        \begin{alignat}{3}
            \min \quad
            & \Tend{m}
            \\
            \st \quad
            \label{eq.chargeseq_Tstart_first}
            & \Tstart{1} = 0
            \\
            \label{eq.chargeseq_Tstart_TstartMax}
            & \Tstart{k} \leq \TstartMax{s_k}
            &\quad& \forall \ k \in \{1, \dots, m\}
            \\
            \label{eq.chargeseq_Tend_last}
            & \Tend{m} \leq T
            \\
            \label{eq.chargeseq_Tstart}
            & \Tstart{k+1} = \Tend{k} + \tau_k
            &\quad& \forall \ k \in \{1, \dots, m-1\}
            \\
            \label{eq.chargeseq_Tend}
            & \Tend{k} = \max\left\{\Tstart{k} + \Ttravel{s_k},\TendMin{s_k}\right\}
            &\quad& \forall \ k \in \{1, \dots, m\}
            \\
            \label{eq.chargeseq_Bstart_first}
            & \Bstart{1} = B
            \\
            \label{eq.chargeseq_Bend}
            & \Bend{k} = \Bstart{k} - B(s_k)
            &\quad& \forall \ k \in \{1, \dots, m\}
            \\
            \label{eq.chargeseq_Bstart}
            & \Bstart{k+1} = \gC{\Bend{k}}{\tau_k}
            &\quad& \forall \ k \in \{1, \dots, m-1\}
            \\
            \label{eq.chargeseq_nonneg}
            & \Tstart{k}, \Tend{k}, 
            \Bstart{k}, \Bend{k} \geq 0
            &\quad& \forall \ k \in \{1, \dots, m\}
            \\
            \label{eq.chargeseq_tau_nonneg}
            & \tau_k \geq 0
            &\quad& \forall \ k \in \{1, \dots, m-1\}
        \end{alignat}
        \item $\BendMax{p}$: maximum ending charge at the minimum ending time $\TendMin{p}$, defined as the solution of the following problem $\text{CS}^B(p)$:
        \begin{equation}
            \label{eq.chargeseq_problem_Bend}
            \max \quad B^{\text{end}}_m;
            \quad
            \st \quad
            \text{Equations~\eqref{eq.chargeseq_Tstart_first}--\eqref{eq.chargeseq_tau_nonneg}};
            \quad\text{and}\quad
            \Tend{m} = \TendMin{p}
        \end{equation}
    \end{enumerate}
\end{definition}

This definition involves a non-linear optimization problem---where the non-linearities stem from charging dynamics---to determine the feasibility of a partial path and the corresponding resources. We prove in Proposition~\ref{prop.chargeseq} that we can characterize the optimal charging sequence for each subpath sequence via a greedy procedure. This result exploits the concavity of the charging schedule $g(\cdot)$, which guarantees that it is at least as fast to gain some amount of charge when starting with a lower charge, i.e.: if $b_1 \leq b_2$, then $\gT{b_1}{b_1 + b} \leq \gT{b_2}{b_2 + b}$ (Proposition 2.1 of \cite{lee_exact_2021}).

\begin{proposition}
    \label{prop.chargeseq}
    Let $p$ be a feasible partial path with subpath sequence $(s_1, \dots, s_m)$. The following recursive updates produce a feasible solution and return optimal values $\Tend{k}$ and $\Bend{k}$ for $\text{CS}^T(s_1, \dots, s_k)$ and $\text{CS}^B(s_1, \dots, s_k)$ respectively, for all $k\in\{1, \dots, m\}$.
    
    Initialize $\Tend{1} = f_{s_1}(0)$ and $\Bend{1} = B - B(s_1)$. Repeat, for $k \in \{1, \dots, m-1\}$:
    \begin{align}
    \label{eq.chargeseq_tau_tilde_alg}
    \widetilde{\tau}_k & \leftarrow \max \left\{ 0, \gT{ \Bend{k} }{ B(s_{k+1}) } \right\}
    \\
    \label{eq.chargeseq_tau_alg}
    \tau_k & \leftarrow  \max \left\{ 
        \widetilde{\tau}_k, \
        P(s_{k+1}) - \Tend{k}
    \right\}
    \\
    \label{eq.chargeseq_Tstart_alg}
    \Tstart{k+1} & \leftarrow  \Tend{k} + \tau_k
    \\
    \label{eq.chargeseq_Bstart_alg}
    \Bstart{k+1} & \leftarrow  \gC{\Bend{k}}{\tau_k}
    \\
    \label{eq.chargeseq_Tend_alg}
    \Tend{k+1} & \leftarrow f_{s_{k+1}}(\Tstart{k+1}) = \max \{ \Tstart{k+1} + \Ttravel{s_{k+1}}, \TendMin{s_{k+1}} \}
    \\
    \label{eq.chargeseq_Bend_alg}
    \Bend{k+1} & \leftarrow \Bstart{k+1} - B(s_{k+1})
\end{align}
\end{proposition}

Specifically, for any subpath sequence $(s_1, \dots, s_m)$, Proposition~\ref{prop.chargeseq} computes the minimum charging time $\widetilde{\tau}_k$ at the end of each subpath $s_k$ to ensure the feasibility of the subsequent subpath $s_{k+1}$. In the absence of time windows, it would be without loss of generality to always charge to the minimum required $\widetilde{\tau}_k$ (Figure~\ref{fig.extension_noTW}). However, with time windows, it might be beneficial to charge for longer than necessary and avoid wasteful waiting along subpath $s_{k+1}$, and the corresponding optimal charging time is denoted by $\tau_k\geq\widetilde{\tau}_k$ (Figure~\ref{fig.extension_TW}). These charging times can then be plugged into the definition of the timing and charge resources (Definition~\ref{def.path_resources}), as guaranteed by Proposition~\ref{prop.chargeseq}.

\begin{figure}[h!]
    \centering
    \begin{subfigure}[t]{0.48\textwidth}
    \centering
    \begin{tikzpicture}
        \begin{axis}[
            axis equal image,
            width=9cm,
            xmin=0, xmax=28,
            ymin=0, ymax=16,
            xlabel={\small Time},
            ylabel={\small Charge},
            axis lines=middle,
            xtick=\empty,
            ytick=\empty,
        ]

        \addplot [black, thick, mark=*] coordinates {(0,14) (8,5)};
        \node[anchor=south west, inner sep=2pt] at (axis cs:4,9.5) {$s_1$};
        \node[anchor=north, inner sep=2pt] at (axis cs:8,5) {\small $(\TendMin{p_1}, \BendMin{p_1})$};

        \addplot [purple, thick, domain=8:10] {-(x-14)^2/4 + 14};
        \addplot [purple, thick, dotted, domain=10:14] {-(x-14)^2/4 + 14};

        \draw[latex-latex, purple] (8,1)
                                    node[left] {\small \color{purple} $\widetilde{\tau}$} 
                                    --(10,1);
        
        \addplot [blue!50!black, thick, mark=*] coordinates {(10,10) (18,0)};
        \node[blue!50!black, anchor=south west, inner sep=2pt] at (axis cs:14,5) {$s_2$};
        \node[blue!50!black, anchor=south, inner sep=2pt] at (axis cs:18,0) {\small $(\TendMin{p_2}, \BendMin{p_2})$};
        
        \end{axis}
    \end{tikzpicture}
    \caption{\small Without time windows: charging to minimum $\widetilde{\tau}$}
    \label{fig.extension_noTW}
    \end{subfigure}%
    \hfill
    \begin{subfigure}[t]{0.48\textwidth}
    \centering
    \begin{tikzpicture}
        \begin{axis}[
            axis equal image,
            width=9cm,
            xmin=0, xmax=28,
            ymin=0, ymax=16,
            xlabel={\small Time},
            ylabel={\small Charge},
            axis lines=middle,
            xtick=\empty,
            ytick=\empty,
        ]

        \addplot [black, thick, mark=*] coordinates {(0,14) (8,5)};
        \node[anchor=south west, inner sep=2pt] at (axis cs:4,9.5) {$s_1$};
        \node[anchor=north, inner sep=2pt] at (axis cs:8,5) {\small $(\TendMin{p_1}, \BendMin{p_1})$};

        \addplot [purple, thick, domain=8:12] {-(x-14)^2/4 + 14};
        \addplot [purple, thick, dotted, domain=12:14] {-(x-14)^2/4 + 14};
        \draw[latex-latex, purple] (8,1)
                                    node[left] {\small \color{purple} $\widetilde{\tau}$} 
                                    --(10,1);
        \draw[latex-latex, purple] (8,2.5)
                                    node[left] {\small \color{purple} ${\tau}$}
                                    --(12,2.5);

        \addplot [blue!50!black, dashed] coordinates {(12,0) (12,16)};
        \node[blue!50!black, anchor=west, inner sep=2pt] at (axis cs:12,15) {\small $P(s_2)$};
        
        \addplot [blue!50!black, thick, mark=*] coordinates {(12,13) (20,3)};
        \node[blue!50!black, anchor=south west, inner sep=2pt] at (axis cs:16,8) {$s_2$};
        \node[blue!50!black, anchor=north, inner sep=2pt] at (axis cs:20,3) {\small {$(\TendMin{p_2}, \BendMin{p_2})$}};

        \addplot [blue!50!black, thick, dotted] coordinates {(10,10) (12,10) (20,0)};
        
        \end{axis}
    \end{tikzpicture}
    \caption{\small With time windows: charging to $\tau$ while waiting}
    \label{fig.extension_TW}
    \end{subfigure}
    \caption{Extension of $p = \{s_1\}$ along subpath $s_2$. Without time windows, charging to minimum, $\widetilde{\tau}$, is optimal due to concavity of charging schedule. With time windows, waiting time might be forced on $s_2$ and it would be optimal to charge for $\tau > \widetilde{\tau}$ until time $P(s_2) = \min \{ \TendMin{s_2} - \Ttravel{s_2}, \TstartMax{s_2}\}$.}
    \label{fig.extension_EVRPTWNL}
\end{figure}

Proposition~\ref{prop.path_REFs} leverages these insights to define recursive updates of the path-based resources along subpath extensions. This follows from Proposition~\ref{prop.chargeseq}, which shows that the definition of $\tau$ recovers the ``best'' charging time between partial path $p$ and subpath $s$ to compute $\TendMin{p \oplus s}$ and $\BendMax{p \oplus s}$ as a function of $\TendMin{p}$ and $\BendMax{p}$, regardless of the future subpaths.

\begin{proposition}
    \label{prop.path_REFs}
    Let $p$ be a feasible partial path and $s$ a subpath extending $p$. 
    Define $\widetilde{\tau} = \max \left\{ 0, \gT{ \BendMax{p} }{ B(s) } \right\}$ and $\tau = \max \{ \widetilde{\tau}, P(s) - \TendMin{p}\}$.
    Then:
    \begin{enumerate*}[(i),noitemsep]
        \item 
        \label{prop.path_REFs_rc}
        $\cbar{p \oplus s} = \cbar{p} + \cbar{s}$;
        \item 
        \label{prop.path_REFs_load}
        $D(p \oplus s) = D(p) + D(s)$;
        \item 
        \label{prop.path_REFs_TendMin}
        $\TendMin{p \oplus s} = f_s \left( \TendMin{p} + \tau \right)$; and
        \item 
        \label{prop.path_REFs_BendMax}
        $\BendMax{p \oplus s} = \gC{ \BendMax{p} }{ \tau } - B(s)$.
    \end{enumerate*}
    This extension is feasible if:
    \begin{enumerate*}[(a),noitemsep]
        \item 
        \label{prop.path_REFs_feasible_load}
        $D(p \oplus s) \leq D$ (vehicle capacity);
        \item
        \label{prop.path_REFs_feasible_TstartMax}
        $\TendMin{p} + \tau \leq \TstartMax{s}$ (battery); and
        \item 
        \label{prop.path_REFs_feasible_T}
        $\TendMin{p \oplus s} \leq T$ (time horizon).
    \end{enumerate*}
\end{proposition}

This allows us to define path-level domination criteria, again using non-linear transformations of the resources from Definition~\ref{def.path_resources}. Interestingly, criterion \ref{def.path_domination_BendMax} allows the maximum charge of $p_1$ to be smaller than that of $p_2$, as long as there is sufficient time to charge to $\BendMax{p_2}$ (Figure~\ref{fig.domination_EVRPTWNL}).
\begin{definition}[Path domination]
    \label{def.path_domination}
    Let $p_1$, $p_2$ be feasible partial paths starting and ending at the same nodes. If the following hold:
    \begin{enumerate*}[(i),noitemsep]
        \item 
        \label{def.path_domination_rc}
        $\cbar{p_1} \leq \cbar{p_2}$;
        \item 
        \label{def.path_domination_load}
        $D(p_1) \leq D(p_2)$; and
        \item 
        \label{def.path_domination_BendMax}
        $\TendMin{p_1} +  \gT{ \BendMax{p_1} }{ \max \left\{ \BendMax{p_1}, \BendMax{p_2} \right\} } \leq \TendMin{p_2}$;
    \end{enumerate*}
    then $p_1 \succeq p_2$.
\end{definition}

\begin{figure}[h!]
    \centering
    \begin{tikzpicture}
        \begin{axis}[
            axis equal image,
            width=9cm,
            xmin=0, xmax=28,
            ymin=0, ymax=16,
            xlabel={\small Time},
            ylabel={\small Charge},
            axis lines=middle,
            xtick=\empty,
            ytick=\empty,
        ]

        \addplot [black, thick, mark=*] coordinates {(13,5)};
        \addplot [black, dashed] coordinates {(13,0) (13,5)};
        \addplot [black, dashed] coordinates {(0,5) (28,5)};
        \node[anchor=south east, inner sep=2pt] at (axis cs:13,5) {\small $(\TendMin{p_1}, \BendMin{p_1})$};

        \addplot [name path=curve, purple, thick, domain=13:19] {-(x-19)^2/4 + 14};
        \addplot [name path=bottom, draw=none, domain=13:25] {0};
        \addplot [name path=top, draw=none, domain=19:25] {14};

        \addplot [fill=red!50, opacity=0.2] fill between[of=curve and bottom, soft clip={domain=13:19}];
        \addplot [fill=red!50, opacity=0.2] fill between[of=top and bottom, soft clip={domain=19:25}];
        \end{axis}
    \end{tikzpicture}
    \caption{Domination criteria (Definition~\ref{def.path_domination}\ref{def.path_domination_BendMax}): if $p_1\succeq p_2$, then $(\TendMin{p_2}, \BendMax{p_2})$ falls in the shaded area.}
    \label{fig.domination_EVRPTWNL}
\end{figure}

Algorithm~\ref{alg.FindNonDominatedPaths} presents the second-phase label-setting procedure to find non-dominated paths. Starting at a depot, it iteratively extends partial paths along subpaths while checking feasibility and pruning dominated partial paths, until reaching a depot. It is analogous to Algorithm~\ref{alg.FindNonDominatedSubpaths}, except that it takes place on a multigraph due to multiple non-dominated subpaths extending the same partial path $p$ to the next node $n' \in \Stations \cup \Depots$, each covering different sets of customers.

\begin{algorithm} [h!]
\caption{\textsc{FindNonDominatedPaths}$(\calS^*)$}\small
\label{alg.FindNonDominatedPaths}
\begin{algorithmic}
\item \textbf{Initialization:} $\Pgen = \varnothing$; store in $\Pqueue$ all single-node partial paths starting in $i \in \Depots$.

\begin{itemize}[left=0pt..1em]
    \item[] \textbf{Step 1.} 
    Select $p \in \argmin \Set{ \TendMin{p} | p \in \Pqueue}$. 
    Remove $p$ from $\Pqueue$, and add it to $\Pgen$. 
    If $p$ is a (complete) path, \texttt{GOTO} Step 3. Else, \texttt{GOTO} Step 2.
    \item[] \textbf{Step 2.} 
    For each subpath $s \in \calS^*$ extending $p$:
    \begin{enumerate}[label*=\arabic*.]
        \item If $p \oplus s$ is not feasible, 
        or dominated by some $p' \in \Pqueue \cup \Pgen$,
        \texttt{CONTINUE}.
        \item Otherwise, 
        remove any $p' \in \Pqueue$ such that $p \oplus s \succeq p'$, 
        and add $p \oplus s$ to $\Pqueue$.
    \end{enumerate}
    \item[] \textbf{Step 3.} 
    If $\Pqueue = \varnothing$, \texttt{STOP}: return $\Presult = \Set{p \in \Pgen | p \text{ is a path} }$. Otherwise, \texttt{GOTO} Step 1.
\end{itemize}
\end{algorithmic}
\end{algorithm}

\subsection{Correctness of the subpath-based label-setting algorithm for the EVRPTWNL}
\label{ssec.labelsetting_correctness}

We now prove the correctness of the subpath-based label-setting algorothm, both in general terms and in the context of the EVRPTWNL. We first show the critical result in the context of the EVRPTWNL, that is, that the resources, REFs and dominance criteria defined in Section~\ref{ssec.labelsetting_subpath} at the subpath level and in Section~\ref{ssec.labelsetting_path} at the path level satisfy Property~\ref{property.domination_ours}.
\begin{theorem}
    \label{thm.EVRPTWNL_correctness}
    The resources (Definitions
    \ref{def.subpath_resources} and~\ref{def.path_domination}), REFs (Propositions
    \ref{prop.subpath_REFs} and~\ref{prop.path_REFs}) and domination criteria (Definitions~\ref{def.subpath_domination} and~\ref{def.path_domination})
    satisfy Property~\ref{property.domination_ours} for the EVRPTWNL.
\end{theorem}

Next, Theorem~\ref{thm.exact_PP} provides the main result on our subpath-based label-setting algorithm, by showing that Property~\ref{property.domination_ours} is sufficient to establish its exactness in the pricing problem. This result makes use of two technical assumptions stating that one resource tracks the reduced cost and one resource is nonnegative and increasing to guarantee finite termination. In practice, these two conditions are immediately satisfied by including reduced cost and arrival time in the resources.

\begin{assumption}\label{ass.technical}
    For both subpaths and paths, one resource must sum to the reduced cost, and one resource must be nonnegative, bounded above, and strictly increasing along extensions.
\end{assumption}

\begin{theorem}\label{thm.exact_PP}
    Under Property~\ref{property.domination_ours} and Assumption~\ref{ass.technical}, 
    Algorithms~\ref{alg.FindNonDominatedSubpaths} and~\ref{alg.FindNonDominatedPaths} 
    together find paths of negative reduced cost in the pricing problem of the EVRPTWNL, or proves that none exists.
\end{theorem}

In fact, this result hinges on Property~\ref{property.domination_ours}\ref{property.domination_pss} to propagate domination from subpaths ($s_1 \succeq s_2$) to paths ($p \oplus s_1 \succeq p \oplus s_2$). In fact, if we defined resources and REFs that merely satisfied Properties~\ref{property.domination_ours}\ref{property.domination_ssa}--\ref{property.domination_pps}, then our subpath-based label-setting algorithm would no longer be exact and could lead to arbitrarily large errors in the pricing problem. This is proved in Proposition~\ref{prop.arbitrarily_large}.

\begin{proposition}
    \label{prop.arbitrarily_large}
    Under Properties~\ref{property.domination_ours}\ref{property.domination_ssa}--\ref{property.domination_pps}, 
    Algorithms~\ref{alg.FindNonDominatedSubpaths} and~\ref{alg.FindNonDominatedPaths} 
    can fail to find paths of arbitrarily large negative reduced cost in the pricing problem of the EVRPTWNL. 
\end{proposition}

Note, also, that this property is only needed for ``complex'' resources and REFs in the pricing problem. In the case where subpath-level and path-level resources satisfy monotonicity conditions, Properties~\ref{property.domination_ours}\ref{property.domination_ssa}--\ref{property.domination_pps} are sufficient and Property~\ref{property.domination_ours}\ref{property.domination_pss} is automatically satisfied (Remark~\ref{remark.property_i_ii_sufficient}). This would be the case for instance for the capacitated EVRP, as reduced costs, travel times and vehicle loads are additive and monotonic along each extension. However, Properties~\ref{property.domination_ours}\ref{property.domination_ssa}--\ref{property.domination_pps} are no longer sufficient to handle time windows---as well as \textit{ng}-sets, as discussed in the next section.

\begin{remark}
    \label{remark.property_i_ii_sufficient}
    Property~\ref{property.domination_ours}\ref{property.domination_pss} is automatically satisfied under monotonicity assumptions on subpath-level and path-level resources.
\end{remark}

Altogether, our domination framework in Property~\ref{property.domination_ours} establishes the exactness of the subpath-based label-setting algorithm for hierarchical problems combining discrete routing decisions and continuous resource management decisions. Guided by Property~\ref{property.domination_ours}, Theorems~\ref{thm.EVRPTWNL_correctness} and~\ref{thm.exact_PP} imply the correctness of our subpath-based label-setting algorithm in the pricing problem. Given the finite number of paths (Definition~\ref{def.path}), this implies the correctness of the resulting column generation algorithm to solve the EVRPTWNL linear relaxation. This is stated in the corollary below.

\begin{corollary}
    \label{cor.cg_finite_correct}
    Algorithm~\ref{alg.CG} converges finitely to an optimal solution of $\EVRPLP(\calP)$, when Step 2. is solved via the subpath-based label-setting algorithm 
    (Algorithms~\ref{alg.FindNonDominatedSubpaths} and~\ref{alg.FindNonDominatedPaths}) 
    with the corresponding resources 
    (Definitions
    \ref{def.subpath_resources},~\ref{def.path_domination}), 
    REFs (Propositions
    \ref{prop.subpath_REFs},~\ref{prop.path_REFs}) 
    and domination criteria (Definitions~\ref{def.subpath_domination},~\ref{def.path_domination}).
\end{corollary}
\section{Tighter relaxations via adaptive ng-relaxations and cutting planes}
\label{sec.tighter}

We augment the algorithm with adaptive \textit{ng}-relaxations and subset-row cuts (Algorithm~\ref{alg.AdaptiveColumnGenerationWithCuts}) to tighten the linear relaxation. For both extensions, we develop dedicated domination criteria in our two-phase label-setting algorithm and prove the convergence of the resulting column generation schemes---again, by showing that Property~\ref{property.domination_ours} is satisfied. We also add a restricted master heuristic (RMH) to retrieve a feasible integer solution. The full procedure is detailed in Algorithm~\ref{alg.AdaptiveColumnGenerationWithCuts}.

\begin{algorithm}
\caption{Augmented column generation with adaptive \textit{ng}-relaxations and subset-row cuts.}\small
\label{alg.AdaptiveColumnGenerationWithCuts}
\begin{algorithmic}
\item \textbf{Initialization:} 
\textit{ng}-neighborhood $\calN^0$; 
set of cuts $\calQ = \varnothing$; 
set of paths $\Pinit^0 \subseteq \calP(\calN^0)$; 
$t = 0$.

\item Iterate between Steps 1-4.

\begin{itemize}[left=0pt..1em]
    \item[] \textbf{Step 1: Column generation.} 
    Solve $\EVRPLP(\calP(\calN^{t}))$ via 
    \textsc{ColumnGeneration}$(\calP(\calN^{t}))$ with cuts $\calQ$, 
    starting with $\Pinit^{t}$. Obtain solution $\bz^t$, set of \textit{ng}-feasible paths $\calP^t \subseteq \calP(\calN^t)$ with respect to $\calN^t$, and objective $\OPTLP(\calP(\calN^{t})) = \OPTLP(\calP^t)$.
    \item[] \textbf{Step 2: Solve RMH.} 
    Solve $\EVRP(\calP^t)$ with integrality. Obtain solution $\bz^t_{I}$, objective $\OPT(\calP^t)$.
    \item[] \textbf{Step 3: Termination.} If $\bz^t$ uses elementary paths and no subset-row cut is violated, \texttt{STOP}; return dual and primal solutions and bounds $(\bz^t, \bz^t_{I}, \OPTLP(\calP^t), \OPT(\calP^t))$ respectively.
    \item[] \textbf{Step 4: Elementarity.} If the support of $\bz^t$ uses non-elementary paths: 
    for each such path $p$,
    and for all cycles $\{i, n_1, \cdots, n_m, i\}$ (with $i\in\Custs$) in its node sequence, define $\calN^{t+1}$ by adding $i$ to the subsets $N_{n_1}, \cdots, N_{n_m}$. 
    Define $\Pinit^{t+1} = \calP^t \cap \calP(\calN^{t+1})$, increment $t\gets t+1$, and \texttt{GOTO} Step 1. Else \texttt{GOTO} Step 5.
    \item[] \textbf{Step 5: Integrality.} 
    Find $S \subseteq \Custs$ and $\Set{ w_i | i \in S}$ 
    such that $\bz^t$ violates Equation~\eqref{eq.SRC} over $\calP^t$. 
    Add $(S, \bw)$ to $\calQ$, 
    define $\Pinit^{t+1} = \calP^t$, 
    increment $t\gets t+1$, and \texttt{GOTO} Step 1.
\end{itemize}
\end{algorithmic}
\end{algorithm}

\subsection{Adaptive \textit{ng}-relaxations for elementarity constraints}
\label{ssec.ng}

The linear relaxation $\EVRPLP(\calP)$ depends on the path set $\calP$. Considering the full set of paths $\Pnone$ would lead to a weaker linear relaxation, with short cycles. Vice versa, considering the set of elementary paths $\Pelem$ would tighten the relaxation but would require one resource per customer \citep{beasley1989algorithm}. In between, solving $\EVRPLP(\calP)$ over a restricted set of paths 
$\Pnone \subseteq \calP \subseteq \Pelem$ tightens the relaxation without imposing full elementarity.

\begin{definition}[\textit{ng}-neighborhood]
    An \textit{ng}-neighborhood is a collection of subsets $\calN = \Set{ N_i \subseteq \calV | i \in \calV} $ where:
    \begin{enumerate*}[(i)]
        \item $i \in N_i, \ \forall i \in \calV$;
        \item $N_i \subseteq \Custs, \forall i \in \Custs$; and
        \item $N_i \subseteq \Custs \cup \{ i \}, \forall i \in \Depots \cup \Stations$.
    \end{enumerate*}
\end{definition}

\begin{definition}[\textit{ng}-feasibility]
    A (partial) path or subpath is \textit{ng}-feasible if its node sequence $(n_0, \dots, n_m)$ is \textit{ng}-feasible, i.e.,: for every $j < k$ with $n_j = n_k$, there exists $j < \ell < k$ with $n_j \notin N_{n_\ell}$. Let $\calP(\calN)$ store the \textit{ng}-feasible paths with respect to $\calN$.
\end{definition}

Intuitively, \textit{ng}-feasible paths are ``locally elementary'', in that customer $i$ can only be visited multiple times if a node whose \textit{ng}-neighborhood does not contain $i$ is visited in between. Larger \textit{ng}-neighborhoods therefore make paths with small cycles \textit{ng}-infeasible, thus tightening the $\EVRPLP(\calP)$ relaxation \citep{baldacci2011new}. We adopt the adaptive \textit{ng}-relaxation approach from \cite{martinelli2014efficient}, which alternates between solving $\EVRPLP(\calP(\calN))$ and expanding $\calN$ to eliminate non-elementary paths (Steps 1--4 of Algorithm~\ref{alg.AdaptiveColumnGenerationWithCuts}). This solves $\EVRPLP(\Pelem)$ without imposing elementarity. The main challenge involves defining resources, REFs and domination criteria to incorporate the \textit{ng}-relaxations into our two-phase label-setting algorithm. We address this question below.

\subsubsection*{\textit{ng}-relaxations in our two-phase label-setting algorithm.}

In path-based label-setting algorithms, \textit{ng}-feasible paths are computed by tracking the forward \textit{ng}-set $\Pi(p)$, defined as the set of nodes that cannot be appended to a path $p$ while retaining \textit{ng}-feasibility. Specifically, a partial path $p$ can be extended along arc $a = (\enode{p}, \nnode)$ if $\nnode \notin \Pi(p)$. However, $\Pi(p)$ is no longer sufficient to propagate domination patterns in our subpath-based label-setting algorithm.

We augment our algorithm with three resource definitions, formalized in Definition~\ref{def.ng_resources}: (i) forward \textit{ng}-set $\Pi$, (ii) backward \textit{ng}-set $\Pi^{-1}$, and (iii) \textit{ng}-residue $\Omega$. Whereas the forward and backward \textit{ng}-sets were initially introduced in the context of bidirectional path-based label-setting, we prove in this paper that they are necessary for the correctness of our (uni-directional) two-phase label-setting algorithm, together with the third resource termed \textit{ng}-residue.

\begin{definition}[\textit{ng}-route resources]
    \label{def.ng_resources}
    For a node sequence $N = (n_0, \dots, n_m)$, we define:
    \begin{enumerate}[(i),noitemsep]
        \item 
        \label{def.forward_ngset}
        its forward \textit{ng}-set: 
        $\Pi(N) := \Set{n_r | n_r \in \bigcap_{\rho=r+1}^m N_{n_\rho}, \ r \in \{0, \dots, m-1\}} \cup \{n_m\}$;
        \item 
        \label{def.backward_ngset}
        its backwards \textit{ng}-set: 
        $\Pi^{-1}(N) := \{ n_0 \} \cup \Set{n_r | n_r \in \bigcap_{\rho=0}^{r-1} N_{n_\rho}, \ r \in \{1, \dots, m\}}$; and
        \item 
        \label{def.ng_residue}
        its \textit{ng}-residue: 
        $\Omega(N) := \bigcap_{\rho=0}^m N_{n_\rho}$.
    \end{enumerate}
\end{definition}

All three resources apply to subpaths in the first phase of our algorithm, while only the forward \textit{ng}-set applies to paths in the second phase. Specifically, forward \textit{ng}-sets extend domination forward for subpath extensions or path extensions (Properties~\ref{property.domination_ours}\ref{property.domination_ssa} and \ref{property.domination_pps}) but backward \textit{ng}-sets are needed to extend domination backward in our second phase (Property~\ref{property.domination_ours}\ref{property.domination_pss}). Finally, the \textit{ng}-residue $\Omega(\cdot)$ is required to update $\Pi(p \oplus s)$ recursively from $\Pi(s)$ along subpath extensions.

\begin{proposition}
    \label{prop.subpath_REFs_ngroute}
    Let $s$ be a feasible partial subpath and $a = (\enode{s}, \nnode)$
    be an arc extending $s$. Then:
    \begin{enumerate*}[(i),noitemsep]
        \item $\Pi(s \oplus a) = \big( \Pi(s) \cap N_{\nnode} \big) \cup \{ \nnode \}$;
        \item $\Pi^{-1}(s \oplus a) = \Pi^{-1}(s) \cup \big( \{ \nnode \} \cap \Omega(s) \big)$;
        \item $\Omega(s \oplus a) = \Omega(s) \cap N_{\nnode}$.
    \end{enumerate*}
    This extension is \textit{ng}-feasible if $\nnode \notin \Pi(s)$.
\end{proposition}

\begin{definition}[Subpath \textit{ng}-domination]
    \label{def.subpath_domination_ngroute}
    Let $s_1, s_2$ be partial subpaths between the same nodes. $s_1$ dominates $s_2$ if, additionally,
    \begin{enumerate*}[(i),noitemsep]
        \item $\Pi(s_1) \subseteq \Pi(s_2)$;
        \item $\Pi^{-1}(s_1) \subseteq \Pi^{-1}(s_2)$;
        \item $\Omega(s_1) \subseteq \Omega(s_2)$.
    \end{enumerate*}
\end{definition}

\begin{proposition}
    \label{prop.path_REFs_ngroute}
    Let $p$ be a feasible partial path and $s$ be a subpath extending $p$. Then: $\Pi(p \oplus s) = \big( \Pi(p) \cap \Omega(s) \big) \cup \Pi(s)$.
    This extension is \textit{ng}-feasible if $\Pi(p) \cap \Pi^{-1}(s) \subseteq \{ \snode{s} \}$.
\end{proposition}

\begin{definition}[Path \textit{ng}-domination]
    \label{def.path_domination_ngroute}
    Let $p_1, p_2$ be feasible partial paths starting and ending at the same nodes. $p_1$ dominates $p_2$ if, additionally, $\Pi(p_1) \subseteq \Pi(p_2)$.
\end{definition}

These characterizations enable to verify that Property~\ref{property.domination_ours} is satisfied, and to thereby establish the exactness of the two-phase label-setting algorithm and the resulting column generation scheme.

\begin{proposition}
    \label{prop.ngroute_correctness}
    The definitions and propositions in Theorem~\ref{thm.EVRPTWNL_correctness}, 
    together with 
    Definitions
    \ref{def.subpath_domination_ngroute},
    \ref{def.path_domination_ngroute},
    and Propositions
    \ref{prop.subpath_REFs_ngroute},
    \ref{prop.path_REFs_ngroute},
    satisfy Property~\ref{property.domination_ours} for the
    EVRPTWNL 
    with \textit{ng}-relaxation $\calN$.
\end{proposition}

\tikzset{
    inner sep=2pt,
    sourcesink/.style={rectangle,draw=black,fill=purple!50!white,thick,minimum size=4mm},
    charger/.style={rectangle,draw=black,fill=yellow,thick,minimum size=4mm},
    customer/.style={circle,draw=black,thick,minimum size=4mm},
    >={Stealth},
    legend/.style={right,inner sep=2mm,text height=1.5ex,text depth=0.25ex},
    lasso/.style={dashed,thick,rounded corners,inner sep=3mm},
}
\begin{figure}[h!]
    \centering
    \begin{subfigure}{0.48\textwidth}
        \centering
        \begin{tikzpicture}
            \node[sourcesink]   (r0)    at ( 1.5,-0.5) {\small $a$};
            \node[customer]     (c2)    at ( 2.0, 1.0) {\small $1$};
            \node[customer]     (c3)    at ( 3.0, 0.0) {\small $2$};
            \node[charger]      (r1)    at ( 4.0, 1.0) {\small $b$};
            \node[customer]     (c4)    at ( 5.0, 0.0) {\small $5$};
            \node[customer]     (c5)    at ( 5.5, 1.5) {\small $3$};
            \node[customer]     (c6)    at ( 6.5, 0.5) {\small $4$};
            \node[charger]      (r2)    at ( 6.5,-1.0) {\small $c$};
            \node[customer]     (c7)    at ( 8.0,-1.0) {\small $6$};
            \draw[-]    (r0) to [bend left=15] (c2);
            \draw[-]    (c2) to [bend left=45] node[above,inner sep=2mm] {$p$} (c5);
            \draw[-]    (c5) to [bend left=45] (c4);
            \draw[->]   (c4) to [bend left=45] (r1);
            \draw[-]    (r1) to [bend left=15] (c5);
            \draw[-]    (c5) to [bend left=30] node[above,inner sep=2mm,anchor=south west,yshift=12pt,xshift=-24pt] {$s: \Pi^{-1}(s) = \{ b, 3 \}$} (c6);
            \draw[->]   (c6) to [bend left=15] (r2);
            \node       (t) at ($(c4) + (0,-0.5)$) {\small $\Pi(p) = \{b, 5, 3\}$};
        \end{tikzpicture}%
        \caption{Example where $p \oplus s$ is \textit{ng}-infeasible.}
    \end{subfigure}%
    \hfill
    \begin{subfigure}{0.48\textwidth}
        \centering
        \begin{tikzpicture}
            \node[sourcesink]   (r0)    at ( 1.5,-0.5) {\small $a$};
            \node[customer]     (c2)    at ( 2.0, 1.0) {\small $1$};
            \node[customer]     (c3)    at ( 3.0, 0.0) {\small $2$};
            \node[charger]      (r1)    at ( 4.0, 1.0) {\small $b$};
            \node[customer]     (c4)    at ( 5.0, 0.0) {\small $5$};
            \node[customer]     (c5)    at ( 5.5, 1.5) {\small $3$};
            \node[customer]     (c6)    at ( 6.5, 0.5) {\small $4$};
            \node[charger]      (r2)    at ( 6.5,-1.0) {\small $c$};
            \node[customer]     (c7)    at ( 8.0,-1.0) {\small $6$};
            \draw[-]    (r0) to [bend left=15] (c3);
            \draw[-]    (c3) to [out=0,in=225] node[above,inner sep=2mm] {\small $p$} (c4);
            \draw[->]   (c4) to [out=75,in=0] (r1);
            \draw[-]    (r1) to [bend left=15] (c5);
            \draw[-]    (c5) to [bend left=30] node[above,inner sep=2mm,anchor=south west,yshift=12pt,xshift=-24pt] {$s: \Pi^{-1}(s) = \{ b, 3 \}$} (c6);
            \draw[->]   (c6) to [bend left=15] (r2);
            \node       (t) at ($(c4) + (0,-0.6)$) {\small $\Pi(p) = \{b, 5, 2\}$};
        \end{tikzpicture}%
        \caption{Example where $p \oplus s$ is \textit{ng}-feasible.}
    \end{subfigure}%
    \caption{Illustration of \textit{ng}-feasibility along subpath extensions of partial paths. Numbered and lettered nodes denote customers and non-customers; $N_b=\{b,2,3,5\}$, $N_3=\{3,4,5\}$, $N_4=\{3,4,5\}$, $N_5=\{2,3,5\}$.}
    \label{fig.ngfeas}
    \vspace{-12pt}
\end{figure}

In summary, although our two-phase label-setting algorithm is uni-directional, it requires domination criteria based on forward and backward \textit{ng}-sets to guarantee \textit{ng}-feasibility. This stems from the fact that multiple non-dominated subpaths can extend subpath sequences between the same pair of nodes in our second-level procedure. Computationally, since $\Pi(s) \subseteq N_i$, $\Pi^{-1}(s) \subseteq N_i$, and $\Omega(s) \subseteq N_i$, the state space of \textit{ng}-resources is at most $2^{3|N_i|}$ for $ng$-feasible partial subpaths ending in node $i$, versus $2^{|\Custs|}$ with full elementarity, thus alleviating computational requirements.

\subsection{Tightened relaxations via subset-row cuts}
\label{ssec.cut}

\subsubsection*{Subset-row cuts.}
\cite{jepsen2008subset} defined subset-row cuts (SRCs) as rank-1 Chvátal-Gomory cuts from elementarity constraints (Equation~\eqref{eq.path.serve_customers}). Specifically, for any subset $S \subseteq \Custs$ and weights $\Set{ w_i \in [0, 1) \mid i \in S}$, the following constraints define valid inequalities for $\EVRP(\calP)$:
\begin{align}
    \label{eq.SRC}
    \sum_{i \in S} \sum_{p \in \calP} w_i \gamma_i^p z^p \leq \sum_{i \in S} w_i
    \implies
    \sum_{p \in \calP} \alpha_{S, \bw}(p)  z^p 
    \leq \floor*{\sum_{i \in S} w_i}
    \text{ with }
    \alpha_{S, \bw}(p) = 
    \floor*{\sum_{i \in S} w_i \serve{i}{p}}
\end{align}

In our implementation, to simplify the separation problem, we restrict our attention to SRCs with $|S| = 3$ and $w_i = \frac{1}{2}$ for all $i \in |S|$ (as in \cite{pecin_limited_2017}).

We index the SRCs by $q\in\calQ$, and let $\Set{ (S_q, \bw^q, \lambda_q) | q \in \calQ}$ store the sets $S_q\subseteq\Custs$, the weight vectors $\bw^q$, and the dual variables $\lambda_q \leq 0$ of Equation~\eqref{eq.SRC}. The reduced cost of a path becomes:
\begin{equation}\label{eq.RC_SRC}
    \cbar{p}
    = c^p
    - \kappa_{\snode{p}}
    - \mu_{\enode{p}}
    - \sum_{i \in \Custs} \serve{i}{p} \nu_i
    - \sum_{q \in \calQ} \lambda_q \cdot 
    \alpha_{S_q, \bw^q}(p)
\end{equation}

Note that SRCs are non-robust, as they alter the structure of the pricing problem. Yet, they can be incorporated by tracking the fractional part of $\sum_{i \in S_q} w^q_i \gamma_i^p$, which we denote by $\alpha_q$:
\begin{definition}
    \label{def.subpath_path_resources_SRC}
    Let $(S_q, \bw^q, \lambda_q)$ be a SRC. The SRC resource for a (partial) path or subpath is the one for its node sequence $N = (n_0, \dots, n_m)$:
    $\alpha_q(N) = \fracpart{\sum_{i \in S_q} w^q_i \gamma_i^p}$, where $\fracpart{x} := x - \floor*{x}$. 
\end{definition}

Since $\alpha_q$ tracks the (weighted) number of visits to customers in $S_q$, this can be aggregated at the subpath level and at the path level. Propositions~\ref{prop.subpath_REFs_SRC}--\ref{prop.path_REFs_SRC} define the resources in our subpath-based label-setting algorithm, with domination criteria from Definition~\ref{def.subpath_path_domination_SRC}.
Note that when $\alpha_q$ exceeds 1, $\gamma_q$ is deducted from the reduced cost and $\alpha_q$ is decremented by 1, keeping it in $[0, 1)$.
\begin{proposition}
    \label{prop.subpath_REFs_SRC}
    Let $s$ be a feasible partial subpath and $a = (n, n')$ be an arc extension. Then: $\cbar{s \oplus a}
    = \cbar{s}
    + \overline{c}_{n,n'}
    - \sum_{q \in \calQ}
    \lambda_q \cdot \ind{\alpha_q(s) + \ind{n' \in S_q} w^q_{n'} \geq 1}$ and $\alpha_q(s \oplus a)
    = \fracpart{ \alpha_q(s) + \ind{n' \in S_q} w^q_{n'} },
    \quad \forall \ q \in \calQ$.
\end{proposition}

\begin{proposition}
    \label{prop.path_REFs_SRC}
    Let $p$ be a feasible partial path and $s$ be a subpath extending $p$. Then: $\cbar{p \oplus s}
    =
    \cbar{p}
    + \cbar{s}
    - \sum_{q \in \calQ}
    \lambda_q \cdot \ind{\alpha_q(p) + \alpha_q(s) \geq 1}$ and $\alpha_q(p \oplus s)
    = \fracpart{ \alpha_q(p) + \alpha_q(s) },
    \quad \forall \ q \in \calQ$.
\end{proposition}

\begin{definition}
    \label{def.subpath_path_domination_SRC}
    Let $p_1, p_2$ be feasible partial paths (resp. $s_1, s_2$ feasible partial subpaths) starting and ending at the same nodes.
    $p_1 \succeq p_2$ (resp. $s_1 \succeq s_2$) with subset-row cuts $\calQ$ if the condition $\cbar{p_1} \leq \cbar{p_2}$ (resp. $\cbar{s_1} \leq \cbar{s_2}$)
    is replaced with $\cbar{p_1} - \cbar{p_2} \leq \sum_{q \in \calQ} \lambda_q \ind{\alpha_q(p_1) > \alpha_q(p_2)}$ (resp. $\cbar{s_1} - \cbar{s_2} \leq \sum_{q \in \calQ} \lambda_q \ind{\alpha_q(s_1) > \alpha_q(s_2)}$).
\end{definition}

\begin{proposition}
    \label{prop.SRCs_correctness}
    The definitions and propositions in Theorem~\ref{thm.EVRPTWNL_correctness},
    together with
    Definition
    ~\ref{def.subpath_path_domination_SRC}
    and Propositions
    ~\ref{prop.subpath_REFs_SRC}--\ref{prop.path_REFs_SRC}
    satisfy Property~\ref{property.domination_ours}
    for the EVRPTWNL
    with subset-row cuts $\calQ$.
\end{proposition}

\subsection{Summary}

Algorithm~\ref{alg.AdaptiveColumnGenerationWithCuts} tightens the EVRP relaxation using adaptive \textit{ng}-relaxations to enforce elementarity and subset-row cuts to eliminate fractional solutions. We have developed appropriate resources, REFs, and domination criteria to ensure the validity of our two-phase label-setting algorithm to solve the resulting pricing problems. In particular, resources in traditional (path-based) label-setting algorithms needed to be augmented to propagate domination patterns when combining subpaths into full paths, per Property~\ref{property.domination_ours}. Leveraging these results (Propositions~\ref{prop.ngroute_correctness} and~\ref{prop.SRCs_correctness}) and those from Section~\ref{sec.labelsetting} (Theorem~\ref{thm.EVRPTWNL_correctness}), we obtain the correctness of the overall algorithm below:
\begin{corollary}
    \label{cor.iterative_cg_correctness}
    Algorithm~\ref{alg.AdaptiveColumnGenerationWithCuts} terminates finitely. Steps 1--4 return $\OPTLP(\Pelem)$, and Steps 1--5 return a solution $\texttt{OPT}$ such that $\OPTLP(\Pelem) \leq \texttt{OPT} \leq \OPT(\Pelem)$.
\end{corollary}
\section{Computational results}
\label{sec.numerical_results}

We now implement our methodology to identify and delineate the benefits of the subpath-based label-setting algorithm against path-based benchmarks. We consider both benchmark instances from the literature and new instances developed in this paper to introduce variability across numbers of customers, geographic layouts, charging requirements, and time windows (Section~\ref{subsec.setup}). The results are reported in Sections~\ref{subsec.results_benchmark} and~\ref{subsec.results_new}. All models are solved with Gurobi v12.0, using the JuMP package in Julia v1.11 \citep{dunning2017jump}. All runs are performed using a single thread, on a computing cluster hosting Intel Xeon Platinum 8260 processors \citep{reuther2018interactive}, with a one-hour limit. 
All instances, code and results are made available online 
at \url{github.com/sean-lo/ElectricRouting.jl}
to enable replication.
%

\subsection{Experimental setup}
\label{subsec.setup}

We first use benchmark single-depot instances of the EVRPTW from \cite{schneider2014electric}, also used by \cite{desaulniers2016exact}. We consider the 25-, 50- and 100-customer instances of types C1, R1, and RC1, which correspond to clustered, random, and hybrid locations. We set a fleet of 10, 15 and 25 vehicles, respectively. These instances are typical of last-mile urban logistics, with a dense geographic layout, tight time windows and limited recharging needs. We also scale the width of the time windows and the battery capacity to identify the drivers of performance of the algorithms.

We also create multi-depot EVRPTW instances, inspired by our motivating examples in middle-mile logistics, industrial logistics, and robotic fleet management. We consider a rectangular $x$-by-$y$ grid, with $d$ customers per unit area, depots in the four corners, and charging stations on integer grid points. A fleet of $K = 10$ vehicles starts on the left side ($v_1^\text{start} = v_2^\text{start} = 5$), with some vehicles required to end on the right side ($v_3^\text{end} = v_4^\text{end} = 1$). We consider both tight and loose time windows, both of which increase with the $x$-coordinate. Other parameters are specified as in \cite{schneider2014electric}. This environment define instances where vehicles traverse a widespread service area, with up to 100 customers. We generate five instances for each parameter setting. For both instance types, we consider linear recharging functions: $g(t) = t$ (after rescaling); however, our method applies more generally to nonlinear charging as well. 

\begin{figure}
    \centering
    \begin{subfigure}{0.49\textwidth}
        \centering
        \includegraphics[width=\linewidth]{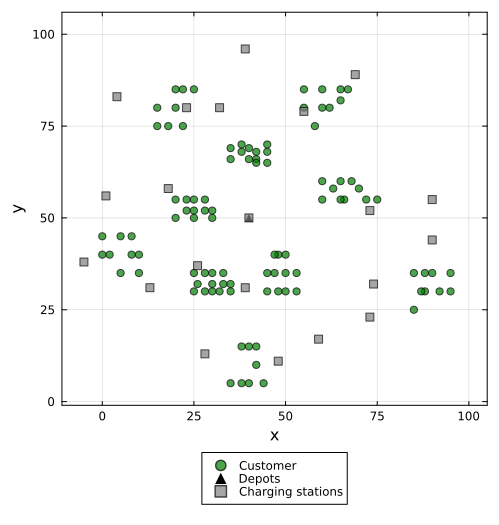}
        \caption{\small C101 instance with 100 customers}
        \label{fig.instance_benchmark}
    \end{subfigure}%
    \hfill
    \begin{subfigure}{0.49\textwidth}
        \centering
        \includegraphics[width=\linewidth]{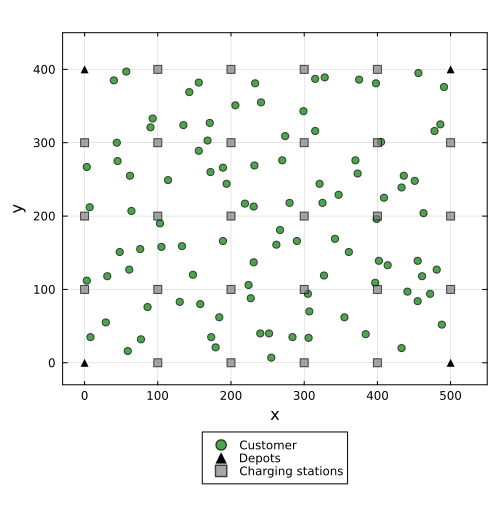}
        \caption{\small Our 100-customer instance \texttt{100\_20\_small\_01}}
        \label{fig.instance_ours}
    \end{subfigure}%
    \caption{Examples of instances from \cite{schneider2014electric} and our new multi-depot instances.}
    \label{fig.instances}
    \vspace{-12pt}
\end{figure}

\begin{table}[h!]
    \centering
    \caption{Parameter settings for the EVRPTW instances from \cite{schneider2014electric} and our new multi-depot EVRPTW instances.}
    \label{tab.instances_schneider_ours}
    \small
    \begin{tabular}{p{6cm} cc}
        \toprule
        Parameter & EVRPTW \citep{schneider2014electric} & Multi-depot EVRPTW \\
        \midrule
        Customer distribution 
        & \begin{tabular}{@{}lr@{}}
            C1 (clustered): & 9 instances \\
            R1 (uniformly random): & 12 instances \\ 
            RC1 (mixed): & 8 instances \\ 
        \end{tabular} 
        & uniformly random \\
        \midrule
        Number of vehicles $K$ & $\{10, 15, 25\}$ & 10 \\
        Number of depots $|\Depots|$ & 1 & 4 \\
        Density $d$ & {--} & $\{2, 3, 4, 5\}$ \\
        Grid size $(x, y)$ & $(1, 1) $ & $\{2, 3, 4, 5\} \times \{4\}$ \\
        Number of customers $|\Custs|$ & $\{25, 50, 100\}$ & $dxy$ \\
        Number of charging stations $|\Stations|$ & 21 & $(x+1)(y+1)-4$ \\
        Time window widths & varied & $\{ \text{narrow}, \text{wide}\}$ \\
        \midrule 
        Number of random instances & 1 & 5 \\
        Total number of instances & 87 & 160 \\
        \bottomrule
    \end{tabular}
\end{table}

Table~\ref{tab.instances_schneider_ours} and Figure~\ref{fig.instances} provide details on these instances. For each one, we apply Algorithm~\ref{alg.AdaptiveColumnGenerationWithCuts} and a corresponding version with a monodirectional implementation of the path-based label setting algorithm from \cite{desaulniers2016exact} to solve the pricing problem. We define the initial \textit{ng}-relaxation with a \textit{ng}-neighborhood for each customer $i$ containing the $5$ closest ones inclusive of $i$ (or $i$ only for a charging station $i$). As \cite{desaulniers2016exact}, we first solve the pricing problem over a restricted graph and resolve it on the full graph if no path of negative reduced cost is found. We evaluate computational time and solution quality acrorss three increasingly tight relaxations: (1) after column generation on the initial \textit{ng}-relaxation (end of Step 1), (2) upon solving the linear relaxation with elementary paths by tightening the \textit{ng}-relaxation (end of Step 4), and (3) upon optimizing over the SRC closure (end of the algorithm or time limit).

\subsection{Results on benchmark instances} 
\label{subsec.results_benchmark}

Table~\ref{tab.ours_benchmark_table} reports results in the EVRPTW instances from \cite{schneider2014electric}. Our two-phase label-setting algorithm is $50\%$ faster than the path-based benchmark on 25-customer instances, while achieving similar integrality gaps. These benefits are consistent across the initial relaxation, the tighter \textit{ng}-relaxation, and the SRC closure. In turn, 27 out of the 29 instances can be solved to optimality. In contrast, our two-phase label-setting algorithm is competitive but no longer provides significant improvements in 50- and 100-customer EVRPTW instances. As described below, this stems from large battery capacity and tight time windows.

\begin{table}
\small
\renewcommand{\arraystretch}{1}
\caption{\footnotesize
Comparison of our two-phase label-setting (2-LS) against the path-based subpath-based label-setting benchmark (LS) across the C1, R1 and RC1 instances: (1) after the initial relaxation, (2) after \textit{ng}-relaxation tightening, and (3) after adding subset-row cuts. We record the number of instances solved to optimality and the geometric means of the computational time, the dual and primal bounds, and the integrality gap.
}
\label{tab.ours_benchmark_table}
\footnotesize
\begin{tabular}{
    @{}
    r@{\hskip6pt}r@{\hskip6pt}l
    *{2}{
        S[table-format=2.0]
        @{\hskip5pt}
        S[table-format=2.0]
        @{\hskip5pt}
        S[table-format=4.2,round-mode=figures,round-precision=3,group-separator={}]
        @{\hskip5pt}
        S[table-format=1.2e1,round-mode=figures,round-precision=3,tight-spacing=true]
        @{\hskip5pt}
        S[table-format=1.2e1,round-mode=figures,round-precision=3,tight-spacing=true]
        @{\hskip4pt}
        S[table-format=2.2\%]
    }
    @{}
}
    \toprule
    & &
    & \multicolumn{6}{c}{LS}
    & \multicolumn{6}{c}{2-LS}
    \\
    \cmidrule(lr){4-9}
    \cmidrule(lr){10-15}
    {Stage} & {$N$} & {Type} 
    & {\#} & {Opt.} & {$t$ (s)} & {Opt. (D)} & {Opt. (P)} & {Gap (\%)}  
    & {\#} & {Opt.} & {$t$ (s)} & {Opt. (D)} & {Opt. (P)} & {Gap (\%)} 
    \\
    \midrule
    \multirow[t]{9}{*}{(1)} 
    & \multirow[t]{3}{*}{25} 
    & {C1} & 9 & 2 & 10.89 & 4.429e+04 & 5.635e+04 & 19.82\% & 9 & 2 & 6.599 & 4.429e+04 & 5.544e+04 & 18.68\% \\ 
    & & {R1} & 12 & 1 & 13.96 & 3.914e+04 & 5.32e+04 & 24.12\% & 12 & 1 & 7.34 & 3.914e+04 & 5.267e+04 & 23.39\% \\ 
    & & {RC1} & 8 & 0 & 11.05 & 4.849e+04 & 6.334e+04 & 23.05\% & 8 & 0 & 4.638 & 4.849e+04 & 6.61e+04 & 25.79\% \\ 
    & \multirow[t]{3}{*}{50} 
    & {C1} & 9 & 1 & 74.97 & 6.605e+04 & 7.21e+04 & 7.79\% & 9 & 1 & 75.4 & 6.605e+04 & 7.228e+04 & 8.01\% \\ 
    & & {R1} & 12 & 0 & 110.3 & 6.731e+04 & 7.947e+04 & 14.87\% & 12 & 0 & 97.04 & 6.731e+04 & 7.958e+04 & 14.90\% \\ 
    & & {RC1} & 8 & 0 & 96.06 & 7.844e+04 & 8.826e+04 & 10.95\% & 8 & 0 & 61.97 & 7.844e+04 & 8.975e+04 & 12.33\% \\ 
    & \multirow[t]{3}{*}{100} 
    & {C1} & 7 & 1 & 844.1 & 9.622e+04 & 1.019e+05 & 5.52\% & 6 & 1 & 1155 & 9.723e+04 & 1.032e+05 & 5.65\% \\ 
    & & {R1} & 8 & 0 & 632 & 1.199e+05 & 1.304e+05 & 7.85\% & 6 & 0 & 695.1 & 1.265e+05 & 1.354e+05 & 6.37\% \\ 
    & & {RC1} & 4 & 0 & 992.3 & 1.411e+05 & 1.534e+05 & 7.99\% & 4 & 0 & 883.5 & 1.411e+05 & 1.535e+05 & 8.03\% \\ 
    \midrule
    \multirow[t]{9}{*}{(2)} 
    & \multirow[t]{3}{*}{25} 
    & {C1} & 9 & 5 & 32.84 & 4.859e+04 & 4.888e+04 & 0.60\% & 9 & 5 & 21.58 & 4.859e+04 & 4.886e+04 & 0.56\% \\ 
    & & {R1} & 12 & 3 & 38.73 & 4.418e+04 & 4.473e+04 & 1.23\% & 12 & 3 & 27.07 & 4.418e+04 & 4.483e+04 & 1.43\% \\ 
    & & {RC1} & 8 & 3 & 35.59 & 5.579e+04 & 5.652e+04 & 1.26\% & 8 & 3 & 19.72 & 5.579e+04 & 5.644e+04 & 1.14\% \\ 
    & \multirow[t]{3}{*}{50} 
    & {C1} & 9 & 5 & 174 & 6.918e+04 & 6.967e+04 & 0.69\% & 9 & 5 & 219 & 6.918e+04 & 6.975e+04 & 0.79\% \\ 
    & & {R1} & 12 & 0 & 332.6 & 7.126e+04 & 7.312e+04 & 2.52\% & 12 & 0 & 318.4 & 7.052e+04 & 7.43e+04 & 4.69\% \\ 
    & & {RC1} & 8 & 1 & 260 & 8.192e+04 & 8.504e+04 & 3.61\% & 8 & 1 & 224.3 & 8.192e+04 & 8.486e+04 & 3.41\% \\ 
    & \multirow[t]{3}{*}{100} 
    & {C1} & 7 & 1 & 1212 & 9.652e+04 & 1.019e+05 & 5.16\% & 6 & 1 & 1447 & 9.75e+04 & 1.03e+05 & 5.28\% \\ 
    & & {R1} & 8 & 0 & 1079 & 1.215e+05 & 1.28e+05 & 4.79\% & 6 & 0 & 907.1 & 1.274e+05 & 1.341e+05 & 4.84\% \\ 
    & & {RC1} & 4 & 0 & 1917 & 1.431e+05 & 1.511e+05 & 5.21\% & 4 & 0 & 2096 & 1.431e+05 & 1.503e+05 & 4.78\% \\ 
    \midrule
    \multirow[t]{9}{*}{(3)} 
    & \multirow[t]{3}{*}{25} 
    & {C1} & 9 & 9 & 41.31 & 4.877e+04 & 4.877e+04 & 0.00\% & 9 & 9 & 29.04 & 4.877e+04 & 4.877e+04 & 0.00\% \\ 
    & & {R1} & 12 & 10 & 113 & 4.462e+04 & 4.467e+04 & 0.12\% & 12 & 10 & 74.01 & 4.461e+04 & 4.467e+04 & 0.14\% \\ 
    & & {RC1} & 8 & 8 & 55.14 & 5.63e+04 & 5.63e+04 & 0.00\% & 8 & 8 & 29.76 & 5.63e+04 & 5.63e+04 & 0.00\% \\ 
    & \multirow[t]{3}{*}{50} 
    & {C1} & 9 & 7 & 252.4 & 6.934e+04 & 6.96e+04 & 0.38\% & 9 & 7 & 288.8 & 6.931e+04 & 6.966e+04 & 0.48\% \\ 
    & & {R1} & 12 & 4 & 1164 & 7.182e+04 & 7.248e+04 & 0.89\% & 12 & 2 & 818 & 7.09e+04 & 7.402e+04 & 3.74\% \\ 
    & & {RC1} & 8 & 2 & 1252 & 8.313e+04 & 8.427e+04 & 1.32\% & 8 & 2 & 851.2 & 8.31e+04 & 8.424e+04 & 1.32\% \\ 
    & \multirow[t]{3}{*}{100} 
    & {C1} & 7 & 1 & 1974 & 9.688e+04 & 1.019e+05 & 4.88\% & 6 & 1 & 1930 & 9.773e+04 & 1.03e+05 & 5.00\% \\ 
    & & {R1} & 8 & 3 & 1780 & 1.218e+05 & 1.277e+05 & 4.27\% & 6 & 1 & 1634 & 1.277e+05 & 1.339e+05 & 4.45\% \\ 
    & & {RC1} & 4 & 0 & 2838 & 1.437e+05 & 1.501e+05 & 4.22\% & 4 & 0 & 3032 & 1.439e+05 & 1.497e+05 & 3.85\% \\ 
    \bottomrule
\end{tabular}
\end{table}

\begin{figure}[h!]
    \centering
    \begin{subfigure}{0.33\textwidth}
        \centering
        \includegraphics[width=\linewidth]{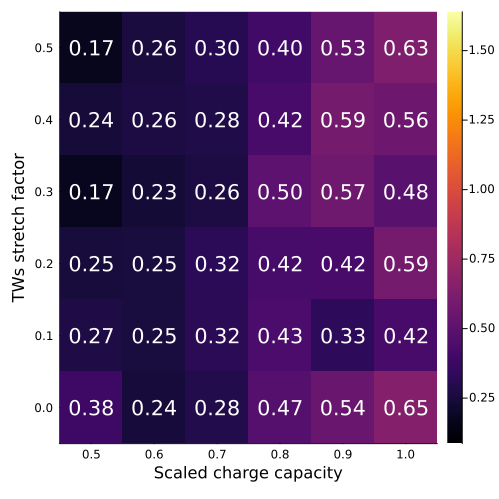}
        \caption{25 customers}
        \label{fig.ours_benchmark_ablation_heatmap_time_taken_first_ratio_n25_c1}
    \end{subfigure}%
    \hfill
    \begin{subfigure}{0.33\textwidth}
        \centering
        \includegraphics[width=\linewidth]{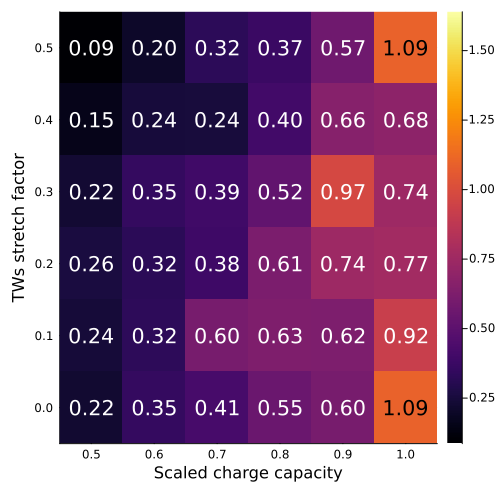}
        \caption{50 customers}
        \label{fig.ours_benchmark_ablation_heatmap_time_taken_first_ratio_n50_c1}
    \end{subfigure}%
    \hfill
    \begin{subfigure}{0.33\textwidth}
        \centering
        \includegraphics[width=\linewidth]{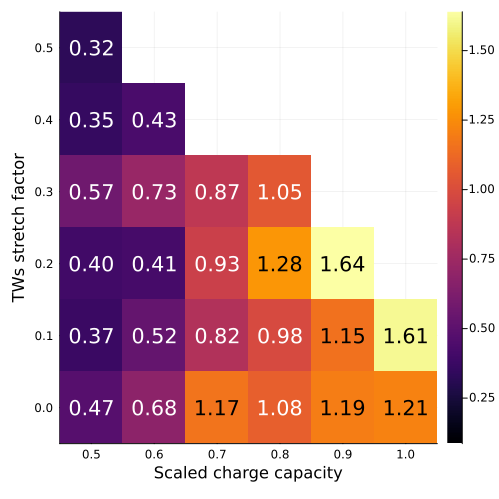}
        \caption{100 customers}
        \label{fig.ours_benchmark_ablation_heatmap_time_taken_first_ratio_n100_c1}
    \end{subfigure}%
    \caption{Ratio of computational times for the initial \textit{ng}-relaxation between our subpath-based label-setting algorithm and the path-based label-setting benchmark, in the C101 instances.}
    \label{fig.ours_benchmark_ablation_heatmap_time_taken_ratio_c1}
    \vspace{-12pt}
\end{figure}

To better understand and characterize these effects, we conduct a sensitivity analysis on battery capacity and time window widths. In each instance, we scale the battery capacity down by a multiplicative factor; smaller battery capacities are also equivalent to longer distances between customers and therefore fewer customer visits per subpath. We also stretch the time windows by a factor $\lambda$ in $[0, 1]$, transforming $[\alpha, \beta]$ into $[(1-\lambda)\alpha, (1-\lambda)\beta + \lambda T]$; wider time windows enable less pruning in the label-setting algorithms. Figure~\ref{fig.ours_benchmark_ablation_heatmap_time_taken_ratio_c1} shows the ratio of the computational times between subpath-based label-setting algorithm and the path-based label-setting benchmark as a function of these two factors.  As noted in Table~\ref{tab.ours_benchmark_table}, our method can lead to performance deterioration with large numbers of customers, large battery capacity and tight time windows. In other settings, our subpath-based label-setting algorithm can lead to significant speedups as:
\begin{itemize}
    \item[--] Battery capacity becomes smaller (or travel distances become longer): in the 100-customer instance, as the battery is scaled down from 100\% to 50\% of its original capacity, our method becomes more than 50\% faster than the path-based benchmark on the initial \textit{ng}-relaxation. Due to more frequent recharging, each path is composed of more subpaths, which strengthens the impact of the subpath-based decomposition in the pricing problem.
    \item[--] Time windows become wider: as the stretch factor $\lambda$ varies from $0$ to $0.5$, the instance becomes more challenging due to more limited pruning in the pricing problem; both methods incur longer computational times, but our subpath-based label-setting algorithm maintains or grows its relative advantage compared to the path-based benchmark. Indeed, generating non-dominated subpaths between charging stations quells the rate of exponential growth, and wider time windows makes subpaths more ``reusable'' in the second phase of the algorithm.
\end{itemize}

\subsection{Results on new multi-depot EVRPTW instances}
\label{subsec.results_new}

\begin{figure}[h!]
    \centering
    \begin{subfigure}{0.50\textwidth}
        \centering
        \includegraphics[width=\linewidth]{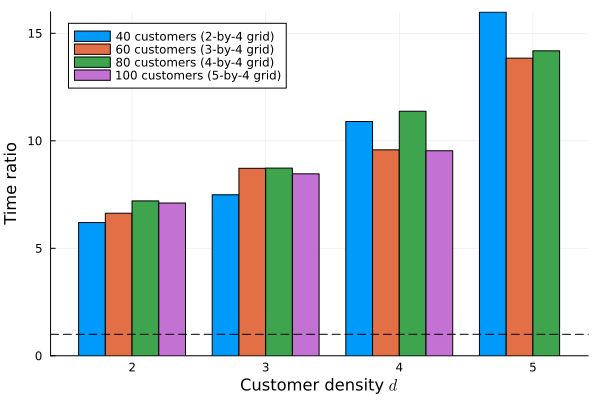}
        \caption{Wider time windows}
        \label{fig.multidepot_barplot_time_taken_first_speedup_big}
    \end{subfigure}%
    \hfill
    \begin{subfigure}{0.50\textwidth}
        \centering
        \includegraphics[width=\linewidth]{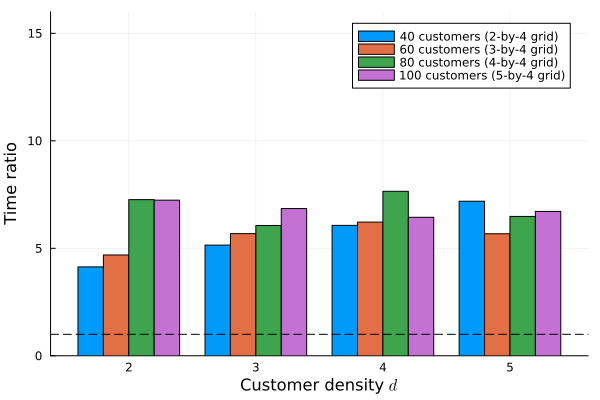}
        \caption{Narrower time windows}
        \label{fig.multidepot_barplot_time_taken_first_speedup_small}
    \end{subfigure}%
    \caption{Computational speedup for the initial \textit{ng}-relaxation between our subpath-based label-setting algorithm and the path-based label-setting benchmark (values higher than one indicate faster solutions).}
    \label{fig.multidepot_barplot_time_taken_first_speedup}
    \vspace{-12pt}
\end{figure}

Figure~\ref{fig.multidepot_barplot_time_taken_first_speedup} reports the computational speedup of our subpath-based label-setting algorithm versus the path-based benchmark on our new multi-depot EVRPTW instances, upon solving the initial \textit{ng}-relaxation (end of Step 1 in the algorithm). Note that our method terminates 5--7 faster than the benchmark with narrower time windows, and 5-15 faster with wider time windows. These benefits increase as the problem size---encoded by customer density---becomes larger. The path-based label-setting benchmark remains tractable with narrow time windows due to extensive pruning, but becomes much more computationally intensive as the number of paths increases due to more customers or wider time windows. In comparison, our subpath-based label-setting algorithm remains much more scalable in challenging instances, underscoring the role of subpath-based decomposition to retain a tractable structure when the number of paths grows exponentially.

\begin{figure}[h!]
    \centering
    \begin{subfigure}{0.50\textwidth}
        \centering
        \includegraphics[width=\linewidth]{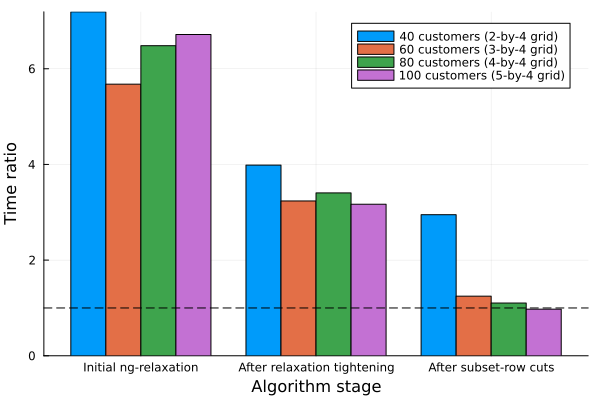}
        \caption{\footnotesize Computational speedup (narrower time windows)}
        \label{fig.multidepot_barplot_time_taken_speedup_small_density5}
    \end{subfigure}%
    \hfill
    \begin{subfigure}{0.50\textwidth}
        \centering
        \includegraphics[width=\linewidth]{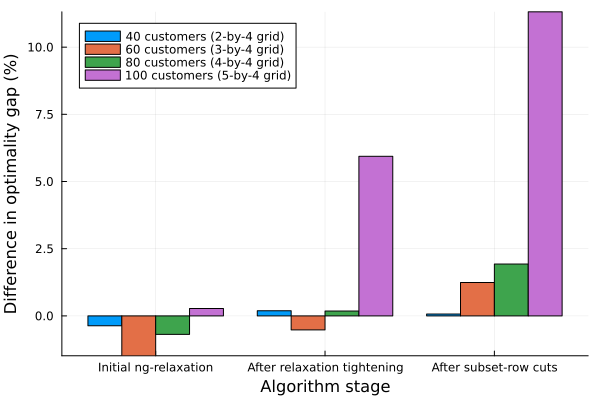}
        \caption{\footnotesize Optimality gap difference (narrower time windows)}
        \label{fig.multidepot_barplot_LP_IP_gap_diff_small_density5}
    \end{subfigure}%
    \newline
    \begin{subfigure}{0.50\textwidth}
        \centering
        \includegraphics[width=\linewidth]{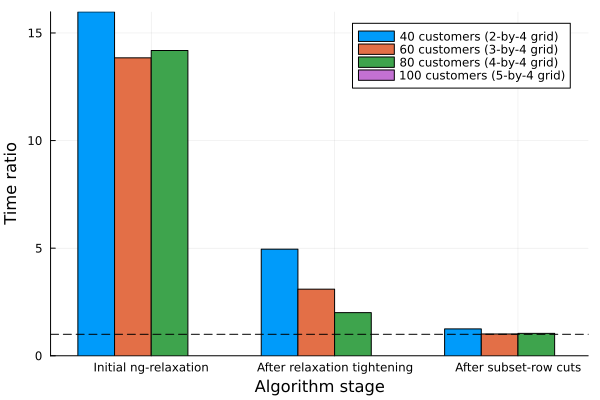}
        \caption{\footnotesize Computational speedup (wider time windows)}
        \label{fig.multidepot_barplot_time_taken_speedup_big_density5}
    \end{subfigure}%
    \hfill
    \begin{subfigure}{0.50\textwidth}
        \centering
        \includegraphics[width=\linewidth]{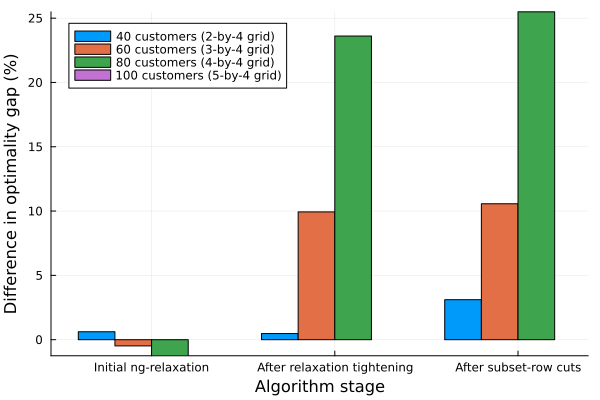}
        \caption{\footnotesize Optimality gap difference (wider time windows)}
        \label{fig.multidepot_barplot_LP_IP_gap_diff_big_density5}
    \end{subfigure}%
    \caption{Computational speedup and difference in optimality gap upon convergence with our subpath-based label-setting algorithm and the path-based label-setting benchmark.}
    \label{fig.multidepot_barplot_time_taken_speedup_LP_IP_gap_diff}
    \vspace{-12pt}
\end{figure}

Finally, we implement the full column generation algorithm to solve the initial \textit{ng}-relaxation, the tighter \textit{ng}-relaxation, and the SRC closure. Results are summarized in Figure~\ref{fig.multidepot_barplot_time_taken_speedup_LP_IP_gap_diff} and detailed in Tables~\ref{tab.multidepot_small}--\ref{tab.multidepot_big}. As expected from Figure~\ref{fig.multidepot_barplot_time_taken_first_speedup}, our method results in a significant speedup in the first step of the algorithm to solve the initial \textit{ng}-relaxation. These benefits translate into stronger eventual performance. Upon tightening the \textit{ng}-relaxation and adding SRCs, computational times remain faster in smaller instances. In larger instances, the speedups enable more iterations and therefore more progress toward tighter primal and dual bounds, leading to a much tighter optimality gap. In particular, our method leads to up to a 11 percentage-point improvement with narrower time windows and to a 25 percentage-point improvement with wider time windows. Importantly, these benefits are driven by both improvements in primal solutions (up to 9\% and 18\% cost reductions, respectively) and dual bounds (3\% and 13\% tighter bounds, respectively).

\begin{figure}[h!]
    \centering
    \begin{subfigure}{0.50\textwidth}
        \centering
        \includegraphics[width=\linewidth]{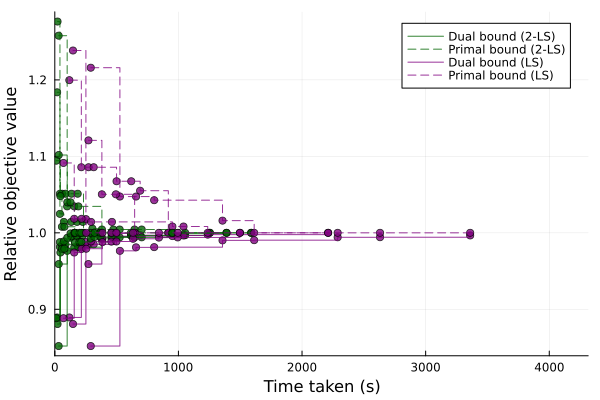}
        \caption{40 customers (narrower time windows)}
        \label{fig.multidepot_bounds_n40_2.0_5_small}
    \end{subfigure}%
    \begin{subfigure}{0.50\textwidth}
        \centering
        \includegraphics[width=\linewidth]{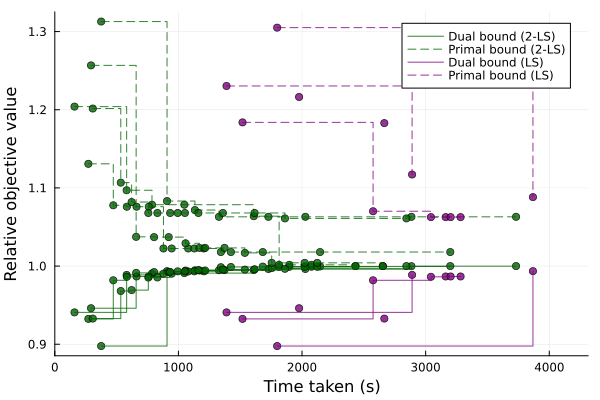}
        \caption{100 customers (narrower time windows)}
        \label{fig.multidepot_bounds_n100_5.0_5_small}
    \end{subfigure}%
    \newline
    \begin{subfigure}{0.50\textwidth}
        \centering
        \includegraphics[width=\linewidth]{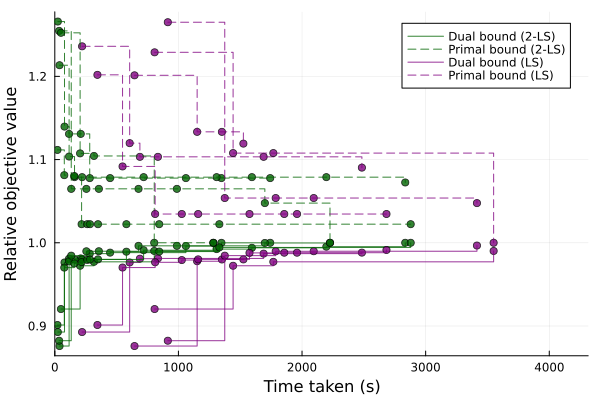}
        \caption{40 customers (wider time windows)}
        \label{fig.multidepot_bounds_n40_2.0_5_big}
    \end{subfigure}%
    \hfill
    \begin{subfigure}{0.50\textwidth}
        \centering
        \includegraphics[width=\linewidth]{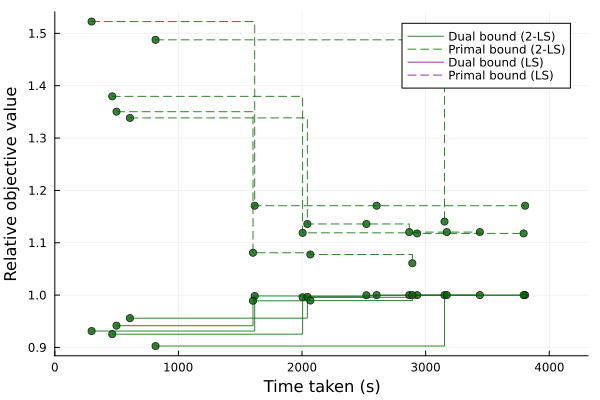}
        \caption{100 customers (wider time windows)}
        \label{fig.multidepot_bounds_n100_5.0_5_big}
    \end{subfigure}%
    \caption{Convergence of the subpath-based label-setting algorithm and the path-based label-setting benchmark, normalized to the best known dual bound.}
    \label{fig.multidepot_bounds}
    \vspace{-12pt}
\end{figure}

These benefits are further observed in Figure~\ref{fig.multidepot_bounds}, which tracks the evolution of the primal and dual bounds over time with both methods. These visualizations underscore the speedups from our subpath-based label-setting algorithm in the smaller 40-customer instances, and the improvements in the primal solutions and dual bounds---hence, in optimality gaps---in the 100-customer instances.

\subsection{Summary}

These computational results compare our subpath-based label-setting method against the path-based label-setting benchmark from \cite{desaulniers2016exact}. Our method is competitive in the original EVRPTW instances from \cite{schneider2014electric} in both computational time and solution quality. At the same time, these instances featured limited recharging (e.g., a vehicle can perform its daily duty with one or two recharging stops) and narrow time windows (e.g., urban logistics where each customer can be served during a short interval), which leads to extensive pruning in path-based label-setting algorithms. Our method becomes more beneficial with smaller battery capacities---or, equivalently, with a lower spatial density---and wider time windows. Finally, our new multi-depot instances underscore that our method substantially outperforms the benchmark in both runtime and solution quality when customers' spatial locations are more dispersed and time windows allow subpaths to be reused across paths---e.g., to travel from one area to another with charging stations in between. Altogether, the proposed method can lead to performance improvements and scales to otherwise intractable instances with up to 100 customers (Table~\ref{tab.arc}).

Altogether, our method generates strongest comparative benefits when time windows become wider, paths become longer (many customers within the planning horizon), each path is constituted by a large enough number of subpaths (charging stations visited frequently), and subpaths become more ``reusable'' across paths. These conditions arise in emerging electrified operations from our motivating examples. For instance, regional distribution and industrial logistics typically feature geographically dispersed facilities, limited charging infrastructure, and loose delivery windows, allowing routing decisions to interact strongly with charging opportunities. Similarly, energy-constrained last-mile systems and autonomous robotic fleets operate under tight battery limitations, requiring frequent recharging. Ultimately, our two-phase label-setting algorithm can provide a scalable solution method across multiple EVRP regimes, and support electrified logistics and robotic fleet management systems where routing and charging decisions must be optimized jointly.
\section{Conclusion}
\label{sec.conclusion}

This paper considers a multi-depot electric vehicle routing problem with time windows and nonlinear charging (EVRPTWNL). The problem jointly optimizes routing and charging decisions for an electric fleet with flexibility regarding where, when and for how long to charge. We develop a column generation methodology based on a two-phase, subpath-based label-setting algorithm for the pricing problem. In the first phase, subpaths are generated between charging stations; in the second phase, subpaths are combined into full paths while optimizing charging decisions in between. We design a new framework of subpath and path domination through the design of resources, REFs, and domination criteria to accommodate elements such as customer time windows and nonlinear charging, as well as relaxation-tightening strategies such as \textit{ng}-routes and subset-row cuts. We prove that the resulting subpath-based label-setting algorithm solves the pricing problem exactly, thanks in part to the new domination framework developed in this paper, thereby enabling finite convergence of the column generation algorithm for EVRPTWNL relaxations.

We perform extensive computational experiments on publicly available benchmark instances and new multi-depot instances. First, the two-phase label-setting algorithm can provide improvements in both computation time and solution quality over a path-based label-setting benchmark. Second, improvements are stronger as time windows become wider, battery capacities become smaller, paths become longer, and paths are decomposed into more subpaths. Third, our proposed method can scale to otherwise intractable instances with up to 100 customers. In turn, the proposed methodology can support emerging electrified operations across logistics domains and robotics systems.

This paper also opens opportunities for future research to extend the subpath-based label-setting algorithm to broader classes of problems with hierarchical structures. In particular, the algorithm could be augmented with various acceleration strategies, such as bidirectional label-setting (both in the first phase and the second phase) and multi-threading, as well as with branch-and-price-and-cut schemes. It could also be applied to other vehicle routing and machine scheduling problems to capture, for instance, load-dependent or speed-dependent charge consumptions, capacitated charging stations, hybrid fleets of electric and gasoline vehicles, or heterogeneity in charging technologies across charging stations. More broadly, the electric routing algorithms developed in this paper could also be integrated into upstream strategic planning to support the development of charging infrastructure or the acquisition of electrified fleets. The subpath-based decomposition method developed in this paper provides methodological foundations to address these broader problems.

\ACKNOWLEDGMENT{
The authors would like to thank the review team for their comprehensive comments, 
as well as Xinwu Qian and Bart van Rossum for their helpful suggestions.
}

{\scriptsize
\bibliographystyle{informs2014}
\bibliography{ref.bib}

@article{
    schneider2014electric,
    author = {Schneider, Michael and Stenger, Andreas and Goeke, Dominik},
    title = {The Electric Vehicle-Routing Problem with Time Windows and Recharging Stations},
    journal = {Transportation Science},
    volume = {48},
    number = {4},
    pages = {500-520},
    year = {2014},
}

@article{
    desaulniers2016exact,
    author = {Desaulniers, Guy and Errico, Fausto and Irnich, Stefan and Schneider, Michael},
    title = {Exact Algorithms for Electric Vehicle-Routing Problems with Time Windows},
    journal = {Operations Research},
    volume = {64},
    number = {6},
    pages = {1388-1405},
    year = {2016},
}

@article{
    baldacci2011new,
    author = {Roberto Baldacci and Aristide Mingozzi and Roberto Roberti},
    journal = {Operations Research},
    number = {5},
    pages = {1269--1283},
    publisher = {INFORMS},
    title = {New Route Relaxation and Pricing Strategies for the Vehicle Routing Problem},
    volume = {59},
    year = {2011}
}

@article{
    kullman2021electric,
    author = {Kullman, Nicholas D. and Goodson, Justin C. and Mendoza, Jorge E.},
    title = {Electric Vehicle Routing with Public Charging Stations},
    journal = {Transportation Science},
    volume = {55},
    number = {3},
    pages = {637-659},
    year = {2021},
}

@article{
    froger2022electric,
    author = {Froger, Aur\'{e}lien and Jabali, Ola and Mendoza, Jorge E. and Laporte, Gilbert},
    title = {The Electric Vehicle Routing Problem with Capacitated Charging Stations},
    journal = {Transportation Science},
    volume = {56},
    number = {2},
    pages = {460-482},
    year = {2022},
}

@article{
    dror1994note,
    author = {Dror, Moshe},
    title = {Note on the Complexity of the Shortest Path Models for Column Generation in VRPTW},
    journal = {Operations Research},
    volume = {42},
    number = {5},
    pages = {977-978},
    year = {1994},
}

@article{
    jepsen2008subset, 
    title={Subset-Row Inequalities Applied to the Vehicle-Routing Problem with Time Windows},
    volume={56}, 
    number={2}, 
    journal={Operations Research}, 
    publisher={INFORMS}, 
    author={Jepsen, Mads and Petersen, Bjørn and Spoorendonk, Simon and Pisinger, David}, 
    year={2008}, 
    month={Apr}, 
    pages={497–511},
}

@article{
    martinelli2014efficient, 
    title={Efficient elementary and restricted non-elementary route pricing}, volume={239}, 
    number={1}, 
    journal={European Journal of Operational Research}, 
    author={Martinelli, Rafael and Pecin, Diego and Poggi, Marcus}, 
    year={2014}, 
    month=nov, 
    pages={102–111} 
}

@article{
    beasley1989algorithm, 
    title={An algorithm for the resource constrained shortest path problem}, 
    volume={19}, 
    rights={Copyright © 1989 Wiley Periodicals, Inc., A Wiley Company}, 
    number={4}, 
    journal={Networks}, 
    author={Beasley, J. E. and Christofides, N.}, 
    year={1989}, 
    pages={379–394}, 
    language={en},
}

@article{
    ralphs2003capacitated,
    title={On the capacitated vehicle routing problem},
    author={Ralphs, Ted K and Kopman, Leonid and Pulleyblank, William R and Trotter, Leslie E},
    journal={Mathematical programming},
    volume={94},
    pages={343--359},
    year={2003},
    publisher={Springer}
}

@inproceedings{
    reuther2018interactive,
    title={Interactive supercomputing on 40,000 cores for machine learning and data analysis},
    author={Reuther, Albert and Kepner, Jeremy and Byun, Chansup and Samsi, Siddharth and Arcand, William and Bestor, David and Bergeron, Bill and Gadepally, Vijay and Houle, Michael and Hubbell, Matthew and Jones, Michael and Klein, Anna and Milechin, Lauren and Mullen, Julia and Prout, Andrew and Rosa, Antonio and Yee, Charles and Michaleas, Peter},
    booktitle={2018 IEEE High Performance extreme Computing Conference (HPEC)},
    pages={1--6},
    year={2018},
    organization={IEEE}
}

@misc{
    iea,
    title="{Energy Technology Perspectives 2020}",
    author="{International Energy Agency}",
    year="2020",
    note={\url{https://iea.blob.core.windows.net/assets/7f8aed40-89af-4348-be19-c8a67df0b9ea/Energy_Technology_Perspectives_2020_PDF.pdf}}
}

@article{
    fernandez2022arc, 
    title={Arc Routing with Electric Vehicles: Dynamic Charging and Speed-Dependent Energy Consumption}, 
    volume={56}, 
    number={5}, 
    journal={Transportation Science}, 
    publisher={INFORMS}, 
    author={Fernández, Elena and Leitner, Markus and Ljubić, Ivana and Ruthmair, Mario}, 
    year={2022}, 
    month=sep, 
    pages={1219–1237}
}

@article{
    nejad2016optimal, 
    title={Optimal Routing for Plug-In Hybrid Electric Vehicles}, 
    volume={51}, 
    number={4}, 
    journal={Transportation Science}, 
    publisher={INFORMS}, 
    author={Nejad, Mark M. and Mashayekhy, Lena and Grosu, Daniel and Chinnam, Ratna Babu}, 
    year={2017}, 
    month=nov, 
    pages={1304–1325}
}

@article{
    sweda2017adaptive, 
    title={Adaptive Routing and Recharging Policies for Electric Vehicles}, 
    volume={51}, 
    number={4}, 
    journal={Transportation Science}, 
    publisher={INFORMS}, 
    author={Sweda, Timothy M. and Dolinskaya, Irina S. and Klabjan, Diego}, 
    year={2017}, 
    month=nov, 
    pages={1326–1348}
}

@article{parmentier2023electric,
  title={Electric vehicle fleets: Scalable route and recharge scheduling through column generation},
  author={Parmentier, Axel and Martinelli, Rafael and Vidal, Thibaut},
  journal={Transportation Science},
  volume={57},
  number={3},
  pages={631--646},
  year={2023},
  publisher={Informs}
}

@article{
    andelmin2017exact, 
    title={An Exact Algorithm for the Green Vehicle Routing Problem}, 
    volume={51}, 
    number={4}, 
    journal={Transportation Science}, 
    publisher={INFORMS}, 
    author={Andelmin, Juho and Bartolini, Enrico}, year={2017}, 
    month=nov, 
    pages={1288–1303} 
}

@article{
    kinay2021charging, 
    title={Charging Station Location and Sizing for Electric Vehicles Under Congestion}, 
    volume={57}, 
    journal={Transportation Science},
    publisher={INFORMS}, 
    author={Kınay, {\"O}mer Burak and Gzara, Fatma and Alumur, Sibel A.}, 
    year={2023}, 
    month=nov, 
    pages={1433–1451} 
}

@article{
    molenbruch2023electric, 
    title={The Electric Dial-a-Ride Problem on a Fixed Circuit}, 
    volume={57}, 
    journal={Transportation Science},
    publisher={INFORMS}, 
    author={Molenbruch, Yves and Braekers, Kris and Eisenhandler, Ohad and Kaspi, Mor}, 
    year={2023}, 
    month=may, 
    pages={594–612} 
}

@article{
    devos2022electric, 
    title={Electric Vehicle Scheduling in Public Transit with Capacitated Charging Stations}, 
    volume={58}, 
    number={2}, 
    journal={Transportation Science}, 
    publisher={INFORMS}, 
    author={de Vos, Marelot H. and van Lieshout, Rolf N. and Dollevoet, Twan}, 
    year={2024}, 
    month=mar, 
    pages={279–294} 
}

@article{
    arslan2019branch, 
    title={A Branch-and-Cut Algorithm for the Alternative Fuel Refueling Station Location Problem with Routing}, 
    volume={53}, 
    number={4}, 
    journal={Transportation Science}, publisher={INFORMS}, 
    author={Arslan, Okan and Karaşan, Oya Ekin and Mahjoub, A. Ridha and Yaman, Hande}, 
    year={2019}, 
    month=jul, 
    pages={1107–1125} 
}

@article{
    mak2013infrastructure, 
    title={Infrastructure Planning for Electric Vehicles with Battery Swapping}, 
    volume={59}, 
    number={7}, 
    journal={Management Science}, 
    publisher={INFORMS}, 
    author={Mak, Ho-Yin and Rong, Ying and Shen, Zuo-Jun Max}, 
    year={2013}, 
    month=jul, 
    pages={1557–1575} 
}

@article{
    brandstatter2020location, 
    title={Location of Charging Stations in Electric Car Sharing Systems}, 
    volume={54}, 
    number={5}, 
    journal={Transportation Science}, 
    publisher={INFORMS}, 
    author={Brandstätter, Georg and Leitner, Markus and Ljubić, Ivana}, 
    year={2020}, 
    month=sep, 
    pages={1408–1438} 
}

@article{
    qi2023scaling, 
    title={Scaling Up Electric-Vehicle Battery Swapping Services in Cities: A Joint Location and Repairable-Inventory Model}, 
    volume={69}, 
    number={11}, 
    journal={Management Science}, 
    publisher={INFORMS}, 
    author={Qi, Wei and Zhang, Yuli and Zhang, Ningwei}, 
    year={2023}, 
    month=nov, 
    pages={6855–6875} }

@article{
    felipe2014heuristic, 
    title={A heuristic approach for the green vehicle routing problem with multiple technologies and partial recharges}, 
    volume={71}, 
    journal={Transportation Research Part E: Logistics and Transportation Review},
    author={Felipe, {\'A}ngel and Ortu{\~n}o, M. Teresa and Righini, Giovanni and Tirado, Gregorio}, 
    year={2014}, 
    month=nov, 
    pages={111–128}, 
    language={en} 
}

@article{
    goeke2015routing, 
    title={Routing a mixed fleet of electric and conventional vehicles}, 
    volume={245}, 
    number={1}, 
    journal={European Journal of Operational Research}, 
    author={Goeke, Dominik and Schneider, Michael}, 
    year={2015}, 
    month=aug, 
    pages={81–99}, 
    language={en} 
}

@article{
    he2013optimal, 
    title={Optimal deployment of public charging stations for plug-in hybrid electric vehicles}, 
    volume={47}, 
    journal={Transportation Research Part B: Methodological}, 
    author={He, Fang and Wu, Di and Yin, Yafeng and Guan, Yongpei}, 
    year={2013}, 
    month=jan, 
    pages={87–101}, 
    language={en} 
}

@article{
    montoya2017electric, 
    series={Green Urban Transportation}, 
    title={The electric vehicle routing problem with nonlinear charging function}, 
    volume={103}, 
    journal={Transportation Research Part B: Methodological}, 
    author={Montoya, Alejandro and Guéret, Christelle and Mendoza, Jorge E. and Villegas, Juan G.}, 
    year={2017}, 
    month=sep, 
    pages={87–110}, 
    collection={Green Urban Transportation}, 
    language={en} 
}

@inbook{
    kallehauge2005vehicle, 
    address={Boston, MA}, 
    title={Vehicle Routing Problem with Time Windows}, 
    booktitle={Column Generation}, 
    publisher={Springer US}, 
    author={Kallehauge, Brian and Larsen, Jesper and Madsen, Oli B.G. and Solomon, Marius M.}, 
    editor={Desaulniers, Guy and Desrosiers, Jacques and Solomon, Marius M.}, 
    year={2005}, 
    pages={67–98} 
}

@article{
    dunning2017jump,
    title={{JuMP: A modeling language for mathematical optimization}},
    author={Dunning, Iain and Huchette, Joey and Lubin, Miles},
    journal={SIAM Review},
    volume={59},
    number={2},
    pages={295--320},
    year={2017},
    publisher={SIAM}
}

@article{
    erdogan2012green, 
    series={Select Papers from the 19th International Symposium on Transportation and Traffic Theory}, 
    title={A Green Vehicle Routing Problem}, 
    volume={48}, 
    number={1}, 
    journal={Transportation Research Part E: Logistics and Transportation Review}, 
    author={Erdoğan, Sevgi and Miller-Hooks, Elise}, 
    year={2012}, 
    month=jan, 
    pages={100–114}, 
    collection={Select Papers from the 19th International Symposium on Transportation and Traffic Theory} 
}

@inbook{
    irnich2005shortest, 
    address={Boston, MA}, 
    title={Shortest Path Problems with Resource Constraints}, 
    booktitle={Column Generation}, 
    publisher={Springer US}, 
    author={Irnich, Stefan and Desaulniers, Guy}, 
    editor={Desaulniers, Guy and Desrosiers, Jacques and Solomon, Marius M.}, 
    year={2005}, 
    pages={33–65}, 
    language={en} 
}

@article{
    barnhart1998branch, 
    title={Branch-and-Price: Column Generation for Solving Huge Integer Programs}, 
    volume={46}, 
    number={3}, 
    journal={Operations Research}, 
    publisher={INFORMS}, 
    author={Barnhart, Cynthia and Johnson, Ellis L. and Nemhauser, George L. and Savelsbergh, Martin W. P. and Vance, Pamela H.}, 
    year={1998}, 
    pages={316–329} 
}

@article{
    righini2006symmetry, 
    series={Graphs and Combinatorial Optimization}, 
    title={Symmetry helps: Bounded bi-directional dynamic programming for the elementary shortest path problem with resource constraints}, 
    volume={3}, 
    number={3}, 
    journal={Discrete Optimization}, 
    author={Righini, Giovanni and Salani, Matteo}, 
    year={2006}, 
    month=sep, 
    pages={255–273}, 
    collection={Graphs and Combinatorial Optimization}, 
    language={en} 
}

@article{
    hasan2021benefits,
    title={The benefits of autonomous vehicles for community-based trip sharing},
    author={Hasan, Mohd Hafiz and van Hentenryck, Pascal},
    journal={Transportation Research Part C: Emerging Technologies},
    volume={124},
    pages={102929},
    year={2021},
    publisher={Elsevier}
}

@article{
    alyasiry2019exact,
    title={An exact algorithm for the pickup and delivery problem with time windows and last-in-first-out loading},
    author={Alyasiry, Ali Mehsin and Forbes, Michael and Bulmer, Michael},
    journal={Transportation Science},
    volume={53},
    number={6},
    pages={1695--1705},
    year={2019},
    publisher={INFORMS}
}

@article{
    rist2021new,
    title={A new formulation for the dial-a-ride problem},
    author={Rist, Yannik and Forbes, Michael A},
    journal={Transportation Science},
    volume={55},
    number={5},
    pages={1113--1135},
    year={2021},
    publisher={INFORMS}
}

@article{
    zhang2023routing,
    title={Routing optimization with vehicle--customer coordination},
    author={Zhang, Wei and Jacquillat, Alexandre and Wang, Kai and Wang, Shuaian},
    journal={Management Science},
    volume={69},
    number={11},
    pages={6876--6897},
    year={2023},
    publisher={INFORMS}
}

@article{
    cummings2024deviated,
    title={Deviated Fixed-route Microtransit: Design and Operations},
    author={Cummings, Kayla and Jacquillat, Alexandre and Martin-Iradi, Bernardo},
    journal={arXiv preprint arXiv:2402.01265},
    year={2024}
}

@article{
    rist2022column,
    title={A column generation and combinatorial benders decomposition algorithm for the selective dial-a-ride-problem},
    author={Rist, Yannik and Forbes, Michael A},
    journal={Computers \& Operations Research},
    volume={140},
    pages={105649},
    year={2022},
    publisher={Elsevier}
}

@article{
    gschwind2019stabilized,
    title = {Stabilized branch-price-and-cut for the commodity-constrained split delivery vehicle routing problem},
    journal = {European Journal of Operational Research},
    volume = {278},
    number = {1},
    pages = {91-104},
    year = {2019},
    issn = {0377-2217},
    doi = {https://doi.org/10.1016/j.ejor.2019.04.008},
    url = {https://www.sciencedirect.com/science/article/pii/S0377221719303285},
    author = {Timo Gschwind and Nicola Bianchessi and Stefan Irnich},
}

@article{nafstad_branch-price-and-cut_2025,
    title = {Branch-{Price}-and-{Cut} for the {Electric} {Vehicle} {Routing} {Problem} with {Heterogeneous} {Recharging} {Technologies} and {Nonlinear} {Recharging} {Functions}},
    volume = {59},
    issn = {0041-1655},
    url = {https://pubsonline.informs.org/doi/10.1287/trsc.2024.0725},
    doi = {10.1287/trsc.2024.0725},
    number = {3},
    urldate = {2025-11-11},
    journal = {Transportation Science},
    publisher = {INFORMS},
    author = {Nafstad, Gaute Messel and Desaulniers, Guy and Stålhane, Magnus},
    month = may,
    year = {2025},
    pages = {628--646},
}

@article{irnich_resource_2008,
    title = {Resource extension functions: properties, inversion, and generalization to segments},
    volume = {30},
    issn = {1436-6304},
    shorttitle = {Resource extension functions},
    url = {https://doi.org/10.1007/s00291-007-0083-6},
    doi = {10.1007/s00291-007-0083-6},
    language = {en},
    number = {1},
    urldate = {2026-02-18},
    journal = {OR Spectrum},
    author = {Irnich, Stefan},
    month = jan,
    year = {2008},
    pages = {113--148},
}

@article{lam_branch-and-cut-and-price_2022,
    title = {Branch-and-cut-and-price for the {Electric} {Vehicle} {Routing} {Problem} with {Time} {Windows}, {Piecewise}-{Linear} {Recharging} and {Capacitated} {Recharging} {Stations}},
    volume = {145},
    issn = {0305-0548},
    url = {https://www.sciencedirect.com/science/article/pii/S0305054822001423},
    doi = {https://doi.org/10.1016/j.cor.2022.105870},
    journal = {Computers \& Operations Research},
    author = {Lam, Edward and Desaulniers, Guy and Stuckey, Peter J.},
    year = {2022},
    pages = {105870},
}

@article{lera-romero_branch-cut-and-price_2023,
    title = {A {Branch}-{Cut}-and-{Price} {Algorithm} for the {Time}-{Dependent} {Electric} {Vehicle} {Routing} {Problem} with {Time} {Windows}},
    volume = {312},
    doi = {10.1016/j.ejor.2023.06.037},
    journal = {European Journal of Operational Research},
    author = {Lera-Romero, Gonzalo and Bront, Juan and Soulignac, Francisco},
    month = jul,
    year = {2023},
}

@article{lee_exact_2021,
    title = {An exact algorithm for the electric-vehicle routing problem with nonlinear charging time},
    volume = {72},
    issn = {0160-5682},
    url = {https://doi.org/10.1080/01605682.2020.1730250},
    doi = {10.1080/01605682.2020.1730250},
    number = {7},
    urldate = {2025-11-10},
    journal = {Journal of the Operational Research Society},
    publisher = {Taylor \& Francis},
    author = {Lee, Chungmok},
    month = jul,
    year = {2021},
    note = {\_eprint: https://doi.org/10.1080/01605682.2020.1730250},
    pages = {1461--1485},
}

@article{ceselli_branch-and-cut-and-price_2021,
    title = {A {Branch}-and-{Cut}-and-{Price} {Algorithm} for the {Electric} {Vehicle} {Routing} {Problem} with {Multiple} {Technologies}},
    volume = {2},
    issn = {2662-2556},
    url = {https://doi.org/10.1007/s43069-020-00052-x},
    doi = {10.1007/s43069-020-00052-x},
    language = {en},
    number = {1},
    urldate = {2025-11-11},
    journal = {Operations Research Forum},
    author = {Ceselli, Alberto and Felipe, Angel and Ortu{~n}o, M. Teresa and Righini, Giovanni and Tirado, Gregorio},
    month = jan,
    year = {2021},
    pages = {8},
}

@article{bezzi_route-based_2023,
    title = {A route-based algorithm for the electric vehicle routing problem with multiple technologies},
    volume = {157},
    issn = {0968-090X},
    url = {https://www.sciencedirect.com/science/article/pii/S0968090X23003649},
    doi = {10.1016/j.trc.2023.104374},
    urldate = {2025-11-11},
    journal = {Transportation Research Part C: Emerging Technologies},
    author = {Bezzi, Dario and Ceselli, Alberto and Righini, Giovanni},
    month = dec,
    year = {2023},
    pages = {104374},
}

@techreport{
    icct_2022_charging,
    author       = {{International Council on Clean Transportation}},
    title        = {Charging Solutions for Battery-Electric Trucks},
    year         = {2022},
    month        = {dec},
    address      = {Washington, DC},
    url          = {https://theicct.org/publication/charging-infrastructure-trucks-zeva-dec22/},
}

@article{boysen_energy_2025,
    title = {Energy management for electric vehicles in facility logistics: {A} survey from an operational research perspective},
    issn = {0377-2217},
    shorttitle = {Energy management for electric vehicles in facility logistics},
    url = {https://www.sciencedirect.com/science/article/pii/S0377221725009993},
    doi = {10.1016/j.ejor.2025.12.031},
    urldate = {2026-03-01},
    journal = {European Journal of Operational Research},
    author = {Boysen, Nils and Schneider, Michael and Žulj, Ivan},
    month = dec,
    year = {2025},
}

@techreport{
    boston2025cargo,
    author       = {{Boston Region MPO}},
    title        = {Exploring Cargo E-Bikes for Last-Mile Deliveries},
    institution  = {Boston MPO},
    year         = {2025},
    month        = {oct},
    url          = {https://test.bostonmpo.org/data/pdf/studies/bikeped/Cargo-Ebikes-Report.pdf}
}

@article{boysen_last-mile_2021,
    title = {Last-mile delivery concepts: a survey from an operational research perspective},
    volume = {43},
    issn = {1436-6304},
    shorttitle = {Last-mile delivery concepts},
    url = {https://doi.org/10.1007/s00291-020-00607-8},
    doi = {10.1007/s00291-020-00607-8},
    language = {en},
    number = {1},
    urldate = {2026-03-01},
    journal = {OR Spectrum},
    author = {Boysen, Nils and Fedtke, Stefan and Schwerdfeger, Stefan},
    month = mar,
    year = {2021},
    pages = {1--58},
}

@misc{
    changi2026autonomous,
    author      = {{Changi Airport Group}},
    title       = {Changi Airport rolls out autonomous tractors in major step towards airside automation},
    year        = {2026},
    month       = {jan},
    url         = {https://www.changiairport.com/en/corporate/our-media-hub/newsroom/press-releases.autonomous-tractors.2026.all.html}
}

@article{valner_scalable_2022,
    title = {Scalable and heterogenous mobile robot fleet-based task automation in crowded hospital environments—a field test},
    volume = {9},
    issn = {2296-9144},
    url = {https://pmc.ncbi.nlm.nih.gov/articles/PMC9445435/},
    doi = {10.3389/frobt.2022.922835},
    urldate = {2026-03-01},
    journal = {Frontiers in Robotics and AI},
    author = {Valner, Robert and Masnavi, Houman and Rybalskii, Igor and Põlluäär, Rauno and Kõiv, Erik and Aabloo, Alvo and Kruusamäe, Karl and Singh, Arun Kumar},
    month = aug,
    year = {2022},
    pages = {922835},
}

@article{harik_design_2021,
    title = {Design and {Implementation} of an {Autonomous} {Charging} {Station} for {Agricultural} {Electrical} {Vehicles}},
    volume = {11},
    copyright = {http://creativecommons.org/licenses/by/3.0/},
    issn = {2076-3417},
    url = {https://www.mdpi.com/2076-3417/11/13/6168},
    doi = {10.3390/app11136168},
    language = {en},
    number = {13},
    urldate = {2026-03-01},
    journal = {Applied Sciences},
    publisher = {publisher},
    author = {Harik, El Houssein Chouaib},
    month = jul,
    year = {2021},
}

@article{birolini2023day,
  title={Day-ahead aircraft routing with data-driven primary delay predictions},
  author={Birolini, Sebastian and Jacquillat, Alexandre},
  journal={European Journal of Operational Research},
  volume={310},
  number={1},
  pages={379--396},
  year={2023},
  publisher={Elsevier}
}

@article{gopalan_aircraft_1998,
    title = {The {Aircraft} {Maintenance} {Routing} {Problem}},
    volume = {46},
    issn = {0030-364X},
    url = {https://www.jstor.org/stable/222864},
    number = {2},
    urldate = {2026-03-06},
    journal = {Operations Research},
    publisher = {INFORMS},
    author = {Gopalan, Ram and Talluri, Kalyan T.},
    year = {1998},
    pages = {260--271},
}

@article{
    rahimi-vahed_2015_fleet,
    author  = {Rahimi-Vahed, Alireza and Crainic, Teodor Gabriel and Gendreau, Michel and Rei, Walter},
    title   = {Fleet-sizing for multi-depot and periodic vehicle routing problems using a modular heuristic algorithm},
    journal = {Computers \& Operations Research},
    volume  = {53},
    pages   = {9--23},
    year    = {2015},
    doi     = {10.1016/j.cor.2014.07.004}
}

@article{
    hemmelmayr_2013_waste,
    author  = {Hemmelmayr, Vera and Doerner, Karl F. and Hartl, Richard F. and Rath, Stefan},
    title   = {A heuristic solution method for node routing based solid waste collection problems},
    journal = {Journal of Heuristics},
    volume  = {19},
    pages   = {129--156},
    year    = {2013},
    doi     = {10.1007/s10732-011-9188-9}
}

@article{pecin_limited_2017,
    title = {Limited memory {Rank}-1 {Cuts} for {Vehicle} {Routing} {Problems}},
    volume = {45},
    issn = {0167-6377},
    url = {https://www.sciencedirect.com/science/article/pii/S0167637717301049},
    doi = {10.1016/j.orl.2017.02.006},
    number = {3},
    urldate = {2023-09-05},
    journal = {Operations Research Letters},
    author = {Pecin, Diego and Pessoa, Artur and Poggi, Marcus and Uchoa, Eduardo and Santos, Haroldo},
    month = may,
    year = {2017},
    pages = {206--209},
}

@misc{eveborn_branch-price-and-cut_2026,
    title = {Branch-{Price}-and-{Cut} {Accelerated} with a {Pricing} for {Integrality} {Heuristic} for the {Electrical} {Vehicle} {Routing} {Problem} with {Time} {Windows} and {Charging} {Time} {Slots}},
    url = {http://arxiv.org/abs/2602.08673},
    doi = {10.48550/arXiv.2602.08673},
    urldate = {2026-03-11},
    publisher = {arXiv},
    author = {Eveborn, Lukas and Rönnberg, Elina},
    month = feb,
    year = {2026},
    note = {arXiv:2602.08673 [math]},
}

@article{jacquillat2022optimizing,
  title={Optimizing relay operations toward sustainable logistics},
  author={Jacquillat, Alex and Schmid, Alexandria and Wang, Kai},
  journal={Available at SSRN 4241031},
  year={2022}
}

@article{jacquillat2026iterative,
  title={An Iterative Network Flow Algorithm for Online Pickup
and Delivery},
  author={Jacquillat, Alex and Luo, Jason},
  journal={Proceedings of the 2026 INFORMS Optimization Society Conference},
  year={2026}
}

@misc{rossum_adaptive_2026,
    title = {Adaptive {Partitioning} in {Column} {Generation} for {Nested} {Paths}},
    url = {http://arxiv.org/abs/2605.28692},
    doi = {10.48550/arXiv.2605.28692},
    urldate = {2026-09-16},
    publisher = {arXiv},
    author = {van Rossum, Bart and van Lieshout, Rolf and Jacquillat, Alexandre},
    month = may,
    year = {2026},
    note = {arXiv:2605.28692 [math.OC]},
}
}

\newpage
\ECSwitch
\ECHead{Subpath-Based Column Generation for Electric Vehicle Routing Problems: Electronic Companion}
\renewcommand{\ECRUNAUTHOR}{}
{
\small
\section{Column generation for the EVRPTWNL}
\label{app.cg}

Algorithm~\ref{alg.CG} solves the linear relaxation of the EVRPTWNL, parameterized by a set of paths $\calP$. This parameterization allows us to express the linear relaxation solution after \textit{ng}-relaxations and subset-row cuts, both which affect the feasible paths $\calP$. 

\begin{algorithm}
\caption{\textsc{ColumnGeneration}$(\calP)$.}\small
\label{alg.CG}
\begin{algorithmic}
\item \textbf{Initialization:} Construct a set of paths $\calP_{0} \subset \calP$ such that $\EVRPLP(\calP_{0})$ is feasible. Initialize $\ell = 0$.

\item Iterate between Steps 1-3, until termination.

\begin{itemize}[left=0pt..1em]
    \item[] \textbf{Step 1.}  Solve $\EVRPLP(\calP{}_{\ell})$, obtaining primal solution $\bz_\ell$ and dual variables $\bkappa$, $\bmu$, and $\bnu$.
    \item[] \textbf{Step 2.}  Solve pricing problem to generate paths $p \in \Pnew \subseteq \calP$ with negative reduced cost (Equation~\eqref{eq.RC}).
    \item[] \textbf{Step 3.} If $|\Pnew|=0$, \texttt{STOP}: return solution $\bz_\ell$. Otherwise, update $\calP_{\ell+1} \gets \calP_{\ell} \cup \Pnew$ and $\ell \gets \ell+1$.
\end{itemize}
\end{algorithmic}
\end{algorithm}

\section{Arc-based optimization formulation for the EVRPTWNL}
\label{app.arc_formulation}

We provide the arc-based integer optimization formulation for the EVRPTWNL. We define $\Depots^\text{start}$ and $\Depots^\text{end}$ as two copies of the set of depots $\Depots$. We define the following decision variables:
\begin{itemize}[noitemsep]
    \item $x_{i,j}^k$: binary variable denoting if arc $(i, j)\in \calA$ is traversed by vehicle $k \in [K]$.
    \item $T^k_i$ denoting the time that vehicle $k \in [K]$ arrives at node $i \in \calV$;
    \item $B^k_i$ denoting the charge that vehicle $k \in [K]$ arrives at node $i \in \calV$ with;
    \item $D^k_i$ denoting the load that vehicle $k \in [K]$ arrives at node $i \in \calV$ with;
    \item $\tau^k_i$ denoting the time spent charging at node $i \in \calV$ for vehicle $k \in [K]$ ($0$ for non-charging stations $i \notin \Stations$)
\end{itemize}

Equation~\eqref{eq.arc.objective} minimizes travel costs, subject to start and end requirements (Equations~\eqref{eq.arc.start_depots_count}--\eqref{eq.arc.end_depots}), flow conservation for each vehicle (Equation~\eqref{eq.arc.flow_conservation}), and customer service requirements (Equation~\eqref{eq.arc.serve_customers}). 
Equations~\eqref{eq.arc.load}--\eqref{eq.arc.load_bounds} track the accumulated load for each vehicle at each node.
Equations~\eqref{eq.arc.time}--\eqref{eq.arc.time_initial} track the arrival time for each vehicle at each node, and Equations~\eqref{eq.arc.time_window_start}--\eqref{eq.arc.time_window_end} impose time window requirements.
Equations~\eqref{eq.arc.charge}--\eqref{eq.arc.charge_bound} track the arrival state-of-charge for each vehicle at each node. Note that Equations~\eqref{eq.arc.time} and~\eqref{eq.arc.charge} are coupled via the charging time decision $\tau^k_i$.

\begin{alignat}{3}
    \label{eq.arc.objective}
    \min \quad
    & 
    \sum_{k \in [K]} \sum_{(i, j) \in \calA} 
    c_{i,j} x_{i,j}^k
    \\
    \text{such that} \quad
    \label{eq.arc.start_depots_count}
    & 
    \sum_{k \in [K]}
    \sum_{j: (i, j) \in \calA} 
    x_{i,j}^k = v^\text{start}_i
    &\quad& 
    \forall \ i \in \Depots^\text{start}
    \\
    \label{eq.arc.end_depots_count}
    & 
    \sum_{k \in [K]}
    \sum_{i: (i, j) \in \calA} 
    x_{i,j}^k \geq v^\text{end}_j
    &\quad& 
    \forall \ j \in \Depots^\text{end}
    \\
    \label{eq.arc.start_depots}
    & \sum_{i \in \Depots^\text{start}} \sum_{j: (i, j) \in \calA} x_{i,j}^k = 1
    &\quad&
    \forall \ k \in [K]
    \\
    \label{eq.arc.end_depots}
    & \sum_{j \in \Depots^\text{end}} \sum_{i: (i, j) \in \calA} x_{i,j}^k = 1
    &\quad&
    \forall \ k \in [K]
    \\
    \label{eq.arc.flow_conservation}
    & \sum_{j: (j, i) \in \calA} x_{j,i}^k 
    = \sum_{j: (i, j) \in \calA} x_{i,j}^k
    &\quad&
    \forall \ k \in [K], 
    \ \forall \ i \in \Custs \cup \Stations
    \\
    \label{eq.arc.serve_customers}
    & \sum_{k \in [K]} \sum_{j: (i, j) \in \calA} x_{i,j}^k = 1
    &\quad&
    \forall \ i \in \Custs 
    \\
    \label{eq.arc.load}
    & D^k_j
    \geq D^k_i + d_i - M \cdot (1 - x_{i,j}^k)
    &\quad&
    \forall \ k \in [K], 
    \ \forall \ (i, j) \in \calA 
    \\
    \label{eq.arc.load_initial}
    & D^k_i = 0
    &\quad&
    \forall \ k \in [K], 
    \ \forall \ i \in \Depots^\text{start}
    \\
    \label{eq.arc.load_bounds}
    & D^k_i \leq D
    &\quad&
    \forall \ k \in [K], 
    \ \forall \ i \in \calV
    \\
    \label{eq.arc.time}
    & T^k_j
    \geq T^k_i + t_{i,j} + \tau^k_i - M \cdot (1 - x_{i,j}^k)
    &\quad&
    \forall \ k \in [K], 
    \ \forall \ (i, j) \in \calA 
    \\
    \label{eq.arc.time_charging}
    & \tau^k_i = 0
    &\quad&
    \forall \ k \in [K], 
    \ \forall \ i \in \calV \setminus \Stations
    \\
    \label{eq.arc.time_initial}
    & T^k_i = 0
    &\quad&
    \forall \ k \in [K], 
    \ \forall \ i \in \Depots^\text{start}
    \\
    \label{eq.arc.time_window_start}
    & T^k_i \geq \alpha_i
    &\quad&
    \forall \ k \in [K], 
    \ \forall \ i \in \calV
    \\
    \label{eq.arc.time_window_end}
    & T^k_i \leq \beta_i
    &\quad&
    \forall \ k \in [K], 
    \ \forall \ i \in \calV
    \\
    \label{eq.arc.charge}
    & B^k_j
    \leq \gC{B^k_i}{\tau^k_i} - b_{i,j} 
    + M \cdot (1 - x_{i,j}^k)
    &\quad&
    \forall \ k \in [K], 
    \ \forall \ (i, j) \in \calA 
    \\
    \label{eq.arc.charge_initial}
    & B^k_i = B
    &\quad&
    \forall \ k \in [K], 
    \ \forall \ i \in \Depots^\text{start}
    \\
    \label{eq.arc.charge_bound}
    & B^k_i \leq B
    &\quad&
    \forall \ k \in [K], 
    \ \forall \ i \in \calV
    \\
    \label{eq.arc.x_binary}
    & x_{i,j}^k \in \{0, 1\}
    &\quad&
    \forall \ k \in [K],
    \ \forall \ (i, j) \in \calA
    \\
    \label{eq.arc.load_nonnegative}
    & D_i^k \geq 0
    &\quad&
    \forall \ k \in [K],
    \ \forall \ i \in \calV
    \\
    \label{eq.arc.time_charging_nonnegative}
    & \tau_i^k \geq 0
    &\quad&
    \forall \ k \in [K],
    \ \forall \ i \in \Stations 
    \\
    \label{eq.arc.time_nonnegative}
    & T_i^k \geq 0
    &\quad&
    \forall \ k \in [K],
    \ \forall \ i \in \calV
    \\
    \label{eq.arc.charge_nonnegative}
    & B_i^k \geq 0
    &\quad&
    \forall \ k \in [K],
    \ \forall \ i \in \calV
\end{alignat}

Note that this formulation cannot model a vehicle visiting the same charging station more than once; this can be mitigated by duplicating charging station nodes.

\subsubsection*{Computational results.}

We consider EVRPTW instance C101 from \cite{schneider2014electric} with 25 customers, and investigate the scalability of the arc-based formulation. We consider the following variants:
\begin{enumerate}[(1),noitemsep]
    \item Equations~\eqref{eq.arc.objective}--\eqref{eq.arc.serve_customers}, \eqref{eq.arc.time}--\eqref{eq.arc.time_initial}, \eqref{eq.arc.x_binary}, \eqref{eq.arc.time_charging_nonnegative}--\eqref{eq.arc.time_nonnegative}, modeling the core VRP;
    \item All previous constraints, and Equations~\eqref{eq.arc.time_window_start}--\eqref{eq.arc.time_window_end}, modeling the VRPTW;
    \item All previous constraints, and Equations~\eqref{eq.arc.charge}--\eqref{eq.arc.charge_bound} and~\eqref{eq.arc.charge_nonnegative}, modeling the EVRPTW;
    \item The full formulation, modeling the capacitated EVRPTW (as in Section~\ref{sec.numerical_results}).
\end{enumerate}

All runs are conducted with a 1-hour time limit with a 1\% optimality tolerance. Table~\ref{tab.arc} shows that even the smallest ($N = 6$) instances have large computational times, and the full problem does not terminate within the one-hour limit, leaving large optimality gap, with as few as 12 customers. The results also show the impact of the electric vehicle routing problem---which requires jointly optimizing routing and charging variables---as compared to traditional vehicle routing problems, resulting in orders of magnitude increases in computational times. Ultimately, these results show that the---capacitated---EVRPTW is a highly challenging optimization problem for which off-the-shelf methods exhibit limited tractability.

\begin{table}[H]
    \small
    \centering
    \caption{Computation times and optimality gaps for arc-based formulations incorporating various problem elements of the EVRPTW and varying the number of customers.}
    \label{tab.arc}
    \begin{tabular}{
        S[table-format=2.0]
        *{4}{S[table-format=>4.2,group-separator={}]}
        *{4}{S[table-format=<2.2,group-separator={},table-align-comparator=false]}
    }
        \toprule
        & \multicolumn{4}{c}{Computation time (s)}
        & \multicolumn{4}{c}{Optimality gap (\%)}
        \\
        \cmidrule(lr){2-5}
        \cmidrule(lr){6-9}
        $N$ & 
        {(1)} & {(2)} & {(3)} & {(4)} &
        {(1)} & {(2)} & {(3)} & {(4)} 
        \\
        \midrule
        6 &
        18.05 & 42.57 & 956.21 & 1217.66 & 
        <1.00 & <1.00 & <1.00 & <1.00 \\ 
        8 & 
        53.41 & 898.54 & 2624.24 & >3600.0 & 
        <1.00 & <1.00 & <1.00 & 4.18 \\ 
        10 &
        20.95 & 306.03 & >3600.0 & 2048.21 & 
        <1.00 & <1.00 & 6.09 & <1.00 \\ 
        12 &
        >3600.0 & 544.96 & >3600.0 & >3600.0 & 
        8.27 & <1.00 & 9.04 & 9.96 \\ 
        15 &
        >3600.0 & >3600.0 & >3600.0 & >3600.0 & 
        15.31 & 8.05 & 14.85 & 14.54 \\ 
        18 & 
        >3600.0 & >3600.0 & >3600.0 & >3600.0 & 
        11.96 & 8.72 & 22.24 & 19.82 \\ 
        20 &
        >3600.0 & >3600.0 & >3600.0 & >3600.0 & 
        22.74 & 9.47 & 23.44 & 24.22 \\ 
        25 &
        >3600.0 & >3600.0 & {--} & {--} &         
        21.35 & 25.81 & {--} & {--} \\ 
        \bottomrule
    \end{tabular}
\end{table}
\section{Proofs in Section~\ref{sec.labelsetting}}
\label{app.proofs_4}

Throughout the appendix, we denote by $\Pos{x} := \max \{0, x\}$ the positive part of $x$, and by $\fracpart{x} = x - \floor{x}$ its fractional part.

\subsection{Section 4.2}
\label{app.proofs_42}

We first provide a characterization of subpath-level timing resources $\TendMin{s}, \TstartMax{s}$. Let $s$ have node sequence $(n_0, \dots, n_m)$. Define cumulative travel times:
\begin{equation}
    A_0 = 0, 
    \quad 
    A_k = \sum_{\ell = 0}^{k-1} t_{n_\ell,n_\ell+1} 
    \quad \forall \ k \in \{1, \dots, m\}
\end{equation}
Then, $(T^\text{E}_0, \dots, T^\text{E}_m)$ and $(T^\text{L}_0, \dots, T^\text{L}_m)$ in Definition~\ref{def.subpath_resources}\ref{def.subpath_resources_TendMin}--\ref{def.subpath_resources_TstartMax} can be expressed by unrolling the forward and backward recursions:
\begin{equation}
    \label{eq.tE_tL}
    \forall \ k \in \{0, \dots, m\}, 
    \quad 
    T^\text{E}_k = \max_{0 \leq i \leq k} \left\{ \alpha_{n_i} + A_k - A_i \right\},
    \quad 
    T^\text{L}_k = \min_{k \leq i \leq m} \left\{ \beta_{n_i} - A_i + A_k  \right\}.
\end{equation}
This allows us to show:
\begin{lemma}
    \label{lemma.subpath_TWs}
    Let $s$ be a partial subpath with node sequence $(n_0, \dots, n_m)$. Define $(T^\text{E}_0, \dots, T^\text{E}_m)$ and $(T^\text{L}_0, \dots, T^\text{L}_m)$ per Definition~\ref{def.subpath_resources}\ref{def.subpath_resources_TendMin}--\ref{def.subpath_resources_TstartMax}. Suppose $s$ is feasible and has service time sequence $(t_0, \dots, t_m)$. Then, $\alpha_{n_k} \leq T^\text{E}_k \leq t_k \leq T^\text{L}_k \leq \beta_{n_k}$ for all $k \in \{0, \dots, m\}$.
\end{lemma}
\proof{Proof of Lemma~\ref{lemma.subpath_TWs}.}
Feasibility of the service time sequence implies, for any $k \in \{0, \dots, m\}$:
\begin{alignat}{3}
    t_k & \geq \alpha_{n_i} + A_k - A_i 
    &\quad& \forall \ i \in \{0, \dots, k\}
    \\
    t_k & \leq \beta_{n_i} - A_i + A_k
    &\quad& \forall \ i \in \{k, \dots, m\}
\end{alignat}
Taking the maximum and minimum respectively over $i$ gives $T^\text{E}_k \leq t_k \leq T^\text{L}_k$. Taking $i = k$ in Equation~\eqref{eq.tE_tL} gives $\alpha_{n_k} \leq T^\text{E}_k$ and $T^\text{L}_k \leq \beta_{n_k}$. \hfill\Halmos

For a node sequence $(n_0, \dots, n_m)$, $(T^\text{E}_0, \dots, T^\text{E}_m)$ and $(T^\text{L}_0, \dots, T^\text{L}_m)$ can be viewed as ``effective time windows'', that is, the earliest and latest times each node can be visited while maintaining time-window feasibility. Feasibility of a node sequence (without charging) would be equivalent to ensuring $T^\text{E}_k \leq T^\text{L}_k$ for each $k$.

\subsubsection*{Proof of Proposition~\ref{prop.subpath_REFs}.}
    
The correctness of extension rules \ref{prop.subpath_REFs_rc}--\ref{prop.subpath_REFs_TendMin} follow directly from the definitions.
We show that rule \ref{prop.subpath_REFs_TstartMax} is consistent with Definition~\ref{def.subpath_resources}\ref{def.subpath_resources_TstartMax}. Let $s = (n_0, \dots, n_m)$ and $a = (n_m, n_{m+1})$. By definition, $\TstartMax{s} = T^L_0$ and $\TstartMax{s \oplus a} = \widehat{T}^L_0$, where:
\begin{align*}
    & T^L_m = \beta_{n_m}, 
        \quad 
        \forall \ \ell \in \{0, \dots, m-1\}, \
        T^L_\ell = \min \left\{ 
            T^L_{\ell+1} - t_{n_{\ell},n_{\ell+1}}, \
            \beta_{n_\ell}
        \right\}\\
    & \widehat{T}^L_{m+1} = \beta_{n_{m+1}}, 
        \quad 
        \forall \ \ell \in \{0, \dots, m\}, \
        \widehat{T}^L_\ell = \min \left\{ 
            \widehat{T}^L_{\ell+1} - t_{n_{\ell},n_{\ell+1}}, \
            \beta_{n_\ell}
        \right\}
\end{align*}

If there exists some $\ell \in \{0, \dots, m\}$ for which $T^L_{\ell} = \widehat{T}^L_{\ell}$, this must hold for all $0 \leq k \leq \ell$ per the recursive definition. This means that $\widehat{T}^L_0 = T^L_0$, hence $\TstartMax{s \oplus a}=\TstartMax{s}$. Otherwise, $\widehat{T}^L_{\ell} < T^L_\ell$ for all $\ell \in \{0, \dots, m\}$. This implies:
\begin{equation*}
    \beta_{n_{m+1}} - \sum_{\ell=k}^{m} t_{n_\ell,n_{\ell+1}} < \beta_{n_k}, 
    \quad
    \forall \ k \in \{0, \dots, m\}
\end{equation*}

These two cases can be captured succinctly as follows, which shows rule \ref{prop.subpath_REFs_TstartMax}.
\begin{equation*}
    \TstartMax{s \oplus a}
    = \min \Big\{ 
        \beta_{n_{m+1}} - \sum_{\ell=0}^{m} t_{n_\ell,n_{\ell+1}}, \
        \TstartMax{s}
    \Big\}
    = \min \left\{ 
        \beta_{n_{m+1}} - \Ttravel{s \oplus a},
        \TstartMax{s}
    \right\}
\end{equation*}

We next show the criteria for $s \oplus a$ to be feasible. Criteria~\ref{prop.subpath_REFs_feasible_load}--\ref{prop.subpath_REFs_feasible_B} are obvious. Both parts of  \ref{prop.subpath_REFs_feasible_TW} are equivalent due to symmetry. The first requires that starting at node $n_m$ at $\TendMin{s}$, one must reach $n_{m+1}$ before $\beta_{n_{m+1}}$; equivalently, the second requires that reaching node $n_m$ at $\beta_{n_{m+1}}$ must not force one to start subpath $s \oplus a$ before $\alpha_{n_0}$, the start of the time window at the start of $s$.
\hfill\Halmos

\subsubsection*{{Proof of Lemma~\ref{lemma.subpath_fs}.}}
Suppose that one starts subpath $s$ at time $t$. Then by unrolling Equation~\eqref{eq.propagation_time}:
\begin{equation}
    f_s(t) = \max \{ \TendMin{s}, t + \Ttravel{s} \}
    \quad (t \in [0, \TstartMax{s}])
\end{equation}
The two cases can be explained as follows: if some node $n$ is reached before or at the start of its time window $\alpha_n$, then the subpath must arrive at each subsequent node at the earliest possible time, and therefore terminate at time $\TendMin{s}$. Otherwise no waiting is incurred and $s$ terminates at time $t + \Ttravel{s}$.
The domain of $f_s$ is $[0, \TstartMax{s}]$ by definition. 
From the definition of $Q(s)$, $f_s(t)$ is equivalently equal to $\max \{ \TendMin{s}, t + Q(s) \}$:
\begin{alignat}{3}
    \max \{ \TendMin{s}, \ t + Q(s) \}
    & = \max \{ 
        \TendMin{s}, 
        \ t + \Ttravel{s}, 
        \ \TendMin{s} + t - \TstartMax{s} 
    \}
    \\
    & = \max \{ \TendMin{s}, \ t + \Ttravel{s} \}
    \qquad \text{(since $t \in [0, \TstartMax{s}]$)}
    \qquad \Halmos
\end{alignat}
\begin{remark}
    We also note the helpful identity $P(s) + Q(s) = \TendMin{s}$.
\end{remark}
\begin{remark}
    This lemma is also from \cite{irnich_resource_2008} (see Proposition 6); however, as far as we are aware we are the first to use this result to directly define label-setting resources.
\end{remark}

\subsubsection*{{Proof of Proposition~\ref{prop.subpath_fs_domination}.}}

We show that $s_1, s_2$ satisfy Definition~\ref{def.subpath_domination}\ref{def.subpath_domination_TendMin}--\ref{def.subpath_domination_Q},if and only if $\text{dom}(f_{s_2}) \subseteq \text{dom}(f_{s_1})$ and, for any $t \in \text{dom}(f_{s_2})$, we have $f_{s_1}(t) \leq f_{s_2}(t)$.
\begin{itemize}
    \item[$(\Rightarrow)$] 
    Suppose $s_1, s_2$ satisfies Definition~\ref{def.subpath_domination}\ref{def.subpath_domination_TendMin}--\ref{def.subpath_domination_Q}. 
    Then $\TendMin{s_1} \leq \TendMin{s_2}$, $\TstartMax{s_1} \geq \TstartMax{s_2}$, and $Q(s_1) \leq Q(s_2)$, so Equation~\eqref{eq.subpath_fs} gives $\text{dom}(f_{s_2}) \subseteq \text{dom}(f_{s_1})$ and $f_{s_1} \leq f_{s_2}$ on the domain. 
    \item[$(\Leftarrow)$]
    Suppose the latter. Then, $\text{dom}(f_{s_2}) \subseteq \text{dom}(f_{s_1})$ implies \ref{def.subpath_domination_TstartMax} $\TstartMax{s_1} \geq \TstartMax{s_2}$. Also, $f_{s_1}(0) \leq f_{s_2}(0)$ implies \ref{def.subpath_domination_TendMin} $\TendMin{s_1} \leq \TendMin{s_2}$. Finally:
    \begin{alignat}{3}
        Q(s_1) 
        & \leq f_{s_1}(\TstartMax{s_2}) - \TstartMax{s_2}
        &\quad& \text{(evaluating $f_{s_1}$ at $\TstartMax{s_2}$)}
        \\
        & \leq f_{s_2}(\TstartMax{s_2}) - \TstartMax{s_2}
        = Q(s_2).
        &\quad& \text{(since $f_{s_1} \leq f_{s_2}$)}
        \qquad \Halmos
    \end{alignat}
\end{itemize}
\subsection{Section 4.3}
\label{app.proofs_43}

\subsubsection*{Nonlinear charging preliminaries}

We note the following lemma (Proposition 2.1 of \cite{lee_exact_2021}):
\begin{lemma}
\label{lemma.gT}
If $0 \leq b_1 \leq b_2$ and $0 \leq b \leq B - b_2$, then $\gT{b_1}{b_1 + b} \leq \gT{b_2}{b_2 + b}$.
\end{lemma}

\subsubsection*{Path domination}

We first observe that the nonlinear comparison of path resources for domination can be restated as two componentwise comparisons: 
\begin{lemma}
\label{lemma.path_domination_componentwise}
    Let $p_1, p_2$ be two feasible partial paths starting and ending at the same nodes. Define $\Teff{p} := \TendMin{p} - \gT{0}{\BendMax{p}}$. Then Definition~\ref{def.path_domination}\ref{def.path_domination_BendMax} is equivalent to:
    \begin{enumerate*}[(a)]
        \item 
        \label{lemma.path_domination_componentwise_T}
        $\TendMin{p_1} \leq \TendMin{p_2}$; and 
        \item 
        \label{lemma.path_domination_componentwise_Teff}
        $\Teff{p_1} \leq \Teff{p_2}$.
    \end{enumerate*}
\end{lemma}
\begin{proof}{Proof of Lemma~\ref{lemma.path_domination_componentwise}.}
\begin{itemize}
    \item[$(\Rightarrow)$] Assume Definition~\ref{def.path_domination}\ref{def.path_domination_BendMax}. 
    Then:
    \begin{alignat}{3}
        &\quad& 
        \TendMin{p_1} +  \gT{ \BendMax{p_1} }{ \max \left\{ \BendMax{p_1}, \BendMax{p_2} \right\} } 
        & \leq \TendMin{p_2}
        \\
        \iff 
        &\quad& 
        \TendMin{p_1} + \ginv{ \max \left\{ \BendMax{p_1}, \BendMax{p_2} \right\} } - \ginv{\BendMax{p_1}}
        & \leq \TendMin{p_2}
        \\
        \iff 
        &\quad& 
        \TendMin{p_1} + \Pos{ \ginv{\BendMax{p_2}} - \ginv{\BendMax{p_1}} }
        & \leq \TendMin{p_2}
    \end{alignat}
    \ref{lemma.path_domination_componentwise_T} and \ref{lemma.path_domination_componentwise_Teff} follow from the fact that the second term is greater than or equal to $0$ and $\ginv{\BendMax{p_2}} - \ginv{\BendMax{p_1}}$ respectively.
    \item[$(\Leftarrow)$] Assume \ref{lemma.path_domination_componentwise_T} and \ref{lemma.path_domination_componentwise_Teff}. Then $\TendMin{p_1} + \Pos{ \ginv{\BendMax{p_2}} - \ginv{\BendMax{p_1}} }
    \leq \TendMin{p_2}$, from which Definition~\ref{def.path_domination}\ref{def.path_domination_BendMax} holds.
    \hfill\Halmos
\end{itemize}
\end{proof}

To aid in the proofs that follow, we generalize $\Teff{\cdot}$ to time-charge pairs that might not be equal to $(\TendMin{p}, \BendMax{p})$ (such that $\Teff{p} = \eta(\TendMin{p}, \BendMax{p})$):
\begin{definition}
    Define $\eta(t, b)$ as the charge-adjusted time:
    $\eta(t, b) := t - \ginv{b}$.
\end{definition}

We show the following lemma:
\begin{lemma}
    \label{lemma.chargeseq_onestep}
    Let $(t_1, b_1)$ and $(t_2, b_2)$ be two time-charge pairs arriving at the start of subpath $s$. Assume: \begin{enumerate*}[(a)]
        \item 
        $t_1 \leq t_2$; and 
        \item 
        $\eta(t_1, b_1) \leq \eta(t_2, b_2)$.
    \end{enumerate*}
    Suppose charging for $\tau_2 \geq 0$ at $(t_2, b_2)$ allows feasible service of subpath $s$, i.e.
    $\gC{b_2}{\tau_2} \geq B(s)$ and $t_2 + \tau_2 \leq \TstartMax{s}$. Define:
    \begin{equation}
        \tau_1 := \max \left\{
            0, \ 
            \ginv{B(s)} - \ginv{b_1}, \ 
            P(s) - t_1
        \right\}
    \end{equation}
    Then, charging for $\tau_1$ at $(t_1, b_1)$ allows feasible service of subpath $s$. The resulting states at the end of $s$, given by $(t'_1, b'_1)$ and $(t'_2, b'_2)$, satisfy \begin{enumerate*}[(a)]
        \item $t'_1 \leq t'_2$; and
        \item $\eta(t'_1, b'_1) \leq \eta(t'_2, b'_2)$. 
    \end{enumerate*}
    \begin{alignat}{3}
        t'_1 & := f_s(t_1 + \tau_1),
        &\quad
        b'_1 & := \gC{b_1}{\tau_1} - B(s)
        \\
        t'_2 & := f_s(t_2 + \tau_2),
        &\quad 
        b'_2 & := \gC{b_2}{\tau_2} - B(s)
        &\quad
    \end{alignat}
\end{lemma}
\begin{proof}{Proof of Lemma~\ref{lemma.chargeseq_onestep}.} First, charging feasibility: by definition of $\tau_1$ and the fact that $g$ is nondecreasing,
\begin{equation}
    \gC{b_1}{\tau_1}
    = g\left(\ginv{b_1} + \tau_1 \right)
    \geq g\left( \ginv{B(s)} \right) = B(s)
\end{equation}

Second, time window feasibility. We wish to show
\begin{alignat}{3}
    t_1 + \tau_1 & = \max \left\{
        t_1, \ 
        \eta(t_1, b_1) + \ginv{B(s)}, \ 
        P(s)
    \right\} 
    \leq \TstartMax{s}
\end{alignat}
The first is true since $t_1 \leq t_2 \leq t_2 + \tau_2 \leq \TstartMax{s}$. The third is true since $P(s) := \min \left\{ \TstartMax{s}, \TendMin{s} - \Ttravel{s} \right\} \leq \TstartMax{s}$ by definition. To show the second:
\begin{alignat}{3}
    \eta(t_1, b_1) + \ginv{B(s)}
    & \leq \eta(t_2, b_2) + \ginv{B(s)}
    \\
    & = t_2 + \ginv{B(s)} - \ginv{b_2}
    &\quad& 
    \text{(because $\gC{b_2}{\tau_2} \geq B(s)$)}
    \\
    & \leq t_2 + \tau_2 \leq \TstartMax{s}
\end{alignat}
Indeed, from the above, we can conclude that while $t_1 + \tau_1$ might be greater than $t_2 + \tau_2$, it is not greater than $\max \{ t_2 + \tau_2, P(s) \}$.

Third, $t'_1 \leq t'_2$. Since $t_1 + \tau_1 \geq P(s)$, $f_s(t_1 + \tau_1) = t_1 + \tau_1 + Q(s)$. Then:
\begin{equation}
    t'_1 
    = f_s(t_1 + \tau_1)
    = t_1 + \tau_1 + Q(s)
    \leq \max \left\{ t_2 + \tau_2, P(s) \right\} + Q(s)
    \overset{\text{(i)}}{=} f_s( \max \left\{ t_2 + \tau_2, P(s) \right\} )
    \overset{\text{(ii)}}{=} f_s(t_2 + \tau_2) = t'_2
\end{equation}
where (i) is because $\max \left\{ t_2 + \tau_2, P(s) \right\} \geq P(s)$ and (ii) is by definition of $f_s$ and the fact that $P(s) + Q(s) = \TendMin{s}$.

Fourth, $\eta(t'_1, b'_1) \leq \eta(t'_2, b'_2)$. We define: $r(x) := \ginv{g(x) - B(s)}$ for $x \geq \ginv{B(s)}$. Since $g$ is nondecreasing and concave, $\ginv{\cdot}$ is nondecreasing and convex, and for $\ginv{B(s)} 
\leq x \leq y$:
\begin{alignat}{3}
    0 \leq r(y) - r(x) 
    & = \ginv{g(y) - B(s)} - \ginv{g(x) - B(s)}
    \\
    & \leq \ginv{g(y)} - \ginv{g(x)} 
    &\quad& \text{(by Lemma~\ref{lemma.gT})}
    \\
    & \leq y - x
\end{alignat}
Therefore our defined function $r$ is nondecreasing and 1-Lipschitz. Therefore, $\Phi(d, \eta) := d + Q(s) - r(d - \eta)$ is nondecreasing in each argument on its domain $d - \eta \geq \ginv{B(s)}$. We therefore have:
\begin{alignat}{3}
    \eta(t'_1, b'_1)
    & = t_1 + \tau_1 + Q(s) - \ginv{b'_1}
    \\
    & = t_1 + \tau_1 + Q(s) - \ginv{g\left( \ginv{b_1} + \tau_1 \right) - B(s)}
    \\
    & = \Phi(t_1 + \tau_1, \ \eta(t_1, b_1))
    \qquad \text{(since $t_1 + \tau_1 - \eta(t_1, b_1) = \ginv{b_1} + \tau_1$)}
    \\
    & \leq \Phi( \max \left\{ t_2 + \tau_2, P(s) \right\}, \ \eta(t_1, b_1))
    \\
    & \leq \Phi( \max \left\{ t_2 + \tau_2, P(s) \right\}, \ \eta(t_2, b_2))
    \qquad 
    \text{(since $\Phi$ nondecreasing)}
    \\
    & = \max \left\{ t_2 + \tau_2, P(s) \right\} + Q(s) - r\left( \max \left\{ t_2 + \tau_2, P(s) \right\} - \eta(t_2, b_2) \right) 
    \\
    & \leq \max \left\{ t_2 + \tau_2, P(s) \right\} + Q(s) - r\left( t_2 + \tau_2 - \eta(t_2, b_2) \right) 
    \qquad 
    \text{(since $r$ nondecreasing)}
    \\
    & = t'_2 - \ginv{g(\ginv{b_2} + \tau_2) - B(s)}
    = \eta(t'_2, b'_2).
    \hfill\Halmos
\end{alignat}
\end{proof}

The statement of Lemma~\ref{lemma.chargeseq_onestep} serves two purposes. On one hand, it justifies the definition of $\tau$ in Equation~\eqref{eq.chargeseq_tau_alg} as dominating any other choice of $\tau$, in terms of the two conditions in Lemma~\ref{lemma.path_domination_componentwise}---by setting $(t_1, b_1) = (t_2, b_2)$. On the other hand, when $t_1 \leq t_2$ and $\eta(t_1, b_1) \leq \eta(t_2, b_2)$, it shows that this relationship propagates along subpath extensions, using the greedy charging time rule for $(t_1, b_1)$ and \textit{any} charging time rule for $(t_2, b_2)$. This allows us to prove Proposition~\ref{prop.chargeseq} (and subsequently Theorem~\ref{thm.EVRPTWNL_correctness}).

\subsubsection*{Proof of Proposition~\ref{prop.chargeseq}.}

We prove the following by induction. For feasible partial path $p$ with subpath sequence $(s_1, \dots, s_m)$, let $(\Tend{k}, \Bend{k})$ be generated by Equations~\eqref{eq.chargeseq_tau_tilde_alg}--\eqref{eq.chargeseq_Bend_alg} after subpath $s_k$, and let $(\hat{\Tend{k}}, \hat{\Bend{k}})$ denote any other feasible charging schedule up to $s_k$. We show by induction that for all $k \in \{1, \dots, m\}$:
\begin{equation}
    \label{eq.prop.chargeseq_induct}
    \Tend{k} \leq \hat{\Tend{k}},
    \quad 
    \text{ and }
    \quad
    \eta(\Tend{k}, \Bend{k}) 
    \leq \eta(\hat{\Tend{k}}, \hat{\Bend{k}})
\end{equation}

The base case $k = 1$ is true, since $\hat{\Tend{1}} = \Tend{1} = f_{s_1}(0)$ and $\hat{\Bend{1}} = \Bend{1} = B - B(s_1)$.

Suppose Equation~\eqref{eq.prop.chargeseq_induct} holds for $k$. Let $(\hat{\Tend{k}}, \hat{\Bend{k}})$ be the end time and charge of some feasible charging schedule after subpath $s_k$, and $\hat{\tau_k}$ be the amount charged before subpath $s_{k+1}$. By the induction hypothesis, $\Tend{k} \leq \hat{\Tend{k}}$ and $\Tend{k} - \ginv{\Bend{k}} \leq \hat{\Tend{k}} - \ginv{\hat{\Bend{k}}}$. Applying Lemma~\ref{lemma.chargeseq_onestep} with $(t_1, b_1) = (\Tend{k}, \Bend{k})$, $(t_2, b_2) = (\hat{\Tend{k}}, \hat{\Bend{k}})$, $\tau_2 = \hat{\tau_k}$ and $s = s_{k+1}$ gives:
\begin{equation}
    \tau_k = \max \left\{ 
        0, \ 
        \ginv{B(s_{k+1})} - \ginv{\Bend{k}}, \ 
        P(s_{k+1}) - \Tend{k}
    \right\}
\end{equation}
This coincides with $\tau_k$ defined in Equations~\eqref{eq.chargeseq_tau_tilde_alg}--\eqref{eq.chargeseq_tau_alg}. The lemma therefore shows $\Tend{k+1} \leq \hat{\Tend{k+1}}$ and $\eta(\Tend{k+1}, \Bend{k+1}) \leq \eta(\hat{\Tend{k+1}}, \hat{\Bend{k+1}})$ and the induction is complete.

Finally, from the induction, we get:
\begin{enumerate}[(a)]
    \item $\Tend{k}$ is optimal for $\text{CS}^T(s_1, \dots, s_k)$ (from the first part of Equation~\eqref{eq.prop.chargeseq_induct});
    \item $\Bend{k}$ is optimal for $\text{CS}^B(s_1, \dots, s_k)$, since for any feasible solution with $\hat{\Tend{k}} = \Tend{k}$, the second part of Equation~\eqref{eq.prop.chargeseq_induct} gives $\ginv{\Bend{k}} \geq \ginv{\hat{\Bend{k}}} \implies \Bend{k} \geq \hat{\Bend{k}}$.
\end{enumerate} 
This completes the proof. \hfill\Halmos

\subsubsection*{Proof of Proposition~\ref{prop.path_REFs}.}

Extension rules \ref{prop.path_REFs_rc} and \ref{prop.path_REFs_load} are obvious. 
Extension rules \ref{prop.path_REFs_TendMin} and \ref{prop.path_REFs_BendMax} are shown in Proposition~\ref{prop.chargeseq}. 

Suppose feasibility criteria \ref{prop.path_REFs_feasible_load}--\ref{prop.path_REFs_feasible_T}.
\ref{prop.path_REFs_feasible_load} implies load feasibility of $p \oplus s$ (Definition~\ref{def.path}\ref{def.path_feasible_load}). \ref{prop.path_REFs_feasible_TstartMax} gives $\TendMin{p} + \tau \leq \TstartMax{s}$, so Lemma~\ref{lemma.subpath_fs} and the feasibility of subpath $s$ give the feasibility of all intermediate service times (Definition~\ref{def.path}\ref{def.path_feasible_time}). Also, by construction, charging for $\tau$ gives sufficient charge to perform subpath $s$ thereafter requiring charge $B(s)$:
\begin{alignat}{3}
    \gC{\BendMax{p}}{\tau}
    & \geq \gC{\BendMax{p}}{\tilde{\tau}}
    = \max \left\{ \BendMax{p}, \gC{\BendMax{p}}{\gT{\BendMax{p}}{B(s)}}\right\}
    \\
    & = \max \left\{ \BendMax{p}, B(s)
    \right\} \in [B(s), B]
    \quad \text{(since $\gC{b}{\gT{b}{b'}} = b'$ for $b, b' \in [0, B]$)}
\end{alignat}
Hence all intermediate battery levels are feasible (Definition~\ref{def.path}\ref{def.path_feasible_charge}). This shows feasibility of $p \oplus s$. \hfill\Halmos
\subsection{Section 4.4}
\label{app.proofs_44}

\subsection*{{Proof of Theorem~\ref{thm.EVRPTWNL_correctness}.}}

\subsubsection*{Property~\ref{property.domination_ours}\ref{property.domination_ssa}.}

Let $s_1 \succeq s_2$ be feasible partial subpaths, $a = (n, n')$ be a common arc extension, and suppose $s_2 \oplus a$ is feasible. 
We verify criteria \ref{def.subpath_domination_rc}--\ref{def.subpath_domination_Q} of Definition~\ref{def.subpath_domination} for $s_1 \oplus a \succeq s_2 \oplus a$:
\begin{enumerate}[(i),noitemsep]
    \item[(i)--(iii)] These are trivially true by Proposition~\ref{prop.subpath_REFs}\ref{prop.subpath_REFs_rc}--\ref{prop.subpath_REFs_B} and Definition~\ref{def.subpath_domination}\ref{def.subpath_domination_rc}--\ref{def.subpath_domination_B}.

    \item[(iv)--(vi)] We prove that $s_1 \oplus a$ and $s_2 \oplus a$ satisfy the equivalent Proposition~\ref{prop.subpath_fs_domination}.
    We have the following expressions and domains of $f_{s_1 \oplus a}$ and $f_{s_2 \oplus a}$:
    \begin{alignat}{3}
        \text{dom}(f_{s_1 \oplus a})
        & = \Set{ u \in \text{dom}(f_{s_1}) 
            | 
            f_{s_1}(u) + t_{n,n'} \leq \beta_{n'}
        };
        \\
        \forall \ u \in \text{dom}(f_{s_1 \oplus a}), 
        \quad 
        f_{s_1 \oplus a}(u) & = \max \{ f_{s_1}(u) + t_{n,n'}, \alpha_{n'} \}
        \\
        \text{dom}(f_{s_2 \oplus a})
        & = \Set{ u \in \text{dom}(f_{s_2}) 
            | 
            f_{s_2}(u) + t_{n,n'} \leq \beta_{n'}
        };
        \\
        \forall \ u \in \text{dom}(f_{s_2 \oplus a}), 
        \quad 
        f_{s_2 \oplus a}(u) & = \max \{ f_{s_2}(u) + t_{n,n'}, \alpha_{n'} \}
    \end{alignat}
    Since $s_1 \succeq s_2$, by Proposition~\ref{prop.subpath_fs_domination}, $\text{dom}(f_{s_2}) = [0, \TstartMax{s_2}] \subseteq \text{dom}(f_{s_1}) = [0, \TstartMax{s_1}]$ and $f_{s_1}(u) \leq f_{s_2}(u)$ for all $u \in \text{dom}(f_{s_2})$.
    These two facts together show $\text{dom}(f_{s_2 \oplus a}) \subseteq \text{dom}(f_{s_1 \oplus a})$. 
    Additionally, since the function $x \mapsto \max \{ x + t_{n,n'}, \alpha_{n'} \}$ is increasing, and $f_{s_1}(u) \leq f_{s_2}(u)$, we have $f_{s_1 \oplus a}(u) \leq f_{s_2 \oplus a}(u)$ too for all $u \in \text{dom}(f_{s_2 \oplus a})$.
    Therefore $s_1 \oplus a$ and $s_2 \oplus a$ satisfy Proposition~\ref{prop.subpath_fs_domination}.
\end{enumerate}

Finally, feasibility criteria \ref{prop.subpath_REFs_feasible_load}, \ref{prop.subpath_REFs_feasible_B}, \ref{prop.subpath_REFs_feasible_TW} of Proposition~\ref{prop.subpath_REFs} for $s_1 \oplus a$ are implied by Definition~\ref{def.subpath_domination}\ref{def.subpath_domination_load}, \ref{def.subpath_domination_B}, \ref{def.subpath_domination_TendMin}, and \ref{def.subpath_domination_TstartMax} respectively for $s_1 \oplus a \succeq s_2 \oplus a$.

\subsubsection*{Property~\ref{property.domination_ours}\ref{property.domination_pps}.}

Let $p_1 \succeq p_2$ be feasible partial paths, $s$ be a common subpath extension, and suppose $p_2 \oplus s$ is feasible. 
We aim to show that $p_1 \oplus s$ is also feasible (according to Proposition~\ref{prop.path_REFs}\ref{prop.path_REFs_feasible_load}--\ref{prop.path_REFs_feasible_T}) and $p_1 \oplus s \succeq p_2 \oplus s$  (according to Definition~\ref{def.path_domination}\ref{def.path_domination_rc}--\ref{def.path_domination_BendMax}).

\paragraph{Reduced cost.} Definition~\ref{def.path_domination}\ref{def.path_domination_rc} holds for $p_1 \oplus s \succeq p_2 \oplus s$, based on the additive update rule  (Proposition~\ref{prop.path_REFs}\ref{prop.path_REFs_rc}) and Definition~\ref{def.path_domination}\ref{def.path_domination_rc} for $p_1 \succeq p_2$.

\paragraph{Load.} Proposition~\ref{prop.path_REFs}\ref{prop.path_REFs_feasible_load} for $p_1 \oplus s$ and Definition~\ref{def.path_domination}\ref{def.path_domination_load} for $p_1 \oplus s \succeq p_2 \oplus s$ are shown by the additive update rule
(Proposition~\ref{prop.path_REFs}\ref{prop.path_REFs_load}) and Definition~\ref{def.path_domination}\ref{def.path_domination_load} for $p_1 \succeq p_2$.

\paragraph{Time and charge.}
By Lemma~\ref{lemma.path_domination_componentwise} on $p_1 \succeq p_2$, we have:
\begin{equation}
    \TendMin{p_1} \leq \TendMin{p_2},
    \quad 
    \TendMin{p_1} - \ginv{\BendMax{p_1}}
    \leq \TendMin{p_2} - \ginv{\BendMax{p_2}} 
\end{equation}
Apply Lemma~\ref{lemma.chargeseq_onestep} with $(t_1, b_1) = (\TendMin{p_1}, \BendMax{p_1})$, $(t_2, b_2) = (\TendMin{p_2}, \BendMax{p_2})$, 
and $\tau_1, \tau_2$ defined by the greedy charging rule on $p_1, p_2$ respectively (Equations~\eqref{eq.chargeseq_tau_tilde_alg}--\eqref{eq.chargeseq_tau_alg}). This implies that $p_1 \oplus s$ is time-charge feasible (Proposition~\ref{prop.path_REFs}\ref{prop.path_REFs_feasible_TstartMax}--\ref{prop.path_REFs_feasible_T}), and satisfies:
\begin{equation}
    \TendMin{p_1 \oplus s} \leq \TendMin{p_2 \oplus s},
    \quad 
    \TendMin{p_1 \oplus s} - \ginv{\BendMax{p_1 \oplus s}}
    \leq \TendMin{p_2 \oplus s} - \ginv{\BendMax{p_2 \oplus s}} 
\end{equation}
Therefore Lemma~\ref{lemma.path_domination_componentwise} is satisfied for $p_1 \oplus s \succeq p_2 \oplus s$, which is equivalent to Definition~\ref{def.path_domination}\ref{def.path_domination_BendMax}. 

\subsubsection*{Property~\ref{property.domination_ours}\ref{property.domination_pss}.}

Let $p$ be a feasible partial path, and $s_1, s_2$ be feasible subpaths both extending $p$ such that $s_1 \succeq s_2$. Suppose $p \oplus s_2$ is feasible.
We aim to show that $p \oplus s_1$ is also feasible (according to Proposition~\ref{prop.path_REFs}\ref{prop.path_REFs_feasible_load}--\ref{prop.path_REFs_feasible_T}) and $p \oplus s_1 \succeq p \oplus s_2$  (according to Definition~\ref{def.path_domination}\ref{def.path_domination_rc}--\ref{def.path_domination_BendMax}).

\paragraph{Reduced cost.} Definition~\ref{def.path_domination}\ref{def.path_domination_rc} holds for $p \oplus s_1 \succeq p \oplus s_2$, based on the additive update rule (Proposition~\ref{prop.path_REFs}\ref{prop.path_REFs_rc}) and Definition~\ref{def.subpath_domination}\ref{def.subpath_domination_rc} for $s_1 \succeq s_2$.

\paragraph{Load.} Proposition~\ref{prop.path_REFs}\ref{prop.path_REFs_feasible_load} for $p \oplus s_1$ and Definition~\ref{def.path_domination}\ref{def.path_domination_load} for $p \oplus s_1 \succeq p \oplus s_2$ are shown by the additive update rule
(Proposition~\ref{prop.path_REFs}\ref{prop.path_REFs_load}) and Definition~\ref{def.subpath_domination}\ref{def.subpath_domination_load} for $s_1 \succeq s_2$.

\paragraph{Time and charge.}
Define $\tau_1$ (respectively $\tau_2$) the charging times determined in Proposition~\ref{prop.path_REFs} when extending $p$ along subpath $s_1$ and $s_2$ respectively. We define the following candidate partial paths:
\begin{itemize}[noitemsep]
    \item $p'_1 = p \oplus s_1$, which charges for $\tau_1$ before performing $s_1$;
    \item $\hat{p}$, which performs the same charging in $p$, then charges for $\tau_2$ before performing $s_1$;
    \item $p'_2 = p \oplus s_2$, which charges for $\tau_2$ before performing $s_2$;
\end{itemize}

By Lemma~\ref{lemma.path_domination_componentwise}, the goal is to show:
\begin{equation}
    \TendMin{p'_1} \leq \TendMin{p'_2},
    \quad 
    \TendMin{p'_1} - \ginv{\BendMax{p'_1}}
    \leq \TendMin{p'_2} - \ginv{\BendMax{p'_2}} 
\end{equation}

We first compare $\hat{p}$ to $p'_2$. Since $s_1 \succeq s_2$, by Definition~\ref{def.subpath_domination}:
\begin{alignat}{3}
    & \TendMin{p} + \tau_2 \in \text{dom}(f_{s_2}) \subseteq \text{dom}(f_{s_1}),
    \\
    & f_{s_1}(\TendMin{p} + \tau_2) 
    \leq f_{s_2}(\TendMin{p} + \tau_2)
    = \TendMin{p'_2};
    &\quad& \text{(Proposition~\ref{prop.subpath_fs_domination})}
    \\
    & B(s_1) \leq B(s_2) \leq \gC{\BendMax{p}}{\tau_2}.
    &\quad& \text{(Definition~\ref{def.subpath_domination}\ref{def.subpath_domination_B})}
\end{alignat}
Therefore $\hat{p}$ is time-charge feasible  (Proposition~\ref{prop.path_REFs}\ref{prop.path_REFs_feasible_TstartMax}--\ref{prop.path_REFs_feasible_T}). It completes $s_1$ at time-charge pair $\displaystyle (\hat{t}, \hat{b}) := \left( f_{s_1}(\TendMin{p} + \tau_2), \ \gC{\BendMax{p}}{\tau_2} - B(s_1) \right)$. The above yields:
\begin{equation}
    \label{eq.prop2iii_step1}
    \hat{t} \leq \TendMin{p'_2},
    \quad 
    \hat{b} \geq \BendMax{p'_2}
    \implies 
    \eta(\hat{t}, \hat{b}) \leq \eta(\TendMin{p'_2}, \BendMax{p'_2}).
\end{equation}
We next compare $p'_1$ to $\hat{p}$. By Lemma~\ref{lemma.chargeseq_onestep}, using identical incoming states $(\TendMin{p}, \BendMax{p})$ and the benchmark charging time $\tau_2$ against the greedy charging time $\tau_1$, we obtain that $p'_1$ is time-charge feasible and satisfies:
\begin{equation}
    \label{eq.prop2iii_step2}
    \TendMin{p'_1} \leq \hat{t},
    \quad 
    \eta(\TendMin{p'_1}, \BendMax{p'_1}) \leq \eta(\hat{t}, \hat{b}) 
\end{equation}
Putting Equations~\eqref{eq.prop2iii_step1}--\eqref{eq.prop2iii_step2} together yields 
\begin{equation}
    \TendMin{p'_1} \leq \TendMin{p'_2},
    \quad 
    \eta(\TendMin{p'_1}, \BendMax{p'_1}) \leq \eta(\TendMin{p'_2}, \BendMax{p'_2}).
\end{equation}
proving Lemma~\ref{lemma.path_domination_componentwise} for $p'_1 \succeq p'_2$ which is equivalent to Definition~\ref{def.path_domination}\ref{def.path_domination_BendMax}. Also, the time-charge feasibility of $p'_1$ establishes Proposition~\ref{prop.path_REFs}\ref{prop.path_REFs_feasible_TstartMax}--\ref{prop.path_REFs_feasible_T}. 
This completes the proof.
\hfill\Halmos

\subsection*{{Proof of Theorem~\ref{thm.exact_PP}.}}

We make use of the following lemmas, using Bellman's principle of optimality. In particular, Lemma~\ref{lemma.path_truncate_nondom} underscores the need for Property~\ref{property.domination_ours}\ref{property.domination_pss} in our two-phase label-setting algorithm.
\begin{lemma}
    \label{lemma.subpath_truncate_nondom}
    Under Property~\ref{property.domination_ours}, 
    if $s = s' \oplus a$ is a non-dominated partial subpath
    then $s'$ is non-dominated.
\end{lemma}
\proof{Proof of Lemma~\ref{lemma.subpath_truncate_nondom}.}
    If $s'$ is infeasible, then $s$ is infeasible too. 
    If $s'$ is feasible but dominated by $\overline{s}$, then by Property~\ref{property.domination_ours}\ref{property.domination_ssa} and the feasibility of $s' \oplus a$, $\overline{s} \oplus a$ is also feasible and dominates $s' \oplus a$. Therefore $s'$ is non-dominated. 
    \hfill\Halmos

\begin{lemma}
    \label{lemma.path_truncate_nondom}
    Under Property~\ref{property.domination_ours}, 
    if $p = p' \oplus s$ is a non-dominated partial path, then both $p'$ and $s$ are non-dominated.
\end{lemma}
\proof{Proof of Lemma~\ref{lemma.path_truncate_nondom}.}
    If $p'$ or $s$ are infeasible, then $p$ is infeasible too. 
    If $s$ is feasible but dominated by $\overline{s}$, then by Property~\ref{property.domination_ours}\ref{property.domination_pss} and the feasibility of $p' \oplus s$, $p' \oplus \overline{s}$ is also feasible and dominates $p' \oplus s$. 
    Likewise, if $p'$ is feasible but dominated by $\overline{p}$, then by Property~\ref{property.domination_ours}\ref{property.domination_pps} and the feasibility of $p' \oplus s$, $\overline{p} \oplus s$ is also feasible and dominates $p' \oplus s$. 
    Therefore both $p'$ and $s$ must be non-dominated. 
    \hfill\Halmos

\proof{{Proof of Theorem~\ref{thm.exact_PP}.}}
Finite termination of Algorithms~\ref{alg.FindNonDominatedSubpaths} and~\ref{alg.FindNonDominatedPaths} is guaranteed per Assumption~\ref{ass.technical}. At the subpath level, the resource $\TendMin{s}$ is nonnegative, strictly increasing along extensions (since $t_{i,j} > 0$), and bounded above by $T$. At the path level, the resource $\TendMin{p}$ is nonnegative, strictly increasing along extensions (since $\Ttravel{s} > 0$ for all $s$ consisting of more than 1 node), and bounded above by $T$.

We show that $\Sgen = \calS^\text{nd}$, the set of non-dominated partial subpaths, up to dominance equivalence (i.e. when two labels dominate each other, one is retained). Since $\Sresult$ is the subset of partial subpaths in $\Sgen$ that are also full subpaths, it automatically implies that Algorithm~\ref{alg.FindNonDominatedSubpaths} returns $\Sresult = \calS^*$, the set of non-dominated subpaths, up to dominance equivalence.
\begin{itemize}[]
    \item \underline{$\Sgen \supseteq \calS^\text{nd}$.} 
    Suppose for contradiction that there exists $s \in \calS^\text{nd} \setminus \Sgen$. Choose $s$ with the smallest $\TendMin{s}$.
    Define $s'$ such that $s = s' \oplus a$; by Lemma~\ref{lemma.subpath_truncate_nondom}, we have that $s' \in \calS^\text{nd}$. 
    By Assumption~\ref{ass.technical}, $\TendMin{s'} < \TendMin{s}$, so our construction implies that $s' \in \Sgen$. Consider $\Sgen$ and $\Squeue$ at the point in the algorithm where $s'$ is moved from $\Squeue$ to $\Sgen$. Then, $s = s' \oplus a$ is explored in Step 2 of the algorithm, and added to $\Squeue$, a contradiction.
    
    \item \underline{$\Sgen \subseteq \calS^\text{nd}$.}
    Suppose for contradiction that there exists $s \in \Sgen \setminus \calS^\text{nd}$, i.e. there exists a non-dominated partial subpath $s' \in \calS^\text{nd}$ such that $s' \succeq s$.
    As seen earlier, $s'$ is added to $\Squeue$ at some point of the algorithm and remains in it until it gets added to $\Sgen$. 
    Express $s' := s'' \oplus a$; then we have $\TendMin{s} \geq \TendMin{s'}$ (by domination) and $\TendMin{s'} > \TendMin{s''}$ (by Assumption~\ref{ass.technical}). 
    At the iteration where $s'$ is added to $\Squeue$, $\argmin \Set{ \TendMin{s} | s \in \Squeue } = \TendMin{s''}$; 
    at the iteration where $s$ is added to $\Sgen$, $\argmin \Set{ \TendMin{s} | s \in \Squeue } = \TendMin{s}$. 
    By Assumption~\ref{ass.technical}, $\argmin \Set{ \TendMin{s} | s \in \Squeue}$ is nondecreasing over the algorithm, so $s'$ is added to $\Squeue$ before $s$ is added to $\Sgen$. This contradicts that $s \in \Sgen$, since $s' \succeq s$ would either remove $s$ from $\Squeue$ or prevent $s$ from being added to $\Squeue$.
\end{itemize}

Next, Algorithm~\ref{alg.FindNonDominatedPaths}, when initialized with $\calS^*$ (the set of non-dominated subpaths from Algorithm~\ref{alg.FindNonDominatedSubpaths}), returns $\Presult = \calP^*$, the set of non-dominated paths, up to dominance equivalence. Again, we show the stronger claim: $\Pgen = \calP^\text{nd}$, the set of non-dominated partial paths ending at non-customer nodes, up to dominance equivalence. The proof proceeds exactly as before, using $\TendMin{p}$ as the path-level resource satisfying Assumption~\ref{ass.technical}, and replacing Lemma~\ref{lemma.subpath_truncate_nondom} with Lemma~\ref{lemma.path_truncate_nondom}. 

Finally, we show that $\Presult = \calP^*$ contains paths of negative reduced cost if any exist. Suppose $p$ is a path with negative reduced cost. If it is dominated, there exists a non-dominated path $p'\in\calP^*$ such that $p' \succeq p$. This implies that $\cbar{p'} \leq \cbar{p} < 0$. Therefore $\calP^*$ finds paths of negative reduced cost, if any exist.
\hfill\Halmos

\subsection*{{Proof of Proposition~\ref{prop.arbitrarily_large}.}}

Consider our two-phase label-setting approach for the pricing problem, and assume that we define subpath dominance based on Definitions~\ref{def.subpath_domination}\ref{def.subpath_domination_rc}--\ref{def.subpath_domination_TendMin}, thus excluding \ref{def.subpath_domination_TstartMax} and \ref{def.subpath_domination_Q}. From the proof of Theorem~\ref{thm.EVRPTWNL_correctness}, we can prove that Properties~\ref{property.domination_ours}\ref{property.domination_ssa}--\ref{property.domination_pps} are satisfied but Property~\ref{property.domination_ours}\ref{property.domination_pss} is not.

Consider now a problem with linear charging, $T = 10$, $B = 3$ and $D = 5$ where there is a single non-dominated partial path $p$ from depot $0$ to charging station $1$, and two subpaths $s_1, s_2$ starting at $1$ and ending at another depot $4$. The node sequence of $s_1$ is $(1, 2, 3, 4)$, while the node sequence of $s_2$ is $(1, 3, 2, 4)$. Customers 2 and 3 have time windows $[6, 9]$ and $[7, 8]$ respectively, and a load of 1 each. Travel times and charge requirements along arcs are given in Figure~\ref{fig.counterexample}. Reduced costs are given as $\overline{c}_{i,j} = 0$ for $(i,j) \in \{(1,2),(1,3),(2,3),(3,2),(2,4),(3,4)\}$. We can compute the resources for $s_1$ and $s_2$: $(\cbar{s_1}, D(s_1), B(s_1), \TendMin{s_1}, \TstartMax{s_1}, \Ttravel{s_1}) = (0, 2, 3, 8, 6, 3)$ and $(\cbar{s_2}, D(s_2), B(s_2), \TendMin{s_2}, \TstartMax{s_2}, \Ttravel{s_2}) = (0, 2, 3, 9, 7, 3)$. 
Since $s_1$ dominates $s_2$ under Definitions~\ref{def.subpath_domination}\ref{def.subpath_domination_rc}--\ref{def.subpath_domination_TendMin}, $s_2$ would be pruned in the first phase under Properties~\ref{property.domination_ours}\ref{property.domination_ssa}--\ref{property.domination_pps}.

\begin{figure}[h!]
    \centering
    \begin{tikzpicture}
        [
            node distance=1.5mm,
            inner sep=1mm,
            sourcesink/.style={rectangle,draw=black,fill=purple!50!white,thick,minimum size=4mm},
            charger/.style={rectangle,draw=black,fill=yellow,thick,minimum size=4mm},
            customer/.style={circle,draw=black,thick,minimum size=4mm},
            >={Stealth},
            legend/.style={right,inner sep=2mm,text height=1.5ex,text depth=0.25ex},
        ]
        \node[sourcesink]   (d1)    at ( 0.0, 0.0)  {$0$};
        \node[charger]      (r1)    at ( 8.0, 0.0)  {$1$};
        \node[below=of r1]  (r1_l)                  {\small $[0, 10]$};
        \node[customer]     (c1)    at (10.0, 1.0)  {$2$};
        \node[above=of c1]  (c1_l)                  {\small $[6, 9]$};
        \node[customer]     (c2)    at (10.0,-1.0)  {$3$};
        \node[below=of c2]  (c1_l)                  {\small $[7, 8]$};
        \node[sourcesink]   (d2)    at (12.0, 0.0)  {$4$};
        \node[below=of d2]  (d2_l)                  {\small $[0, 10]$};
        
        \node[sourcesink]   (d0)    at (13.5, 1.4)  {};
        \node[charger]      (r0)    at (13.5, 0.7)  {};
        \node[customer]     (c0)    at (13.5, 0.0)  {};
        \node[legend,xshift=3pt]   (d0)    at (d0)  {\small Depot};
        \node[legend,xshift=3pt]   (r0)    at (r0)  {\small Charger};
        \node[legend,xshift=3pt]   (c0)    at (c0)  {\small Customer};

        \draw (d1) [bend right=15] to node[auto, inner sep=1mm] {$p$} (r1);
        \draw (r1) to node[auto, inner sep=1mm] {$1$} (c1);
        \draw (r1) to node[auto, inner sep=1mm] {$1$} (c2);
        \draw (c1) to node[auto, inner sep=1mm] {$1$} (c2);
        \draw (c1) to node[auto, inner sep=1mm] {$1$} (d2);
        \draw (c2) to node[auto, inner sep=1mm] {$1$} (d2);
        
        \end{tikzpicture}%
    \caption{An instance where the absence of Property~\ref{property.domination_ours}\ref{property.domination_pss} results in an incorrect algorithm.}
    \label{fig.counterexample}
\end{figure}

Suppose $(\cbar{p}, D(p), \TendMin{p}, \BendMax{p}) = (-M, 1, 5, 1)$. We remark that $p \oplus s_1$ is not feasible, since $\TendMin{p} + \gT{\BendMax{p}}{\max \{ \BendMax{p}, B(s_1) \}}
= 5 + (3 - 1) = 7 > 6 = \TstartMax{s_1}$, thus violating condition~\ref{prop.path_REFs_feasible_TstartMax} in Proposition~\ref{prop.path_REFs}. Therefore, under Properties~\ref{property.domination_ours}\ref{property.domination_ssa}--\ref{property.domination_pps} alone, there would be no feasible paths from depot to depot, and the column generation algorithm would then terminate at that iteration.

Let us now assume that we apply the two-phase label-setting algorithm with the full Property~\ref{property.domination_ours}, including Property~\ref{property.domination_ours}\ref{property.domination_pss}. Then, none of the subpaths $s_1$ and $s_2$ is dominated in the first phase, and $p \oplus s_2$ is feasible with reduced cost $-M$ in the second phase. Note, in particular, that in this example $s_1 \succeq s_2$ no longer implies $p \oplus s_1 \succeq p \oplus s_2$. Therefore, the two-phase label-setting algorithm returns a path of reduced cost $-M$, and the column generation algorithm continues. As $M$ becomes arbitrarily large, this means that a pricing algorithm that uses only Properties~\ref{property.domination_ours}\ref{property.domination_ssa}--\ref{property.domination_pps} can fail to find paths of arbitrarily negative reduced cost. 
\hfill\Halmos

\subsection*{{Proof of Remark~\ref{remark.property_i_ii_sufficient}.}}

We prove that Property~\ref{property.domination_ours}\ref{property.domination_pss} is automatically satisfied if:
\begin{enumerate}[(i),noitemsep]
    \item 
    \label{remark.sufficient_feasibility}
    (Domination implies feasibility) Given feasible subpaths $s_1 \succeq s_2$ both extending feasible partial path $p$, the feasibility of $p \oplus s_2$
    implies the feasibility of $p \oplus s_1$;
    \item
    \label{remark.sufficient_subpathdom}
    (Subpath domination) Domination at the subpath level is given by $s_1 \succeq s_2 \iff r_\text{S}(s_1) \leq r_\text{S}(s_2)$ for a set of subpath-level resources $r_\text{S} \in \calR_\text{S}$;
    \item 
    \label{remark.sufficient_pathdom}
    (Path domination) domination at the path level is given by $p_1 \succeq p_2 \iff r_\text{P}(p_1) \leq r_\text{P}(p_2)$ for a set of path-level resources $r_\text{P} \in \calR_\text{P}$;
    \item 
    \label{remark.sufficient_monotonicity}
    (Monotonicity of path REFs along subpath extensions) For each path-level resource $r_\text{P}$, its REF along subpath extensions is a monotonic function of the set of subpath-level resources $r_\text{S} \in \calR_\text{S}(r_\text{P})$ involved. Formally, $r_\text{P}(p \oplus s)$ is a function $f$ of $r_\text{P}(p)$, possible other path-level resources $r'_\text{P}(p)$ for $r'_\text{P} \in \calR_\text{P}$, and some subpath-level resources $r_\text{S}(s)$ for $r_\text{S} \in \calR_\text{S}(r_\text{P})$, and:
    \begin{equation}
        r_\text{P}(p \oplus s) 
        = f \Big( 
            r_\text{P}(p), 
            \{r'_\text{P}(p)\}_{r'_\text{P} \in \calR_\text{P}}, 
            \{ r_\text{S}(s) \}_{r_\text{S} \in \calR_\text{S}(r_\text{P})}
        \Big)
        \text{ is increasing in each }
        r_\text{S}(s), \text{ for } r_\text{S} \in \calR_\text{S}(r_\text{P}).
    \end{equation}
\end{enumerate}

Suppose that subpaths $s_1 \succeq s_2$ both extend feasible partial path $p$ and that $p \oplus s_2$ is feasible. This gives feasibility of $p \oplus s_1$ by \ref{remark.sufficient_feasibility}. For each path-level resource $r_\text{P}$, by monotonicity \ref{remark.sufficient_monotonicity} and $s_1 \succeq s_2$ \ref{remark.sufficient_subpathdom}:
\begin{alignat}{3}
    r_\text{P}(p \oplus s_1) 
    &= f \Big( 
        r_\text{P}(p), 
        \{r'_\text{P}(p)\}_{r'_\text{P} \in \calR_\text{P}}, 
        \{ r_\text{S}(s_1) \}_{r_\text{S} \in \calR_\text{S}(r_\text{P})}
    \Big)\\
    &\leq f \Big( 
        r_\text{P}(p), 
        \{r'_\text{P}(p)\}_{r'_\text{P} \in \calR_\text{P}}, 
        \{ r_\text{S}(s_2) \}_{r_\text{S} \in \calR_\text{S}(r_\text{P})}
    \Big)\\
    &= r_\text{P}(p \oplus s_2) 
\end{alignat}
Therefore by \ref{remark.sufficient_pathdom} $p \oplus s_1 \succeq p \oplus s_2$.
\hfill\Halmos

We next provide examples where the conditions are satisfied. Consider the pricing problem for the EVRP with nonlinear charging, no customer time windows, and no elementarity restrictions. Let $\calR_\text{S}$ be the resources $\cbar{\cdot}$, $D(\cdot)$, $B(\cdot)$, and $\Ttravel{\cdot}$. Let $\calR_\text{P}$ be the resources $\cbar{\cdot}$, $D(\cdot)$, $\TendMin{\cdot}$, and $\Teff{\cdot}$. 
Recall that $\TendMin{\cdot}$ and $\Teff{\cdot}$ are so defined, so that $\TendMin{p_1} \leq \TendMin{p_2}$ and $\Teff{p_1} \leq \Teff{p_2}$ together are equivalent to 
Definition~\ref{def.path_domination}\ref{def.path_domination_BendMax} (Lemma~\ref{lemma.path_domination_componentwise}).

Condition~\ref{remark.sufficient_subpathdom} and~\ref{remark.sufficient_pathdom} are satisfied by construction. The REFs for $p \oplus s$ in terms of $p$ are $\TendMin{p \oplus s} := \TendMin{p} + \tau + \Ttravel{s}$ and $\BendMax{p \oplus s} := \Pos{\BendMax{p} - B(s)}$, with $\tau := \max \left\{0, \gT{\BendMax{p}}{B(s)} \right\}$. After some algebra we also have:
\begin{equation}
    \Teff{p \oplus s} = \Teff{p} + \Ttravel{s} + \gT{\Pos{\BendMax{p} - B(s)}}{\Pos{\BendMax{p} - B(s)} + B(s)}
\end{equation}
The path-level REFs $\cbar{p \oplus s}$, $D(p \oplus s)$ are monotonic in terms of subpath-level resources $\cbar{s}$ and $D(s)$ respectively. The path-level REFs $\TendMin{p \oplus s}$ and $\Teff{p \oplus s}$ are also both monotonic in terms of the subpath-level resources $\{ B(s), \Ttravel{s}\}$. This proves Condition~\ref{remark.sufficient_monotonicity}. Finally, feasibility of $p \oplus s_1$ is satisfied automatically for conditions~\ref{prop.path_REFs_feasible_load}, \ref{prop.path_REFs_feasible_T} (by monotonicity of $\TendMin{p \oplus s}$), and \ref{prop.path_REFs_feasible_TstartMax} (which now reads $\TendMin{p} + \tau \leq T - \Ttravel{s}$, by monotonicity of $\TendMin{p \oplus s}$). This proves~\ref{remark.sufficient_feasibility}. 

\section{Proofs in Section 5}
\label{app.proofs_5}

\subsection{Proofs in Section 5.1}
\label{app.proofs_51}

We first state a preparatory lemma:

\begin{lemma}
    \label{lemma.ngroute_ngset}
    Consider an \textit{ng}-neighborhood $\calN$, a node sequence $N = (n_0, \dots, n_m)$, and an arc $(n_m, n_{m+1}) \in \calA$. Then:
    $$
    (n_0, \dots, n_{m+1}) \text{ is \textit{ng}-feasible w.r.t. $\calN$}
    \iff N \text{ is \textit{ng}-feasible w.r.t. $\calN$, and } n_{m+1} \notin \Pi(N).
    $$
\end{lemma}

\begin{proof}{Proof of Lemma~\ref{lemma.ngroute_ngset}.}
\begin{itemize}
    \item[$(\Leftarrow)$] Let $0 \leq j < k \leq m+1$ be such that $n_j = n_k$. If $k \neq m+1$, then there exists $\ell$ with $j < \ell < k$ with $n_j \notin N_{n_\ell}$, because $N$ is \textit{ng}-feasible with respect to $\calN$. If $k = m+1$, then $n_j = n_{m+1}$ and $n_j \notin \Pi(N)$, so $n_j \notin \bigcap_{\rho = j+1}^{m} N_{n_{\rho}}$. Therefore, there exists $\ell$ such that $j+1 \leq l \leq m$ (i.e., $j < \ell < m+1$) such that $n_j \notin N_{n_\ell}$. Thus, $(n_0, \dots, n_{m+1})$ is \textit{ng}-feasible with respect to $\calN$.
    \item[$(\Rightarrow)$] $N$ is clearly \textit{ng}-feasible with respect to $\calN$. Assume by contradiction that $n_{m+1} \in \Pi(N)$. There exists $r \leq m-1$ such that $n_r = n_{m+1}$ and $n_r \in \bigcap_{\rho = r+1}^{m} N_{n_\rho}$. Hence for $j = r$ and $k = m+1$, there does not exist any $j < \ell < k$ such that $n_j \notin N_{n_\ell}$. This implies that $(n_0, \dots, n_{m+1})$ is not \textit{ng}-feasible with respect to $\calN$, leading to a contradiction.
    \hfill\Halmos
\end{itemize}
\end{proof}

\subsection*{{Proof of Proposition~\ref{prop.subpath_REFs_ngroute}.}}

Let $s$ be a \textit{ng}-feasible partial subpath with node sequence $(n_0, \dots, n_m)$ and $a = (n_m, n_{m+1})$ be an arc extension of $s$. 
Since $n_{m+1} \notin \Pi(s)$, 
Lemma~\ref{lemma.ngroute_ngset} implies that $s \oplus a$ is \textit{ng}-feasible. We extend the forward \textit{ng}-set, backward \textit{ng}-set, and \textit{ng}-residue as follows:
\begin{alignat}{3}
    \Pi(s \oplus a) 
    & = 
    \textstyle
    \Set{
        n_r
        |
        n_r \in \bigcap_{\rho = r + 1}^{m+1} N_{n_{\rho}},
        \ r \in \{0, \dots, m\}
    }
    \cup \{ n_{m+1} \}
    \\
    & = \textstyle
    \left( \Set{
        n_r
        |
        n_r \in \bigcap_{\rho = r + 1}^{m} N_{n_{\rho}},
        \ r \in \{0, \dots, m-1\}
    } \cap N_{n_{m+1}} \right)
    \cup \{ n_{m+1} \}
    = (\Pi(s) \cap N_{n_{m+1}}) \cup \{ n_{m+1} \}
    \\
    \Omega(s \oplus a)
    & = 
    \textstyle
    \bigcap_{\rho = 0}^{m+1} N_{n_{\rho}}
    = 
    \bigcap_{\rho = 0}^{m} N_{n_{\rho}}\cap N_{n_{m+1}}
    = \Omega(s)\cap N_{n_{m+1}}
    \\
    \Pi^{-1}(s \oplus a) 
    & = 
    \textstyle
    \{ n_0 \} \cup \Set{
        n_r
        |
        n_r \in \bigcap_{\rho = 0}^{r-1} N_{n_{\rho}},
        \ r \in \{1, \dots, m+1\}
    }
    \\
    & = \begin{cases}
        \{ n_0 \} \cup \Set{
            n_r
            |
            n_r \in \bigcap_{\rho = 0}^{r-1} N_{n_{\rho}},
            \ r \in \{1, \dots, m\}
        } \cup \{n_{m+1}\}
        & \text{ if } n_{m+1} \in \bigcap_{\rho = 0}^{m} N_{n_{\rho}},
        \\
        \{ n_0 \} \cup \Set{
            n_r
            |
            n_r \in \bigcap_{\rho = 0}^{r-1} N_{n_{\rho}},
            \ r \in \{1, \dots, m\}
        } & \text{ otherwise.}
    \end{cases}
    \\
    & = \Pi^{-1}(s) \cup \left( \{n_{m+1}\} \cap \Omega(s) \right)
    \hfill\Halmos
\end{alignat}

\subsection*{{Proof of Proposition~\ref{prop.path_REFs_ngroute}.}}

Let $p$ have node sequence $(n_0, \dots, n_m)$ and $s$ have node sequence $(n_m, \dots, n_M)$. First suppose $\Pi(p) \cap \Pi^{-1}(s) \not\subseteq \{ n_m\}$, i.e. there exists $\overline{n} \in \Pi(p) \cap \Pi^{-1}(s)$ with $\overline{n} \neq n_m$. 
Let $j<m$ and $k>m$ be such that $n_j = n_k = \overline{n}$. Since $n_j \in \Pi(p)$, $n_j\in N_{n_\ell}$ for all $l \in \{j+1, \dots, m\}$. Similarly, since $n_k \in \Pi^{-1}(s)$, $n_k\in N_{n_\ell}$ for all $l \in \{m, \dots, k-1\}$. Therefore, $n_j=n_k$ and $n_j\in N_{n_\ell}$ for all $l \in \{j+1, \dots, k-1\}$, which proves that $p \oplus s$ is not \textit{ng}-feasible. Therefore, $p \oplus s$ is \textit{ng}-feasible iff $\Pi(p) \cap \Pi^{-1}(s) \subseteq \{ n_m \}$.

We next characterize $\Pi(p \oplus s)$. By the definition of $\Pi(\cdot)$:
\begin{alignat}{3}
    \Pi(p \oplus s) 
    & = 
    \textstyle
    \Set{
        n_r 
        |
        n_r \in \bigcap_{\rho = r + 1}^{M} N_{n_{\rho}},
        \ r \in \{0, \dots, m\}
    } \cup \Set{
        n_r 
        |
        n_r \in \bigcap_{\rho = r + 1}^{M} N_{n_{\rho}},
        \ r \in \{m, \dots, M-1\}
    } \cup \{n_M \}
    \\
    & = 
    \textstyle
    \Set{
        n_r 
        |
        n_r \in \bigcap_{\rho = r + 1}^{M} N_{n_{\rho}},
        \ r \in \{0, \dots, m\}
    } \cup \Pi(s)
    \qquad
    \text{(by Proposition~\ref{prop.subpath_REFs_ngroute})}
    \\
    & = 
    \textstyle
    \left( \left( \Set{
        n_r 
        |
        n_r \in \bigcap_{\rho = r + 1}^{m} N_{n_{\rho}},
        \ r \in \{0, \dots, m-1\}
    } \cup \{ n_m \} \right) \cap \bigcap_{\rho = m}^{M} N_{n_{\rho}} \right) \cup \Pi(s)
    \\
    & = \left( \Pi(p) \cap \Omega(s) \right) \cup \Pi(s)
    \hfill\Halmos
\end{alignat}

\subsection*{{Proof of Proposition~\ref{prop.ngroute_correctness}.}}

We focus on the new elements of the proof. 
All other parts follow Theorem~\ref{thm.EVRPTWNL_correctness}.

\subsubsection*{Property~\ref{property.domination_ours}\ref{property.domination_ssa}.}

Let $s_1 \succeq s_2$ be \textit{ng}-feasible partial subpaths.
Therefore, from Definition~\ref{def.subpath_domination_ngroute} we have: 
$\Pi(s_1) \subseteq \Pi(s_2)$, 
$\Omega(s_1) \subseteq \Omega(s_2)$, 
and $\Pi^{-1}(s_1) \subseteq \Pi^{-1}(s_2)$. 
Let $a = (n, n')$ be a common arc extension and suppose that $s_2 \oplus a$ is \textit{ng}-feasible. 
Then $n' \notin \Pi(s_2)$ by Lemma~\ref{lemma.ngroute_ngset}. 
This implies $n' \notin \Pi(s_1)$ too and $s_1 \oplus a$ is also \textit{ng}-feasible.
Moreover, the three updates in Proposition~\ref{prop.subpath_REFs_ngroute} preserve set inclusion:
\begin{alignat}{3}
    \Pi(s_1 \oplus a) 
    = (\Pi(s_1) \cap N_{n'}) \cup \{n'\} 
    & \subseteq (\Pi(s_2) \cap N_{n'}) \cup \{n'\} 
    = \Pi(s_2 \oplus a)
    \\
    \Omega(s_1 \oplus a) 
    = \Omega(s_1) \cap N_{n'} 
    & \subseteq \Omega(s_2) \cap N_{n'} 
    = \Omega(s_2 \oplus a)
    \\
    \Pi^{-1}(s_1 \oplus a) 
    = \Pi^{-1}(s_1) \cup ( \{ n' \} \cap \Omega(s_1) ) 
    & \subseteq \Pi^{-1}(s_2) \cup ( \{ n' \} \cap \Omega(s_2) ) 
    = \Pi^{-1}(s_2 \oplus a)
\end{alignat}

\subsubsection*{Property~\ref{property.domination_ours}\ref{property.domination_pps}.}

Let $p_1 \succeq p_2$ be \textit{ng}-feasible partial paths ending at node $n$.
Therefore, from Definition~\ref{def.path_domination_ngroute}, we have: $\Pi(p_1) \subseteq \Pi(p_2)$. 
Let $s$ be a common subpath extension and suppose that $p_2 \oplus s$ is \textit{ng}-feasible. 
Then $\Pi(p_2) \cap \Pi^{-1}(s) \subseteq \{ n \}$ by Proposition~\ref{prop.path_REFs_ngroute}. 
Therefore $\Pi(p_1) \cap \Pi^{-1}(s) \subseteq \Pi(p_2) \cap \Pi^{-1}(s) \subseteq \{ n \}$, 
and $p_1 \oplus s$ is also \textit{ng}-feasible.
Moreover, per Proposition~\ref{prop.path_REFs_ngroute}:
\begin{alignat}{3}
    \Pi(p_1 \oplus s)
    & = \Pi(s) \cup (\Pi(p_1) \cap \Omega(s)) 
    \subseteq \Pi(s) \cup (\Pi(p_2) \cap \Omega(s))
    = \Pi(p_2 \oplus s)
\end{alignat}

\subsubsection*{Property~\ref{property.domination_ours}\ref{property.domination_pss}.}

Let $s_1 \succeq s_2$ be \textit{ng}-feasible subpaths.
Therefore, from Definition~\ref{def.subpath_domination_ngroute} we have: 
$\Pi(s_1) \subseteq \Pi(s_2)$, 
$\Omega(s_1) \subseteq \Omega(s_2)$, 
and $\Pi^{-1}(s_1) \subseteq \Pi^{-1}(s_2)$. 
Let $p$ be a \textit{ng}-feasible partial path ending at node $n$, 
and $s_1, s_2$ both extend $p$, and suppose $p \oplus s_2$ is \textit{ng}-feasible.
Then, $\Pi(p) \cap \Pi^{-1}(s_2) \subseteq \{ n \}$, by Proposition~\ref{prop.path_REFs_ngroute}.
Therefore $\Pi(p) \cap \Pi^{-1}(s_1) \subseteq \Pi(p) \cap \Pi^{-1}(s_2) \subseteq \{ n \}$, 
and $p \oplus s_1$ is also \textit{ng}-feasible.
Moreover, per Proposition~\ref{prop.path_REFs_ngroute}:
\begin{alignat}{3}
    \Pi(p \oplus s_1)
    & = \Pi(s_1) \cup (\Pi(p) \cap \Omega(s_1)) 
    \subseteq \Pi(s_2) \cup (\Pi(p) \cap \Omega(s_2))
    = \Pi(p \oplus s_2)
    \hfill\Halmos
\end{alignat}
\subsection{Proofs in Section 5.2}
\label{app.proofs_52}

Before we prove the results in this section, we first note that the reduced cost resource for a (partial) subpath or path with node sequence $N = (n_0, \dots, n_m)$ is given by:
\begin{alignat}{3}
    \label{eq.SRC_resources_rc}
    && 
    \cbar{N} 
    & = \sum_{\ell=0}^{m-1} \overline{c}_{n_\ell,n_{\ell+1}}
    - \sum_{q \in \calQ} \lambda_q \cdot \alpha_{S_q,\bw^q}(N)
    \\
    \label{eq.SRC_resources_alpha}
    \text{ where } \quad 
    &&
    \forall \ q \in \calQ, 
    \ 
    \alpha_{S_q,\bw^q}(N) 
    & = 
    \textstyle
    \floor*{
        \sum_{\ell=0}^{m} w^q_{n_\ell} \ind{n_\ell \in S_q}
    }
    \\
    \label{eq.SRC_resources_alpha_frac}
    &&
    \alpha_{q}(N) 
    & = 
    \textstyle
    \fracpart{
        \sum_{\ell=0}^{m} w^q_{n_\ell} \ind{n_\ell \in S_q}
    }
\end{alignat}

\begin{lemma}
    \label{lemma.SRC_resources_concat}
    For two node sequences $N = (n_0, \dots, n_m)$ and $N' = (n_{m+1}, \dots, n_M)$, we have:
    \begin{alignat}{3}
        \label{eq.SRC_resources_concat_rc}
        \cbar{n_0, \dots, n_M}
        & = \cbar{N} 
        + \overline{c}_{n_m, n_{m+1}} 
        + \cbar{N'} 
        - \sum_{q \in \calQ} 
        \lambda_q \cdot \ind{\alpha_q(N) + \alpha_q(N') \geq 1}
        \\
        \label{eq.SRC_resources_concat_alpha}
        \alpha_{S_q, \bw^q}(n_0, \dots, n_M)
        & = \alpha_{S_q, \bw^q}(N)
        + \alpha_{S_q, \bw^q}(N')
        + \ind{ \alpha_q(N) + \alpha_q(N') \geq 1 }
        \\
        \label{eq.SRC_resources_concat_alpha_frac}
        \alpha_{q}(n_0, \dots, n_M)
        & = \fracpart{
            \alpha_q(N) 
            + \alpha_q(N')
        }
    \end{alignat}
\end{lemma}
\begin{proof}{Proof of Lemma~\ref{lemma.SRC_resources_concat}.}
    This is true by defining $x = \sum_{\ell=0}^{m} w^q_{n_\ell} \ind{n_\ell \in S_q}$, $y = \sum_{\ell=m+1}^{M} w^q_{n_\ell} \ind{n_\ell \in S_q}$, and the identities:
    \begin{equation*}
        \floor{x + y} = \floor{x} + \floor{y} + \ind{\fracpart{x} + \fracpart{y} \geq 1};
        \quad 
        \fracpart{x + y} = \fracpart{ \fracpart{x} + \fracpart{y} }
        \hfill\Halmos
    \end{equation*}
\end{proof}

\subsubsection*{{Proof of Propositions~\ref{prop.subpath_REFs_SRC} and~\ref{prop.path_REFs_SRC}.}}
Apply Lemma~\ref{lemma.SRC_resources_concat}:
\begin{itemize}[noitemsep]
    \item Proposition~\ref{prop.subpath_REFs_SRC}: If $s$ has node sequence $(n_0, \dots, n_m)$ and $a = (n_m, n_{m+1})$, set $N$ as the node sequence of $s$ and $N' = (n_{m+1},)$. Note that $\alpha_{S_q, \bw^q}(N') = 0$ since it is the floor of a single summand less than 1.
    \item Proposition~\ref{prop.path_REFs_SRC}: If $p$ has node sequence $(n_0, \dots, n_m)$ and $s$ has node sequence $(n_m, \dots, n_M)$, set $N$ as the node sequence of $p$ and $N' = (n_{m+1}, \dots, n_M)$. 
    Note that $\cbar{s} = \overline{c}_{n_m,n_{m+1}} + \cbar{N'}$.
    Also, $\alpha_q(s) = \alpha_q(N')$ and $\alpha_{S_q, \bw^q}(s) = \alpha_{S_q, \bw^q}(N')$ since $s$ starts at charging station $n_m$ which must not belong in $S_q$.
    \hfill\Halmos
\end{itemize}

\subsubsection*{{Proof of Proposition~\ref{prop.SRCs_correctness}.}}

\begin{lemma}
    \label{lemma.SRCs_fracpart_cyclic}
    Let $a, b, c \in [0, 1)$. 
    Then:
    $
        \ind{a > b} 
        - \ind{a + c \geq 1}
        + \ind{b + c \geq 1}
        = \ind{\fracpart{a + c} > \fracpart{b + c}}
    $.
\end{lemma}
\begin{proof}{Proof of Lemma~\ref{lemma.SRCs_fracpart_cyclic}}
    We consider cases on $a, b, c$, summarized in the table below:
    \hfill\Halmos
    \begin{table}[H]
        \vspace{-12pt}
        \small
        \centering
        \begin{tabular}{cccc}
            \toprule
            $a > b$ 
            & $c \in [0, 1-a)$ & $c \in [1-a, 1-b)$ & $c \in [1-b, 1)$ \\
            & $1 - 0 + 0 = 1$ & $1 - 1 + 0 = 0$ & $1 - 1 + 1 = 1$ \\
            \midrule
            $a \leq b$ & $c \in [0, 1-b)$ & $c \in [1-b, 1-a)$ & $c \in [1-a, 1)$ \\
            & $0 - 0 + 0 = 0$ & $0 - 0 + 1 = 1$ & $0 - 1 + 1 = 0$ \\
            \bottomrule
        \end{tabular}
        \vspace{-12pt}
    \end{table}
\end{proof}

\proof{Proof of Proposition~\ref{prop.SRCs_correctness}.}

We focus on the new elements of the proof. 
All other parts follow Theorem~\ref{thm.EVRPTWNL_correctness}.

\subsubsection*{Property~\ref{property.domination_ours}\ref{property.domination_ssa}.}

Let $s_1 \succeq s_2$ be feasible partial subpaths. Let $a$ be a common arc extension, and suppose that $s_2 \oplus a$ is feasible. $s_1 \oplus a$ is automatically feasible by Theorem~\ref{thm.EVRPTWNL_correctness}. We wish to show:
\begin{equation*}
    \cbar{s_1} - \cbar{s_2} \leq \sum_{q \in \calQ} \lambda_q \ind{\alpha_q(s_1) > \alpha_q(s_2)}
    \implies 
    \cbar{s_1 \oplus a} - \cbar{s_2 \oplus a} \leq \sum_{q \in \calQ} \lambda_q \ind{\alpha_q(s_1 \oplus a) > \alpha_q(s_2 \oplus a)}
\end{equation*}
We have:
\begin{alignat*}{3}
    & \ \cbar{s_1 \oplus a} - \cbar{s_2 \oplus a}
    \\
    = & \ 
    \cbar{s_1} - \cbar{s_2} 
    - \sum_{q \in \calQ} \lambda_q \ind{\alpha_q(s_1) + w^q_{n'} \ind{n' \in S_q} \geq 1}
    + \sum_{q \in \calQ} \lambda_q \ind{\alpha_q(s_2) + w^q_{n'} \ind{n' \in S_q} \geq 1}
    \\
    \leq & \ 
    \sum_{q \in \calQ} \lambda_q \Big[ 
        \ind{\alpha_q(s_1) > \alpha_q(s_2)}
        - \ind{\alpha_q(s_1) + w^q_{n'} \ind{n' \in S_q} \geq 1}
        + \ind{\alpha_q(s_2) + w^q_{n'} \ind{n' \in S_q} \geq 1}
    \Big]
    \\
    = & \ 
    \sum_{q \in \calQ} \lambda_q \ind{\alpha_q(s_1 \oplus a) > \alpha_q(s_2 \oplus a)}
\end{alignat*}
as desired.
The first equality is by Proposition~\ref{prop.subpath_REFs_SRC}, and the second inequality by Definition~\ref{def.subpath_path_domination_SRC}.
The third equality is by Lemma~\ref{lemma.SRCs_fracpart_cyclic} (letting $a = \alpha_q(s_1), b = \alpha_q(s_2), c = w^q_{n'} \ind{n' \in S_q}$) and Equation~\eqref{eq.SRC_resources_concat_alpha_frac}.

\subsubsection*{Property~\ref{property.domination_ours}\ref{property.domination_pps}.}

Let $p_1 \succeq p_2$ be feasible partial paths. Let $s$ be a common subpath extension, and suppose that $p_2 \oplus s$ is feasible. $p_1 \oplus s$ is automatically feasible by Theorem~\ref{thm.EVRPTWNL_correctness}. We wish to show:
\begin{equation*}
    \cbar{p_1} - \cbar{p_2} \leq \sum_{q \in \calQ} \lambda_q \ind{\alpha_q(p_1) > \alpha_q(p_2)}
    \implies 
    \cbar{p_1 \oplus s} - \cbar{p_2 \oplus s} \leq \sum_{q \in \calQ} \lambda_q \ind{\alpha_q(p_1 \oplus s) > \alpha_q(p_2 \oplus s)}
\end{equation*}
The proof is identical to that for Property~\ref{property.domination_ours}\ref{property.domination_ssa}, except that 
Proposition~\ref{prop.path_REFs_SRC} is applied instead of Proposition~\ref{prop.subpath_REFs_SRC}, 
and Lemma~\ref{lemma.SRCs_fracpart_cyclic} is applied by letting $a = \alpha_q(p_1), b = \alpha_q(p_2), c = \alpha_q(s)$ for each $q$.

\subsubsection*{Property~\ref{property.domination_ours}\ref{property.domination_pss}.}

Let $p$ be a feasible partial paths, and let $s_1 \succeq s_2$ be subpaths both extending $p$. Suppose that $p \oplus s_2$ is feasible. $p \oplus s_1$ is automatically feasible by Theorem~\ref{thm.EVRPTWNL_correctness}. We wish to show:
\begin{equation*}
    \cbar{s_1} - \cbar{s_2} 
    \leq \sum_{q \in \calQ} \lambda_q \ind{\alpha_q(s_1) > \alpha_q(s_2)}
    \implies 
    \cbar{p \oplus s_1} - \cbar{p \oplus s_2} 
    \leq \sum_{q \in \calQ} \lambda_q \ind{\alpha_q(p \oplus s_1) > \alpha_q(p \oplus s_2)}
\end{equation*}
The proof is identical to that for Property~\ref{property.domination_ours}\ref{property.domination_pps}, except that 
Lemma~\ref{lemma.SRCs_fracpart_cyclic} is applied by letting $a = \alpha_q(s_1), b = \alpha_q(s_2), c = \alpha_q(p)$ for each $q$.
\hfill\Halmos

We finally remark that the proof for Property~\ref{property.domination_ours}\ref{property.domination_ssa} was first shown in Proposition 5 of \cite{jepsen2008subset}. We show the hidden additive structure of this REF and domination criteria through Lemma~\ref{lemma.SRCs_fracpart_cyclic}, and make the natural generalization to sequences of visits (concatenation of REFs). 
\subsection{Proofs in Section 5.3}
\label{app.proofs_53}

We first show the following lemma: 
\begin{lemma}
    \label{lemma.ng_neighborhood_nested}
    Let $\calN^1$ and $\calN^2$ be two \textit{ng}-neighborhoods such that $N_i^1 \subseteq N_i^2$ for all $i \in \calV$. Then, $\calP(\calN^1) \supseteq \calP(\calN^2)$, and $\OPTLP(\calP(\calN^1)) \leq \OPTLP(\calP(\calN^2))$.
\end{lemma}
\begin{proof}{Proof of Lemma~\ref{lemma.ng_neighborhood_nested}.}
Suppose that $p \in \calP(\calN^2)$. Then $p$ is a feasible path, so $p \in \Pnone$. Additionally, let $N = (n_0, \dots, n_m)$ be the node sequence of $p$. Suppose that $j < k$ with $n_j = n_k$. Since $N$ is \textit{ng}-feasible with respect to $\calN^2$, there exists a $\ell$ with $j < \ell < k$ such that $n_j \notin N^2_{n_\ell}$. Since $N^1_{n_\ell} \subseteq N^2_{n_\ell}$, $n_j \notin N^1_{n_\ell}$. This shows that $N$ is \textit{ng}-feasible with respect to $\calN^1$, and that $p$ is \textit{ng}-feasible with respect to $\calN^1$. Therefore, $\calP(\calN^2) \subseteq \calP(\calN^1)$, and $\OPTLP(\calP(\calN^1)) \leq \OPTLP(\calP(\calN^2))$.
\hfill\Halmos
\end{proof}

\subsection*{Proof of Corollary~\ref{cor.iterative_cg_correctness}.}

We first show that Steps 1--4 of Algorithm~\ref{alg.AdaptiveColumnGenerationWithCuts} returns an optimal solution to $\EVRPLP(\Pelem)$ in a finite number of iterations. Note that in Algorithm~\ref{alg.AdaptiveColumnGenerationWithCuts}, the \textit{ng}-neighborhoods $\calN^t$ used across iterations are nested: for all $t$, $\calN^t$ and $\calN^{t+1}$ satisfy $N_i^t \subseteq N_i^{t+1}$ for all $i \in \calV$, and the inclusion is strict for at least one $i$. 
Per Lemma~\ref{lemma.ng_neighborhood_nested}, $\calP(\calN^t) \supseteq \calP(\calN^{t+1})$ and $\OPTLP(\calP(\calN^t)) \leq \OPTLP(\calP(\calN^{t+1}))$. Next, consider a non-elementary path $p$ in the support of the incumbent solution $\bm{z}^t$. That path admits a cycle $\{i, n_0, \dots, n_m, i\}$ in $N(p)$, with $i \in \Custs$. Then, the addition of $i$ to $N_{n_0}, \dots, N_{n_m}$ results in $p$ no longer being \textit{ng}-feasible for $\calN^{t+1}$ and hence for any subsequent \textit{ng}-neighborhood. Therefore, the quantity $\sum_{i \in \calV} |N_i^t|$ takes integer values, strictly increases whenever Step 4 is reached, and is upper-bounded by $|\calV|^2$. This proves that there exists some iteration $t_1$ at which all paths in the support of $\bm{z}^{t_1}$ are elementary, so that $\bm{z}^{t_1}$ is a feasible solution to $\EVRPLP(\Pelem)$ with optimal value $\OPTLP(\Pelem)$.

Next, let $t_1, t_2, \dots$ indicate the iterations in which Step 5 is reached. Since each cut separates $\bm{z}^{t_k}$ from the feasible set of the relaxation, the sequence of cuts defines a sequence of nested relaxations with objective values $\OPTLP(\Pelem) = \OPT^{t_1} \leq \OPT^{t_2} \leq \dots \leq \OPT(\Pelem)$. Furthermore, the family of subset-row cuts such that $|S| = 3$ and $w_i = \frac{1}{2} \forall \ i \in S$ is finite. Thus, Algorithm~\ref{alg.AdaptiveColumnGenerationWithCuts} terminates in a finite number of iterations and its optimum $\OPT$ satisfies $\OPTLP(\Pelem) \leq \OPT \leq \OPT(\Pelem)$. 
\hfill \Halmos
\section{Computational results}
\label{app.computational_results}

\subsection{Supplementary results on new multi-depot instances}
\label{app.multidepot}

Tables~\ref{tab.multidepot_small}--\ref{tab.multidepot_big} present results on the full column generation algorithm (Steps 1--5 of Algorithm~\ref{alg.AdaptiveColumnGenerationWithCuts}) on our new multi-depot instances, varying the number of customers $N = d x y$ along the customer density $d$ and the grid size $(x, y)$. Each row is the geometric mean of the 5 random instances generated with those parameters. 
For the instances where both our method and the path-based label-setting benchmark terminate within the 1-hour time limit, our method is 60\%--80\% faster and achieves near-zero optimality gaps. 
For the more challenging instances where our method and the path-based benchmark do not terminate within the 1-hour time limit, our method achieves significantly better primal bounds and dual bounds compared to the benchmark; as Figure~\ref{fig.multidepot_bounds} demonstrates, this is due to more column generation iterations enabled by the faster computational times. 
Our method is more advantageous with wider time windows (comparing Table~\ref{tab.multidepot_big} to Table~\ref{tab.multidepot_small}, larger customer density $d$, and larger service area $(x, y)$.

\begin{table}
    \centering
    \footnotesize
    \caption{Geometric means of the computational time, the dual and primal bounds, and the integrality gap. for multi-depot instances with narrower time windows. The values in parentheses indicate the percentage change of 2-LS with respect to LS (or the percentage point change for ``Gap'').
    }
    \label{tab.multidepot_small}
    \centering
    \begin{tabular}{
        S[table-format=3.0] *{3}{S[table-format=1.0]}
        *{4}{
            r
            r
            r
            r
        }
    }
        \toprule
        & & & 
        & \multicolumn{4}{c}{LS}
        & \multicolumn{4}{c}{2-LS}
        \\
        \cmidrule(lr){5-8}
        \cmidrule(lr){9-12}
        {$N$} & {$x$} & {$y$} & {$d$}  
        & {$t$ (s)} & {(D)} & {(P)} & {Gap}
        & {$t$ (s)} & {(D)} & {(P)} & {Gap}
        \\
        \midrule
        16 & 2 & 4 & 2 & 25.59 & 131063 & 131063 & 0.00\% & 6.75 & 131063 & 131063 & 0.00\% \\ 
        & & & & & & & & ($-$73.62\%) & (0.00\%) & (0.00\%) & (0.00\%) \\ 
        24 & 3 & 4 & 2 & 54.34 & 192509 & 192509 & 0.00\% & 13.20 & 192509 & 192509 & 0.00\% \\ 
        & & & & & & & & ($-$75.71\%) & (0.00\%) & (0.00\%) & (0.00\%) \\ 
        32 & 4 & 4 & 2 & 232.92 & 265381 & 265404 & 0.01\% & 39.58 & 265381 & 265966 & 0.22\% \\ 
        & & & & & & & & ($-$83.01\%) & (0.00\%) & (0.21\%) & (0.21\%) \\ 
        40 & 5 & 4 & 2 & 503.19 & 319728 & 319728 & 0.00\% & 93.45 & 319728 & 319728 & 0.00\% \\ 
        & & & & & & & & ($-$81.43\%) & (0.00\%) & (0.00\%) & (0.00\%) \\
        \midrule 
        24 & 2 & 4 & 3 & 88.37 & 159117 & 159117 & 0.00\% & 17.75 & 159117 & 159117 & 0.00\% \\ 
        & & & & & & & & ($-$79.91\%) & (0.00\%) & (0.00\%) & (0.00\%) \\ 
        36 & 3 & 4 & 3 & 363.90 & 239308 & 239308 & 0.00\% & 95.76 & 239308 & 239308 & 0.00\% \\ 
        & & & & & & & & ($-$73.69\%) & (0.00\%) & (0.00\%) & (0.00\%) \\ 
        48 & 4 & 4 & 3 & 741.40 & 326128 & 326758 & 0.19\% & 232.78 & 325958 & 327297 & 0.40\% \\ 
        & & & & & & & & ($-$68.60\%) & ($-$0.05\%) & (0.16\%) & (0.21\%) \\ 
        60 & 5 & 4 & 3 & 1270.03 & 393665 & 393756 & 0.02\% & 269.47 & 393756 & 393756 & 0.00\% \\ 
        & & & & & & & & ($-$78.78\%) & (0.02\%) & (0.00\%) & ($-$0.02\%) \\
        \midrule 
        32 & 2 & 4 & 4 & 544.97 & 188054 & 188054 & 0.00\% & 115.89 & 188054 & 188054 & 0.00\% \\ 
        & & & & & & & & ($-$78.74\%) & (0.00\%) & (0.00\%) & (0.00\%) \\ 
        48 & 3 & 4 & 4 & 1394.32 & 273942 & 275276 & 0.48\% & 379.86 & 273942 & 275276 & 0.48\% \\ 
        & & & & & & & & ($-$72.76\%) & (0.00\%) & (0.00\%) & (0.00\%) \\ 
        64 & 4 & 4 & 4 & 2462.83 & 363931 & 368628 & 1.23\% & 713.40 & 364568 & 364568 & 0.00\% \\ 
        & & & & & & & & ($-$71.03\%) & (0.17\%) & ($-$1.10\%) & ($-$1.23\%) \\ 
        80 & 5 & 4 & 4 & 3189.44 & 469417 & 481779 & 2.49\% & 1654.25 & 470614 & 480031 & 1.91\% \\ 
        & & & & & & & & ($-$48.13\%) & (0.26\%) & ($-$0.36\%) & ($-$0.58\%) \\
        \midrule
        40 & 2 & 4 & 5 & 1458.52 & 200791 & 200929 & 0.07\% & 494.58 & 200929 & 200929 & 0.00\% \\ 
        & & & & & & & & ($-$66.09\%) & (0.07\%) & (0.00\%) & ($-$0.07\%) \\ 
        60 & 3 & 4 & 5 & 3264.05 & 308430 & 322286 & 4.27\% & 2619.04 & 310187 & 320016 & 3.03\% \\ 
        & & & & & & & & ($-$19.76\%) & (0.57\%) & ($-$0.70\%) & ($-$1.24\%) \\ 
        80 & 4 & 4 & 5 & 3037.97 & 409314 & 435965 & 6.07\% & 2754.91 & 411808 & 429839 & 4.14\% \\ 
        & & & & & & & & ($-$9.32\%) & (0.61\%) & ($-$1.41\%) & ($-$1.93\%) \\ 
        100 & 5 & 4 & 5 & 2865.91 & 524102 & 612108 & 13.97\% & 2938.66 & 540740 & 555879 & 2.65\% \\ 
        & & & & & & & & (2.54\%) & (3.17\%) & ($-$9.19\%) & ($-$11.31\%) \\ 
        \bottomrule
    \end{tabular}
\end{table}

\begin{table}
    \centering
    \footnotesize
    \caption{Geometric means of the computational time, the dual and primal bounds, and the integrality gap. for multi-depot instances with wider time windows. The values in parentheses indicate the percentage change of 2-LS with respect to LS (or the percentage point change for ``Gap''). 
    }
    \label{tab.multidepot_big}
    \centering
    \begin{tabular}{
        S[table-format=3.0] *{3}{S[table-format=1.0]}
        *{4}{
            r
            r
            r
            r
        }
    }
        \toprule
        & & & 
        & \multicolumn{4}{c}{LS}
        & \multicolumn{4}{c}{2-LS}
        \\
        \cmidrule(lr){5-8}
        \cmidrule(lr){9-12}
        {$N$} & {$x$} & {$y$} & {$d$}  
        & {$t$ (s)} & {(D)} & {(P)} & {Gap}
        & {$t$ (s)} & {(D)} & {(P)} & {Gap}
        \\
        \midrule
        16 & 2 & 4 & 2 & 44.13 & 124033 & 124542 & 0.40\% & 7.76 & 124033 & 124542 & 0.40\% \\ 
        & & & & & & & & ($-$82.41\%) & (0.00\%) & (0.00\%) & (0.00\%) \\ 
        24 & 3 & 4 & 2 & 204.37 & 182950 & 182950 & 0.00\% & 37.37 & 182950 & 182950 & 0.00\% \\ 
        & & & & & & & & ($-$81.72\%) & (0.00\%) & (0.00\%) & (0.00\%) \\ 
        32 & 4 & 4 & 2 & 486.02 & 240030 & 240030 & 0.00\% & 64.22 & 240030 & 240030 & 0.00\% \\ 
        & & & & & & & & ($-$86.79\%) & (0.00\%) & (0.00\%) & (0.00\%) \\ 
        40 & 5 & 4 & 2 & 769.64 & 287360 & 287495 & 0.05\% & 124.34 & 287360 & 287495 & 0.05\% \\ 
        & & & & & & & & ($-$83.85\%) & (0.00\%) & (0.00\%) & (0.00\%) \\ 
        \midrule
        24 & 2 & 4 & 3 & 338.45 & 147283 & 147724 & 0.30\% & 66.51 & 147498 & 147692 & 0.13\% \\ 
        & & & & & & & & ($-$80.35\%) & (0.15\%) & ($-$0.02\%) & ($-$0.17\%) \\ 
        36 & 3 & 4 & 3 & 1198.88 & 216347 & 218653 & 1.01\% & 334.54 & 216805 & 216828 & 0.01\% \\ 
        & & & & & & & & ($-$72.10\%) & (0.21\%) & ($-$0.83\%) & ($-$1.00\%) \\ 
        48 & 4 & 4 & 3 & 2113.85 & 291105 & 297193 & 2.00\% & 861.45 & 292262 & 292805 & 0.18\% \\ 
        & & & & & & & & ($-$59.25\%) & (0.40\%) & ($-$1.48\%) & ($-$1.82\%) \\ 
        60 & 5 & 4 & 3 & 3308.28 & 348739 & 363604 & 3.96\% & 1261.50 & 349875 & 351954 & 0.58\% \\ 
        & & & & & & & & ($-$61.87\%) & (0.33\%) & ($-$3.20\%) & ($-$3.38\%) \\ 
        \midrule
        32 & 2 & 4 & 4 & 1719.34 & 162740 & 166302 & 2.11\% & 935.07 & 163838 & 164773 & 0.56\% \\ 
        & & & & & & & & ($-$45.61\%) & (0.67\%) & ($-$0.92\%) & ($-$1.55\%) \\ 
        48 & 3 & 4 & 4 & 2953.17 & 243147 & 253454 & 3.99\% & 1999.61 & 245006 & 247778 & 1.10\% \\ 
        & & & & & & & & ($-$32.29\%) & (0.76\%) & ($-$2.24\%) & ($-$2.90\%) \\ 
        64 & 4 & 4 & 4 & 3335.99 & 310804 & 356480 & 11.57\% & 2565.19 & 318837 & 328307 & 2.82\% \\ 
        & & & & & & & & ($-$23.11\%) & (2.58\%) & ($-$7.90\%) & ($-$8.75\%) \\ 
        80 & 5 & 4 & 4 & 2926.38 & 374185 & 526957 & 28.51\% & 3738.04 & 405763 & 450957 & 9.84\% \\ 
        & & & & & & & & (27.74\%) & (8.44\%) & ($-$14.42\%) & ($-$18.67\%) \\ 
        \midrule
        40 & 2 & 4 & 5 & 2619.24 & 174948 & 187018 & 6.29\% & 2097.62 & 176845 & 182850 & 3.18\% \\ 
        & & & & & & & & ($-$19.92\%) & (1.08\%) & ($-$2.23\%) & ($-$3.11\%) \\ 
        60 & 3 & 4 & 5 & 2302.65 & 250593 & 312472 & 18.33\% & 2340.47 & 259769 & 281828 & 7.76\% \\ 
        & & & & & & & & (1.64\%) & (3.66\%) & ($-$9.81\%) & ($-$10.57\%) \\ 
        80 & 4 & 4 & 5 & 3145.27 & 312714 & 472223 & 33.68\% & 3123.98 & 353265 & 385035 & 8.18\% \\ 
        & & & & & & & & ($-$0.68\%) & (12.97\%) & ($-$18.46\%) & ($-$25.49\%) \\ 
        100 & 5 & 4 & 5 & {--} & {--} & {--} & {--} & 3397.09 & 444600 & 498598 & 10.74\% \\ 
        & & & & & & & & {--} & {--} & {--} & {--} \\ 
        \bottomrule
    \end{tabular}
    \centering
    \begin{tablenotes}
        \item A ``--'' denotes parameter settings where no results were obtained within the 1-hour time limit.
    \end{tablenotes}
\end{table}

}

\end{document}